\documentclass[10pt, reqno]{amsart}

\RequirePackage{amsthm,amsmath,amsfonts,amssymb}
\RequirePackage[numbers,sort&compress]{natbib}
\RequirePackage[colorlinks,citecolor=blue,urlcolor=blue]{hyperref}
\RequirePackage{graphicx}
\usepackage{tikz}

\theoremstyle{plain}

\newtheorem{theorem}{Theorem}[section]
\newtheorem{proposition}{Proposition}[section]
\newtheorem{corollary}{Corollary}[section]
\newtheorem{lemma}[theorem]{Lemma}
\theoremstyle{remark}
\newtheorem{definition}[theorem]{Definition}

\newcommand{\ssk}{\smallskip}

\renewcommand{\epsilon}{\varepsilon}

\newcommand\bbE{\mathbb{E}}

\newcommand\bbN{\mathbb{N}}
\newcommand\bbR{\mathbb{R}}

\newcommand{\mcB}{\mathcal{B}}

\newcommand{\mcG}{\mathcal{G}}

\newcommand{\mcI}{\mathcal{I}}

\newcommand{\mcL}{\mathcal{L}}

\newcommand{\mcT}{\mathcal T}

\makeatletter
\newcommand*{\defeq}{\mathrel{\rlap{%
                     \raisebox{0.3ex}{$\m@th\cdot$}}%
                     \raisebox{-0.3ex}{$\m@th\cdot$}}%
                     =}
\makeatother

\makeatletter
\newcommand*{\eqdef}{=
 		    \mathrel{\rlap{%
                     \raisebox{0.3ex}{$\m@th\cdot$}}%
                     \raisebox{-0.3ex}{$\m@th\cdot$}}%
                     }
\makeatother

\def\Labe{\mathfrak{e}}

\def\Labn{\mathfrak{n}}

\def\Labhom{\mathfrak{t}}

\def\CT{T}

\def\CI{\mathcal{I}}
\def\N{\mathbb{N}}
\def\one{\mathbf{1}}
\def\s{\frak{s}}
\def\root{\mathrm{root}}
\def\nonroot{\mathrm{non-root}}

\def\id{\mathrm{id}}

\usetikzlibrary{calc}
\usetikzlibrary{decorations}
\usetikzlibrary{positioning}
\usetikzlibrary{shapes}
\usetikzlibrary{external}

\definecolor{blush}{rgb}{0.87, 0.36, 0.51}
\definecolor{brightcerulean}{rgb}{0.11, 0.67, 0.84}
\definecolor{greenryb}{rgb}{0.4, 0.69, 0.2}

\newif\ifdark
\darkfalse

\ifdark
\definecolor{darkred}{rgb}{0.9,0.2,0.2}
\definecolor{darkblue}{rgb}{0.7,0.3,1}
\definecolor{darkgreen}{rgb}{0.1,0.9,0.1}
\definecolor{franck}{rgb}{0,0.8,1}
\definecolor{pagebackground}{rgb}{.15,.21,.18}
\definecolor{pageforeground}{rgb}{.84,.84,.85}
\pagecolor{pagebackground}
\AtBeginDocument{\globalcolor{pageforeground}}
\definecolor{symbols}{rgb}{0,0.7,1}
\colorlet{connection}{red!80!black}
\colorlet{boxcolor}{blue!50}

\else

\definecolor{darkred}{rgb}{0.7,0.1,0.1}
\definecolor{darkblue}{rgb}{0.4,0.1,0.8}
\definecolor{darkgreen}{rgb}{0.1,0.7,0.1}
\definecolor{franck}{rgb}{0,0,1}
\definecolor{pagebackground}{rgb}{1,1,1}
\definecolor{pageforeground}{rgb}{0,0,0}
\colorlet{symbols}{blue!90!black}
\colorlet{connection}{red!30!black}
\colorlet{boxcolor}{blue!50!black}

\fi

\def\slash{\leavevmode\unskip\kern0.18em/\penalty\exhyphenpenalty\kern0.18em}
\def\dash{\leavevmode\unskip\kern0.18em--\penalty\exhyphenpenalty\kern0.18em}

\DeclareMathAlphabet{\mathbbm}{U}{bbm}{m}{n}

\DeclareFontFamily{U}{BOONDOX-calo}{\skewchar\font=45 }
\DeclareFontShape{U}{BOONDOX-calo}{m}{n}{
	<-> s*[1.05] BOONDOX-r-calo}{}
\DeclareFontShape{U}{BOONDOX-calo}{b}{n}{
	<-> s*[1.05] BOONDOX-b-calo}{}
\DeclareMathAlphabet{\mcb}{U}{BOONDOX-calo}{m}{n}
\SetMathAlphabet{\mcb}{bold}{U}{BOONDOX-calo}{b}{n}

\DeclareMathAlphabet{\mathbbm}{U}{bbm}{m}{n}

\DeclareMathAlphabet{\mcb}{U}{BOONDOX-calo}{m}{n}
\SetMathAlphabet{\mcb}{bold}{U}{BOONDOX-calo}{b}{n}
\DeclareFontFamily{U}{mathx}{\hyphenchar\font45}
\DeclareFontShape{U}{mathx}{m}{n}{
	<5> <6> <7> <8> <9> <10>
	<10.95> <12> <14.4> <17.28> <20.74> <24.88>
	mathx10
}{}
\DeclareSymbolFont{mathx}{U}{mathx}{m}{n}
\DeclareMathSymbol{\bigtimes}{1}{mathx}{"91}

\def\s{\mathfrak{s}}

\providecommand{\figures}{false}
{ \ifthenelse{\equal{\figures}{false}} {#1}{\[ {\rm Figure \ missing !} \]} }{}
\def\id{\mathrm{id}}

\def\CT{\mathcal{T}}

\tikzstyle{tinydots}=[dash pattern=on \pgflinewidth off \pgflinewidth]
\tikzstyle{superdense}=[dash pattern=on 4pt off 1pt]

\def\root{\mathrm{root}}

\def\nonroot{\mathrm{non-root}}

\newcommand{\beq}{\begin{equation}}
	\newcommand{\eeq}{\end{equation}}

\usepackage{empheq}

\def\Labe{\mathfrak{e}}

\def\Labn{\mathfrak{n}}

\def\Labhom{\mathfrak{t}}

\def\${|\!|\!|}

\def\root{\mathrm{root}}

\def\nonroot{\mathrm{non-root}}

\newenvironment{DIFnomarkup}{}{} 

\newcommand{\rrightarrow}{{\to\hskip -4.9mm\raise 1pt\hbox{$\to$}}}

\newfont{\indic}{bbmss12}

\def\Nabla_#1{\nabla_{\!#1}}

\makeatletter
\pgfdeclareshape{crosscircle}
{
	\inheritsavedanchors[from=circle] 
	\inheritanchorborder[from=circle]
	\inheritanchor[from=circle]{north}
	\inheritanchor[from=circle]{north west}
	\inheritanchor[from=circle]{north east}
	\inheritanchor[from=circle]{center}
	\inheritanchor[from=circle]{west}
	\inheritanchor[from=circle]{east}
	\inheritanchor[from=circle]{mid}
	\inheritanchor[from=circle]{mid west}
	\inheritanchor[from=circle]{mid east}
	\inheritanchor[from=circle]{base}
	\inheritanchor[from=circle]{base west}
	\inheritanchor[from=circle]{base east}
	\inheritanchor[from=circle]{south}
	\inheritanchor[from=circle]{south west}
	\inheritanchor[from=circle]{south east}
	\inheritbackgroundpath[from=circle]
	\foregroundpath{
		\centerpoint%
		\pgf@xc=\pgf@x%
		\pgf@yc=\pgf@y%
		\pgfutil@tempdima=\radius%
		\pgfmathsetlength{\pgf@xb}{\pgfkeysvalueof{/pgf/outer xsep}}%
		\pgfmathsetlength{\pgf@yb}{\pgfkeysvalueof{/pgf/outer ysep}}%
		\ifdim\pgf@xb<\pgf@yb%
		\advance\pgfutil@tempdima by-\pgf@yb%
		\else%
		\advance\pgfutil@tempdima by-\pgf@xb%
		\fi%
		\pgfpathmoveto{\pgfpointadd{\pgfqpoint{\pgf@xc}{\pgf@yc}}{\pgfqpoint{-0.707107\pgfutil@tempdima}{-0.707107\pgfutil@tempdima}}}
		\pgfpathlineto{\pgfpointadd{\pgfqpoint{\pgf@xc}{\pgf@yc}}{\pgfqpoint{0.707107\pgfutil@tempdima}{0.707107\pgfutil@tempdima}}}
		\pgfpathmoveto{\pgfpointadd{\pgfqpoint{\pgf@xc}{\pgf@yc}}{\pgfqpoint{-0.707107\pgfutil@tempdima}{0.707107\pgfutil@tempdima}}}
		\pgfpathlineto{\pgfpointadd{\pgfqpoint{\pgf@xc}{\pgf@yc}}{\pgfqpoint{0.707107\pgfutil@tempdima}{-0.707107\pgfutil@tempdima}}}
	}
}
\makeatother

\def\symbol#1{\textcolor{symbols}{#1}}

\def\decorate#1#2{
	\ifnum#2>0
	\foreach \count in {1,...,#2}{
		let
		\p1 = (sourcenode.center),
		\p2 = (sourcenode.east),
		\n1 = {\x2-\x1},
		\n2 = {1mm},
		\n3 = {(1.3+0.6*(\count-1))*\n1},
		\n4 = {0.7*\n1}
		in 
		node[rectangle,fill=symbols,rotate=30,inner sep=0pt,minimum width=0.2*\n2,minimum height=\n2] at ($(sourcenode.center) + (\n3,\n4)$) {}
	}
	\fi
	\ifnum#1>0
	\foreach \count in {1,...,#1}{
		let
		\p1 = (sourcenode.center),
		\p2 = (sourcenode.east),
		\n1 = {\x2-\x1},
		\n2 = {1mm},
		\n3 = {(1.3+0.6*(\count-1))*\n1},
		\n4 = {0.7*\n1}
		in 
		node[rectangle,fill=symbols,rotate=-30,inner sep=0pt,minimum width=0.2*\n2,minimum height=\n2] at ($(sourcenode.center) + (-\n3,\n4)$) {}
	}
	\fi
}

\tikzset{
	dectriangle/.style 2 args={
		triangle,
		alias=sourcenode,
		append after command={\decorate{#1}{#2}}
	},
	dectriangle/.default={0}{0},
}

\tikzset{
	cross/.style={path picture={ 
			\draw[symbols]
			(path picture bounding box.south east) -- (path picture bounding box.north west) (path picture bounding box.south west) -- (path picture bounding box.north east);
	}},
	root/.style={circle,fill=green!50!black,inner sep=0pt, minimum size=1.2mm},
	dot/.style={circle,fill=pageforeground,inner sep=0pt, minimum size=1mm},
	dotred/.style={circle,fill=pageforeground!50!pagebackground,inner sep=0pt, minimum size=2mm},
	var/.style={circle,fill=pageforeground!10!pagebackground,draw=pageforeground,inner sep=0pt, minimum size=3mm},
	kernel/.style={semithick,shorten >=2pt,shorten <=2pt},
	kernels/.style={snake=zigzag,shorten >=2pt,shorten <=2pt,segment amplitude=1pt,segment length=4pt,line before snake=2pt,line after snake=5pt,},
	rho/.style={densely dashed,semithick,shorten >=2pt,shorten <=2pt},
	testfcn/.style={dotted,semithick,shorten >=2pt,shorten <=2pt},
	renorm/.style={shape=circle,fill=pagebackground,inner sep=1pt},
	labl/.style={shape=rectangle,fill=pagebackground,inner sep=1pt},
	xic/.style={very thin,circle,draw=symbols,fill=symbols,inner sep=0pt,minimum size=1.2mm},
	g/.style={very thin,rectangle,draw=symbols,fill=symbols!10!pagebackground,inner sep=0pt,minimum width=2.5mm,minimum height=1.2mm},
	xi/.style={very thin,circle,draw=symbols,fill=symbols!10!pagebackground,inner sep=0pt,minimum size=1.2mm},
	xies/.style={very thin,rectangle,fill=green!50!black!25,draw=symbols,inner sep=0pt,minimum size=1.1mm},
	xiesf/.style={very thin,rectangle,fill=green!50!black,draw=symbols,inner sep=0pt,minimum size=1.1mm},
	xix/.style={very thin,crosscircle,fill=symbols!10!pagebackground,draw=symbols,inner sep=0pt,minimum size=1.2mm},
	X/.style={very thin,cross,rectangle,fill=pagebackground,draw=symbols,inner sep=0pt,minimum size=1.2mm},
	xib/.style={thin,circle,fill=symbols!10!pagebackground,draw=symbols,inner sep=0pt,minimum size=1.6mm},
	xie/.style={thin,circle,fill=green!50!black,draw=symbols,inner sep=0pt,minimum size=1.6mm},
	xid/.style={thin,circle,fill=symbols,draw=symbols,inner sep=0pt,minimum size=1.6mm},
	xibx/.style={thin,crosscircle,fill=symbols!10!pagebackground,draw=symbols,inner sep=0pt,minimum size=1.6mm},
	kernels2/.style={very thick,draw=connection,segment length=12pt},
	keps/.style={thin,draw=symbols,->},
	kepspr/.style={thick,draw=connection,->},
	krho/.style={thin,draw=symbols,superdense,->},
	krhopr/.style={thick,draw=connection,superdense},
	triangle/.style = { regular polygon, regular polygon sides=3},
	not/.style={thin,circle,draw=connection,fill=connection,inner sep=0pt,minimum size=0.5mm},
	diff/.style = {very thin,draw=symbols,triangle,fill=red!50!black,inner sep=0pt,minimum size=1.6mm},
	diff1/.style = {very thin,dectriangle={1}{0},fill=red!50!black,draw=symbols,inner sep=0pt,minimum size=1.6mm},
	diff2/.style = {very thin,dectriangle={1}{1},fill=red!50!black,draw=symbols,inner sep=0pt,minimum size=1.6mm},
	diffmini/.style = {very thin,rectangle,fill=black,draw=black,inner sep=0pt,minimum size=0.75mm},
	kernelsmod/.style={very thick,draw=connection,segment length=12pt},
	rec/.style = {very thin,rectangle,fill=black,draw=black,inner sep=0pt,minimum size=2mm},
	cerc/.style={very thin,circle,draw=black,fill=symbols,inner sep=0pt,minimum size=2mm},
	stars/.style={very thin,star,star points=6,star point ratio=0.5, draw=black,fill=red,inner sep=0pt,minimum size=0.7mm},
	>=stealth,
}
\tikzset{
	root/.style={circle,fill=black!50,inner sep=0pt, minimum size=3mm},
	circ/.style={circle,fill=white,draw=black,very thin,inner sep=.5pt, minimum size=1.2mm},
	round1/.style={fill=white,outer sep = 0,inner sep=2pt,rounded corners=1mm,draw,text=black,thin,minimum size=1.2mm},
	circ1/.style={circle,fill=red!10,draw=red,very thin,inner sep=.5pt, minimum size=1.2mm},
	rect/.style={fill=white,outer sep = 0,inner sep=2pt,rectangle,draw,text=black,thin,minimum size=1.2mm},
	rect1/.style={fill=white,outer sep = 0,inner sep=2pt,rectangle,draw,text=black,thin,minimum size=1.2mm},
	round2/.style={fill=red!10,outer sep = 0,inner sep=2pt,rounded corners=1mm,draw,text=black,thin,minimum size=1.2mm},
	round3/.style={fill=blue!10,outer sep = 0,inner sep=2pt,rounded corners=1mm,draw,text=black,thin,minimum size=1.2mm}, 
	rect2/.style={fill=black!10,outer sep = 0,inner sep=2pt,rectangle,draw,text=black,thin,minimum size=1.2mm},
	dot/.style={circle,fill=black,inner sep=0pt, minimum size=1.2mm},
	dotred/.style={circle,fill=black!50,inner sep=0pt, minimum size=2mm},
	var/.style={circle,fill=black!10,draw=black,inner sep=0pt, minimum size=3mm},
	kernel/.style={semithick,shorten >=2pt,shorten <=2pt},
	diag/.style={thin,shorten >=4pt,shorten <=4pt},
	kernel1/.style={thick},
	kernels/.style={snake=zigzag,shorten >=2pt,shorten <=2pt,segment amplitude=1pt,segment length=4pt,line before snake=2pt,line after snake=5pt,},
	kernels1/.style={snake=zigzag,segment amplitude=0.5pt,segment length=2pt},
	rho1/.style={densely dotted,semithick},
	rho/.style={densely dashed,semithick,shorten >=2pt,shorten <=2pt},
	testfcn/.style={dotted,semithick,shorten >=2pt,shorten <=2pt},
	visible/.style={draw, circle, fill, inner sep=0.25ex},
	renorm/.style={shape=circle,fill=white,inner sep=1pt},
	labl/.style={shape=rectangle,fill=white,inner sep=1pt},
	xic/.style={very thin,circle,fill=symbols,draw=black,inner sep=0pt,minimum size=1.2mm},
	xi/.style={very thin,circle,fill=blue!10,draw=black,inner sep=0pt,minimum size=1.2mm},
	xib/.style={very thin,circle,fill=blue!10,draw=black,inner sep=0pt,minimum size=1.6mm},
	xie/.style={very thin,circle,fill=green!50!black,draw=black,inner sep=0pt,minimum size=1mm},
	xid/.style={very thin,circle,fill=symbols,draw=black,inner sep=0pt,minimum size=1.6mm},
	edgetype/.style={very thin,circle,draw=black,inner sep=0pt,minimum size=5mm},
	nodetype/.style={very thick,circle,draw=black,inner sep=0pt,minimum size=5mm},
	kernels2/.style={very thick,draw=connection,segment length=12pt},
	clean/.style={thin,circle,fill=black,inner sep=0pt,minimum size=1mm},	not/.style={thin,circle,fill=symbols,draw=connection,fill=connection,inner sep=0pt,minimum size=0.8mm},
	>=stealth,
}

\makeatletter
\def\DeclareSymbol#1#2#3{%
	\expandafter\gdef\csname MH@symb@#1\endcsname{\tikzsetnextfilename{symbol#1}%
		\tikz[baseline=#2,scale=0.15,draw=symbols,line join=round]{#3}}%
	\expandafter\gdef\csname MH@symb@#1s\endcsname{\scalebox{0.75}{\tikzsetnextfilename{symbol#1}%
			\tikz[baseline=#2,scale=0.15,draw=symbols,line join=round]{#3}}}%
	\expandafter\gdef\csname MH@symb@#1ss\endcsname{\scalebox{0.65}{\tikzsetnextfilename{symbol#1}%
			\tikz[baseline=#2,scale=0.15,draw=symbols,line join=round]{#3}}}%
}
\def\<#1>{\ifthenelse{\boolean{mmode}}{\mathchoice{\csname MH@symb@#1\endcsname}{\csname MH@symb@#1\endcsname}{\csname MH@symb@#1s\endcsname}{\csname MH@symb@#1ss\endcsname}}{\csname MH@symb@#1\endcsname}}
\makeatother

\DeclareSymbol{Xi22}{0.5}{\draw (0,0) node[xi] {} -- (-1,1) node[not] {} -- (0,2) node[xi] {};} 
\DeclareSymbol{Xi2}{-2}{\draw (-1,-0.25) node[xi] {} -- (0,1) node[xi] {};} 
\DeclareSymbol{Xi2b}{-2}{\draw (-1,-0.25) node[xic] {} -- (0,1) node[xic] {};} 
\DeclareSymbol{Xi2g}{-2}{\draw (-1,-0.25) node[xies] {} -- (0,1) node[xi] {};} 
\DeclareSymbol{Xi2g2}{-2}{\draw (-1,-0.25) node[xi] {} -- (0,1) node[xies] {};} 
\DeclareSymbol{cXi2}{-2}{\draw (0,-0.25) node[xi] {} -- (-1,1) node[xic] {};}
\DeclareSymbol{Xi3}{0}{\draw (0,0) node[xi] {} -- (-1,1) node[xi] {} -- (0,2) node[xi] {};}
\DeclareSymbol{XiIIXi}{0}{\draw (0,0) node[xi] {} -- (-1,1); \draw[kernels2] (-1,1) node[not] {} -- (0,2) node[xi] {};}

\DeclareSymbol{Xi4}{2}{\draw (0,0) node[xi] {} -- (-1,1) node[xi] {} -- (0,2) node[xi] {} -- (-1,3) node[xi] {};}
\DeclareSymbol{Xi4_1}{2}{\draw (0,0) node[xic] {} -- (-1,1) node[xic] {} -- (0,2) node[xi] {} -- (-1,3) node[xi] {};}
\DeclareSymbol{Xi4_2}{2}{\draw (0,0) node[xic] {} -- (-1,1) node[xi] {} -- (0,2) node[xi] {} -- (-1,3) node[xic] {};}
\DeclareSymbol{Xi2X}{-2}{\draw (0,-0.25) node[xi] {} -- (-1,1) node[xix] {};}
\DeclareSymbol{XXi2}{-2}{\draw (0,-0.25) node[xix] {} -- (-1,1) node[xi] {};}
\DeclareSymbol{IIXi}{0}{\draw (0,-0.25) node[not] {} -- (-1,1) node[xi] {} -- (0,2) node[xi] {};}
\DeclareSymbol{IXi^2}{-1}{\draw (-1,1) node[xi] {} -- (0,0) node[not] {} -- (1,1) node[xi] {};}
\DeclareSymbol{IIXi^2}{-4}{\draw (0,-1.5) node[not] {} -- (0,0);
	\draw[kernels2] (-1,1) node[xi] {} -- (0,0) node[not] {} -- (1,1) node[xi] {};}
\DeclareSymbol{XiX}{-2.8}{\node[xibx] {};}
\DeclareSymbol{tauX}{-2.8}{ \node[X] {};}
\DeclareSymbol{Xi}{-2.8}{\node[xib] {};}

\DeclareSymbol{IXiX}{-1}{\draw (0,-0.25) node[not] {} -- (-1,1) node[xix] {};}
\DeclareSymbol{IXi3}{2}{\draw (0,-0.25) node[not] {} -- (-1,1) node[xi] {} -- (0,2) node[xi] {} -- (-1,3) node[xi] {};}
\DeclareSymbol{IXi}{-2}{\draw (0,-0.25) node[not] {} -- (-1,1) node[xi] {};}
\DeclareSymbol{XiI}{-2}{\draw (0,-0.25) node[xi] {} -- (-1,1) node[not] {};}

\DeclareSymbol{Xi4b}{0}{\draw(0,1.5) node[xi] {} -- (0,0); \draw (-1,1) node[xi] {} -- (0,0) node[xi] {} -- (1,1) node[xi] {};}
\DeclareSymbol{Xi4b'}{0}{\draw(0,1.5) node[xi] {} -- (0,-0.2); \draw (-1,1) node[xi] {} -- (0,-0.2) node[not] {} -- (1,1) node[xi] {};}
\DeclareSymbol{Xi4c}{0}{\draw (0,1) -- (0.8,2.2) node[xi] {};\draw (0,-0.25) node[xi] {} -- (0,1) node[xi] {} -- (-0.8,2.2) node[xi] {};}
\DeclareSymbol{Xi4d}{-4.5}{\draw (0,-1.5) node[not] {} -- (0,0); \draw (-1,1) node[xi] {} -- (0,0) node[xi] {} -- (1,1) node[xi] {};}
\DeclareSymbol{Xi4e}{0}{\draw (0,2) node[xi] {} -- (-1,1) node[xi] {} -- (0,0) node[xi] {} -- (1,1) node[xi] {};}
\DeclareSymbol{Xi4e'}{0}{\draw (0,2) node[xi] {} -- (-1,1) node[xi] {} -- (0,-0.2) node[not] {} -- (1,1) node[xi] {};}

\DeclareSymbol{Xitwo}
{0}{\draw[kernels2] (0,0) node[not] {} -- (-1,1) node[not] {}
	-- (-2,2) node[not]{} -- (-3,3) node[xi]  {};
	\draw[kernels2] (0,0) -- (1,1) node[xi] {};
	\draw[kernels2] (-1,1) -- (0,2) node[xi] {};
	\draw[kernels2] (-2,2) -- (-1,3) node[xi] {};}

\DeclareSymbol{IXitwo}
{0}{\draw (-.7,1.2) node[xi] {} -- (0,-0.2) -- (.7,1.2) node[xi] {};}
\DeclareSymbol{I1Xitwo}
{0}{\draw[kernels2] (0,0) node[not] {} -- (-1,1) node[xi] {};
	\draw[kernels2] (0,0) -- (1,1) node[xi] {};}

\DeclareSymbol{I1Xitwobis}
{0}{\draw[kernels2] (0,0) node[not] {} -- (-1,1) node[xies] {};
	\draw[kernels2] (0,0) -- (1,1) node[xies] {};}

\DeclareSymbol{I1Xitwog}
{0}{\draw[kernels2] (0,0) node[not] {} -- (-1,1) node[xies] {};
	\draw[kernels2] (0,0) -- (1,1) node[xi] {};}

\DeclareSymbol{cI1Xitwo}
{0}{\draw[kernels2] (0,0) node[not] {} -- (-1,1) node[xic] {};
	\draw[kernels2] (0,0) -- (1,1) node[xi] {};}

\DeclareSymbol{I1IXi3}{0}{\draw (0,0) node[xi] {} -- (-1,1) ; 
	\draw[kernels2] (-1,1) node[not] {} -- (0,2) node[xi] {};
	\draw[kernels2] (-1,1) node[not] {} -- (-2,2) node[xi] {};}

\DeclareSymbol{I1Xi3c}{-1}{\draw[kernels2](0,1.5) node[xi] {} -- (0,0) node[not] {}; \draw (-1,1) node[xi] {} -- (0,0) ; \draw[kernels2] (0,0) -- (1,1) node[xi] {};}

\DeclareSymbol{I1Xi3cbis}{-1}{\draw[kernels2](0,1.5) node[xies] {} -- (0,0) node[not] {}; \draw (-1,1) node[xies] {} -- (0,0) ; \draw[kernels2] (0,0) -- (1,1) node[xies] {};}

\DeclareSymbol{I1IXi3b}{0}{\draw[kernels2] (0,0) node[not] {} -- (-1,1) ; \draw[kernels2] (0,0)   -- (1,1) node[xi] {} ;
	\draw (-1,1) node[xi] {} -- (0,2) node[xi] {};
}

\DeclareSymbol{I1IXi3c}{0}{\draw[kernels2] (0,0) node[not] {} -- (-1,1) ; \draw[kernels2] (0,0)   -- (1,1) node[xi] {} ;
	\draw[kernels2] (-1,1) node[not] {} -- (0,2) node[xi] {};
	\draw[kernels2] (-1,1) node[not] {} -- (-2,2) node[xi] {};}

\DeclareSymbol{I1IXi3cbis}{0}{\draw[kernels2] (0,0) node[not] {} -- (-1,1) ; \draw[kernels2] (0,0)   -- (1,1) node[xies] {} ;
	\draw[kernels2] (-1,1) node[not] {} -- (0,2) node[xies] {};
	\draw[kernels2] (-1,1) node[not] {} -- (-2,2) node[xies] {};}

\DeclareSymbol{I1Xi}{0}{\draw[kernels2] (0,0) node[not] {} -- (-1,1)  node[xi] {} ;}

\DeclareSymbol{I1Xi4a}{2}{\draw[kernels2] (0,0) node[not] {} -- (-1,1) ; \draw[kernels2] (0,0) node[not] {} -- (1,1) node[xi] {} ;
	\draw (-1,1) node[xi] {} -- (0,2) node[xi] {} -- (-1,3) node[xi] {};}

\DeclareSymbol{cI1Xi4a}{2}{\draw[kernels2] (0,0) node[not] {} -- (-1,1) ; \draw[kernels2] (0,0) node[not] {} -- (1,1) node[xic] {} ;
	\draw (-1,1) node[xic] {} -- (0,2) node[xi] {} -- (-1,3) node[xi] {};}

\DeclareSymbol{I1Xi4b}{2}{\draw (0,0) node[xi] {} -- (-1,1) node[xi] {} -- (0,2) ; \draw[kernels2] (0,2) node[not] {} -- (-1,3) node[xi] {};\draw[kernels2] (0,2)  -- (1,3) node[xi] {};
}

\DeclareSymbol{cI1Xi4b}{2}{\draw (0,0) node[xic] {} -- (-1,1) node[xic] {} -- (0,2) ; \draw[kernels2] (0,2) node[not] {} -- (-1,3) node[xi] {};\draw[kernels2] (0,2)  -- (1,3) node[xi] {};
}

\DeclareSymbol{I1Xi4c}{2}{\draw (0,0) node[xi] {} -- (-1,1) node[not] {}; \draw[kernels2] (-1,1) -- (0,2) ; 
	\draw[kernels2] (-1,1) -- (-2,2) node[xi] {} ;
	\draw (0,2) node[xi] {} -- (-1,3) node[xi] {};}

\DeclareSymbol{cI1Xi4c}{2}{\draw (0,0) node[xic] {} -- (-1,1) node[not] {}; \draw[kernels2] (-1,1) -- (0,2) ; 
	\draw[kernels2] (-1,1) -- (-2,2) node[xic] {} ;
	\draw (0,2) node[xi] {} -- (-1,3) node[xi] {};}

\DeclareSymbol{I1Xi4ab}{2}{\draw[kernels2] (0,0) node[not] {} -- (-1,1) ; \draw[kernels2] (0,0) node[not] {} -- (1,1) node[xi] {};\draw (-1,1) node[xi] {} -- (0,2) ; \draw[kernels2] (0,2) node[not] {} -- (-1,3) node[xi] {};\draw[kernels2] (0,2)  -- (1,3) node[xi] {}; }

\DeclareSymbol{cI1Xi4ab}{2}{\draw[kernels2] (0,0) node[not] {} -- (-1,1) ; \draw[kernels2] (0,0) node[not] {} -- (1,1) node[xic] {};\draw (-1,1) node[xic] {} -- (0,2) ; \draw[kernels2] (0,2) node[not] {} -- (-1,3) node[xi] {};\draw[kernels2] (0,2)  -- (1,3) node[xi] {}; }

\DeclareSymbol{I1Xi4bc}{2}{\draw (0,0) node[xi] {} -- (-1,1) node[not] {}; \draw[kernels2] (-1,1) -- (0,2) ; 
	\draw[kernels2] (-1,1) -- (-2,2) node[xi] {} ; \draw[kernels2] (0,2) node[not] {} -- (-1,3) node[xi] {};\draw[kernels2] (0,2)  -- (1,3) node[xi] {};
}

\DeclareSymbol{cI1Xi4bc}{2}{\draw (0,0) node[xic] {} -- (-1,1) node[not] {}; \draw[kernels2] (-1,1) -- (0,2) ; 
	\draw[kernels2] (-1,1) -- (-2,2) node[xic] {} ; \draw[kernels2] (0,2) node[not] {} -- (-1,3) node[xi] {};\draw[kernels2] (0,2)  -- (1,3) node[xi] {};
}

\DeclareSymbol{I1Xi4abcc1}{2}{\draw[kernels2] (0,0) node[not] {} -- (-1,1) node[not] {}
	-- (-2,2) node[not]{} -- (-3,3) node[xic]  {};
	\draw[kernels2] (0,0) -- (1,1) node[xic] {};
	\draw[kernels2] (-1,1) -- (0,2) node[xi] {};
	\draw[kernels2] (-2,2) -- (-1,3) node[xi] {};
}

\DeclareSymbol{I1Xi4abcc1b}{2}{\draw[kernels2] (0,0) node[not] {} -- (-1,1) node[not] {}
	-- (-2,2) node[not]{} -- (-3,3) node[xi]  {};
	\draw[kernels2] (0,0) -- (1,1) node[xic] {};
	\draw[kernels2] (-1,1) -- (0,2) node[xic] {};
	\draw[kernels2] (-2,2) -- (-1,3) node[xi] {};
}

\DeclareSymbol{I1Xi4abcc2}{2}{\draw[kernels2] (0,0) node[not] {} -- (-1,1) node[not] {}
	-- (-2,2) node[not]{} -- (-3,3) node[xic]  {};
	\draw[kernels2] (0,0) -- (1,1) node[xi] {};
	\draw[kernels2] (-1,1) -- (0,2) node[xi] {};
	\draw[kernels2] (-2,2) -- (-1,3) node[xic] {};
}

\DeclareSymbol{I1Xi4ac}{2}{\draw[kernels2] (0,0) node[not] {} -- (-1,1) ; \draw[kernels2] (0,0) node[not] {} -- (1,1) node[xi] {}; 
	\draw[kernels2] (-1,1) node[not] {} -- (0,2) ; 
	\draw[kernels2] (-1,1) -- (-2,2) node[xi] {} ;
	\draw (0,2) node[xi] {} -- (-1,3) node[xi] {};}

\DeclareSymbol{cI1Xi4ac}{2}{\draw[kernels2] (0,0) node[not] {} -- (-1,1) ; \draw[kernels2] (0,0) node[not] {} -- (1,1) node[xic] {}; 
	\draw[kernels2] (-1,1) node[not] {} -- (0,2) ; 
	\draw[kernels2] (-1,1) -- (-2,2) node[xic] {} ;
	\draw (0,2) node[xi] {} -- (-1,3) node[xi] {};}

\DeclareSymbol{I1Xi4acc1}{2}{\draw[kernels2] (0,0) node[not] {} -- (-1,1) ; \draw[kernels2] (0,0) node[not] {} -- (1,1) node[xic] {}; 
	\draw[kernels2] (-1,1) node[not] {} -- (0,2) ; 
	\draw[kernels2] (-1,1) -- (-2,2) node[xi] {} ;
	\draw (0,2) node[xic] {} -- (-1,3) node[xi] {};}

\DeclareSymbol{I1Xi4acc2}{2}{\draw[kernels2] (0,0) node[not] {} -- (-1,1) ; \draw[kernels2] (0,0) node[not] {} -- (1,1) node[xic] {}; 
	\draw[kernels2] (-1,1) node[not] {} -- (0,2) ; 
	\draw[kernels2] (-1,1) -- (-2,2) node[xi] {} ;
	\draw (0,2) node[xi] {} -- (-1,3) node[xic] {};}

\DeclareSymbol{2I1Xi4}{2}{\draw[kernels2] (0,0) node[not] {} -- (-1,1) node[not] {};
	\draw[kernels2] (0,0) -- (1,1) node[not] {};
	\draw[kernels2] (-1,1) -- (-1.5,2.5) node[xi] {};
	\draw[kernels2] (-1,1) -- (-0.5,2.5) node[xi] {};
	\draw[kernels2] (1,1) -- (0.5,2.5) node[xi] {};
	\draw[kernels2] (1,1) -- (1.5,2.5) node[xi] {};
}

\DeclareSymbol{2I1Xi4dis}{2}{\draw[kernels2] (0,0) node[not] {} -- (-1,1) node[not] {};
	\draw[kernels2] (0,0) -- (1,1) node[not] {};
	\draw[kernels2] (-1,1) -- (-1.5,2.5) node[xies] {};
	\draw[kernels2] (-1,1) -- (-0.5,2.5) node[xies] {};
	\draw[kernels2] (1,1) -- (0.5,2.5) node[xies] {};
	\draw[kernels2] (1,1) -- (1.5,2.5) node[xies] {};
}

\DeclareSymbol{2I1Xi4c1}{2}{\draw[kernels2] (0,0) node[not] {} -- (-1,1) node[not] {};
	\draw[kernels2] (0,0) -- (1,1) node[not] {};
	\draw[kernels2] (-1,1) -- (-1.5,2.5) node[xic] {};
	\draw[kernels2] (-1,1) -- (-0.5,2.5) node[xi] {};
	\draw[kernels2] (1,1) -- (0.5,2.5) node[xic] {};
	\draw[kernels2] (1,1) -- (1.5,2.5) node[xi] {};
}

\DeclareSymbol{2I1Xi4c2}{2}{\draw[kernels2] (0,0) node[not] {} -- (-1,1) node[not] {};
	\draw[kernels2] (0,0) -- (1,1) node[not] {};
	\draw[kernels2] (-1,1) -- (-1.5,2.5) node[xic] {};
	\draw[kernels2] (-1,1) -- (-0.5,2.5) node[xic] {};
	\draw[kernels2] (1,1) -- (0.5,2.5) node[xi] {};
	\draw[kernels2] (1,1) -- (1.5,2.5) node[xi] {};
}

\DeclareSymbol{2I1Xi4b}{2}{\draw[kernels2] (0,0) node[not] {} -- (-1,1) ;
	\draw[kernels2] (0,0) -- (1,1);
	\draw (-1,1) node[xi] {} -- (-1,2.5) node[xi] {};
	\draw (1,1)  node[xi] {} -- (1,2.5) node[xi] {};
}

\DeclareSymbol{2I1Xi4bb}{2}{\draw[kernels2] (0,0) node[not] {} -- (-1,1) ;
	\draw[kernels2] (0,0) -- (1,1);
	\draw (-1,1) node[xi] {} -- (-1,2.5) node[xiesf] {};
	\draw (1,1)  node[xi] {} -- (1,2.5) node[xic] {};
}

\DeclareSymbol{2I1Xi4c}{2}{\draw[kernels2] (0,0) node[not] {} -- (-1,1);
	\draw[kernels2] (0,0) -- (1,1) node[not] {};
	\draw (-1,1)  node[xi] {} -- (-1,2.5) node[xi] {};
	\draw[kernels2] (1,1) -- (0.4,2.5) node[xi] {};
	\draw[kernels2] (1,1) -- (1.6,2.5) node[xi] {};
}

\DeclareSymbol{2I1Xi4cc1}{2}{\draw[kernels2] (0,0) node[not] {} -- (-1,1);
	\draw[kernels2] (0,0) -- (1,1) node[not] {};
	\draw (-1,1)  node[xic] {} -- (-1,2.5) node[xi] {};
	\draw[kernels2] (1,1) -- (0.4,2.5) node[xic] {};
	\draw[kernels2] (1,1) -- (1.6,2.5) node[xi] {};
}

\DeclareSymbol{2I1Xi4cc2}{2}{\draw[kernels2] (0,0) node[not] {} -- (-1,1);
	\draw[kernels2] (0,0) -- (1,1) node[not] {};
	\draw (-1,1)  node[xic] {} -- (-1,2.5) node[xic] {};
	\draw[kernels2] (1,1) -- (0.4,2.5) node[xi] {};
	\draw[kernels2] (1,1) -- (1.6,2.5) node[xi] {};
}

\DeclareSymbol{Xi4ba}{0}{\draw(-0.5,1.5) node[xi] {} -- (0,0); \draw (-1.5,1) node[xi] {} -- (0,0) node[not] {}; \draw[kernels2] (0,0) -- (1.5,1) node[xi] {};
	\draw[kernels2] (0,0) -- (0.5,1.5) node[xi] {} ;}

\DeclareSymbol{Xi4badis}{0}{\draw(-0.5,1.5) node[xies] {} -- (0,0); \draw (-1.5,1) node[xies] {} -- (0,0) node[not] {}; \draw[kernels2] (0,0) -- (1.5,1) node[xies] {};
	\draw[kernels2] (0,0) -- (0.5,1.5) node[xies] {} ;}

\DeclareSymbol{Xi4ba1}{0}{\draw(-0.5,1.5) node[xi] {} -- (0,0); \draw (-1.5,1) node[xi] {} -- (0,0) node[not] {}; \draw[kernels2] (0,0) -- (1.5,1) node[xic] {};
	\draw[kernels2] (0,0) -- (0.5,1.5) node[xic] {} ;}

\DeclareSymbol{Xi4ba1b}{0}{\draw(-0.5,1.5) node[xic] {} -- (0,0); \draw (-1.5,1) node[xic] {} -- (0,0) node[not] {}; \draw[kernels2] (0,0) -- (1.5,1) node[xi] {};
	\draw[kernels2] (0,0) -- (0.5,1.5) node[xi] {} ;}

\DeclareSymbol{Xi4ba1bdiff}{0}{\draw(-0.5,1.5) node[xic] {} -- (0,0); \draw (-1.5,1) node[xic] {} -- (0,0) node[not] {}; \draw (0,0) -- (1.5,1) node[xi] {};
	\draw (0,0) -- (0.5,1.5) node[xi] {};
	\draw(0,0) node[diff] {};}

\DeclareSymbol{Xi4ba1bb}{0}{\draw(-0.5,1.5) node[xic] {} -- (0,0); \draw (-1.5,1) node[xiesf] {} -- (0,0) node[not] {}; \draw[kernels2] (0,0) -- (1.5,1) node[xi] {};
	\draw[kernels2] (0,0) -- (0.5,1.5) node[xi] {} ;}

\DeclareSymbol{Xi4ba2}{0}{\draw(-0.5,1.5) node[xi] {} -- (0,0); \draw (-1.5,1) node[xic] {} -- (0,0) node[not] {}; \draw[kernels2] (0,0) -- (1.5,1) node[xi] {};
	\draw[kernels2] (0,0) -- (0.5,1.5) node[xic] {} ;}

\DeclareSymbol{Xi4ba2b}{0}{\draw(-0.5,1.5) node[xi] {} -- (0,0); \draw (-1.5,1) node[xic] {} -- (0,0) node[not] {}; \draw[kernels2] (0,0) -- (1.5,1) node[xi] {};
	\draw[kernels2] (0,0) -- (0.5,1.5) node[xiesf] {} ;}

\DeclareSymbol{Xi4ca}{0}{\draw (0,1) -- (-1,2.2) node[xi] {};\draw (0,-0.25) node[xi] {} -- (0,1) ; \draw[kernels2] (0,1) node[not] {} -- (1,2.2) node[xi] {};
	\draw[kernels2] (0,1) {} -- (0,2.7) node[xi] {};
}

\DeclareSymbol{Xi4cb}{0}{\draw (-1,1) -- (-2,2) node[xi] {};\draw[kernels2] (0,0)  -- (-1,1) node[xi] {} ; \draw[kernels2] (0,0) node[not] {} -- (1,1) node[xi] {} ; 
	\draw (-1,1) node[xi] {} -- (0,2) node[xi] {};}

\DeclareSymbol{Xi4cbb}{0}{\draw (-1,1) -- (-2,2) node[xiesf] {};\draw[kernels2] (0,0)  -- (-1,1) node[xi] {} ; \draw[kernels2] (0,0) node[not] {} -- (1,1) node[xi] {} ; 
	\draw (-1,1) node[xi] {} -- (0,2) node[xic] {};}

\DeclareSymbol{Xi4cbc1}{0}{\draw (-1,1) -- (-2,2) node[xic] {};\draw[kernels2] (0,0)  -- (-1,1) node[xic] {} ; \draw[kernels2] (0,0) node[not] {} -- (1,1) node[xi] {} ; 
	\draw (-1,1) node[xic] {} -- (0,2) node[xi] {};}

\DeclareSymbol{Xi4cbc2}{0}{\draw (-1,1) -- (-2,2) node[xi] {};\draw[kernels2] (0,0)  -- (-1,1) node[xi] {} ; \draw[kernels2] (0,0) node[not] {} -- (1,1) node[xic] {} ; 
	\draw (-1,1) node[xic] {} -- (0,2) node[xi] {};}

\DeclareSymbol{Xi4cab}{0}{\draw (-1,1) -- (-2,2) node[xi] {};\draw[kernels2] (0,0)  -- (-1,1); \draw[kernels2] (0,0) node[not] {} -- (1,1) node[xi] {} ; 
	\draw[kernels2] (-1,1)  {} -- (0,2) node[xi] {};
	\draw[kernels2] (-1,1) node[not] {} -- (-1,2.5) node[xi] {};
}

\DeclareSymbol{Xi4cabdis}{0}{\draw (-1,1) -- (-2,2) node[xies] {};\draw[kernels2] (0,0)  -- (-1,1); \draw[kernels2] (0,0) node[not] {} -- (1,1) node[xies] {} ; 
	\draw[kernels2] (-1,1)  {} -- (0,2) node[xies] {};
	\draw[kernels2] (-1,1) node[not] {} -- (-1,2.5) node[xies] {};
}

\DeclareSymbol{Xi4cabc1}{0}{\draw (-1,1) -- (-2,2) node[xi] {};\draw[kernels2] (0,0)  -- (-1,1); \draw[kernels2] (0,0) node[not] {} -- (1,1) node[xic] {} ; 
	\draw[kernels2] (-1,1)  {} -- (0,2) node[xic] {};
	\draw[kernels2] (-1,1) node[not] {} -- (-1,2.5) node[xi] {};
}

\DeclareSymbol{Xi4cabc2}{0}{\draw (-1,1) -- (-2,2) node[xic] {};\draw[kernels2] (0,0)  -- (-1,1); \draw[kernels2] (0,0) node[not] {} -- (1,1) node[xic] {} ; 
	\draw[kernels2] (-1,1)  {} -- (0,2) node[xi] {};
	\draw[kernels2] (-1,1) node[not] {} -- (-1,2.5) node[xi] {};
}

\DeclareSymbol{Xi4ea}{1.5}{\draw (-1,2.5) node[xi] {} -- (-1,1) node[xi] {} -- (0,0); 
	\draw[kernels2] (0,0)  -- (1,1) node[xi] {};
	\draw[kernels2] (0,0) node[not] {} -- (0,1.5) node[xi] {}; }

\DeclareSymbol{Xi4eac1}{1.5}{\draw (-1,2.5) node[xic] {} -- (-1,1) node[xi] {} -- (0,0); 
	\draw[kernels2] (0,0)  -- (1,1) node[xic] {};
	\draw[kernels2] (0,0) node[not] {} -- (0,1.5) node[xi] {}; }

\DeclareSymbol{Xi4eac1b}{1.5}{\draw (-1,2.5) node[xic] {} -- (-1,1) node[xi] {} -- (0,0); 
	\draw[kernels2] (0,0)  -- (1,1) node[xiesf] {};
	\draw[kernels2] (0,0) node[not] {} -- (0,1.5) node[xi] {}; }

\DeclareSymbol{Xi4eac2}{1.5}{\draw (-1,2.5) node[xic] {} -- (-1,1) node[xic] {} -- (0,0); 
	\draw[kernels2] (0,0)  -- (1,1) node[xi] {};
	\draw[kernels2] (0,0) node[not] {} -- (0,1.5) node[xi] {}; }

\DeclareSymbol{Xi4eact1}{1.5}{\draw (-1,2.5) node[xic] {} -- (-1,1) node[xi] {} -- (0,0); 
	\draw (0,0)  -- (1,1) node[xic] {};
	\draw[rho] (0,0) node[not] {} -- (0,1.5) node[xi] {}; }

\DeclareSymbol{Xi4eact2}{1.5}{\draw[rho] (-1,2.5) node[xic] {} -- (-1,1) node[xi] {} -- (0,0); 
	\draw (0,0)  -- (1,1) node[xic] {};
	\draw (0,0) node[not] {} -- (0,1.5) node[xi] {}; }

\DeclareSymbol{Xi4eabis}{1.5}{\draw (-1,2.5) node[xi] {} -- (-1,1) ; \draw[kernels2] (-1,1) node[xi] {} -- (0,0); 
	\draw (0,0)  -- (1,1) node[xi] {};
	\draw[kernels2] (0,0) node[not] {} -- (0,1.5) node[xi] {}; }

\DeclareSymbol{Xi4eabisc1}{1.5}{\draw (-1,2.5) node[xic] {} -- (-1,1) ; \draw[kernels2] (-1,1) node[xi] {} -- (0,0); 
	\draw (0,0)  -- (1,1) node[xi] {};
	\draw[kernels2] (0,0) node[not] {} -- (0,1.5) node[xic] {}; }

\DeclareSymbol{Xi4eabisc1b}{1.5}{\draw (-1,2.5) node[xic] {} -- (-1,1) ; \draw[kernels2] (-1,1) node[xi] {} -- (0,0); 
	\draw (0,0)  -- (1,1) node[xi] {};
	\draw[kernels2] (0,0) node[not] {} -- (0,1.5) node[xiesf] {}; }

\DeclareSymbol{Xi4eabisc1bis}{1.5}{\draw (-1,2.5) node[xi] {} -- (-1,1) ; \draw[kernels2] (-1,1) node[xi] {} -- (0,0); 
	\draw (0,0)  -- (1,1) node[xi] {};
	\draw[kernels2] (0,0) node[not] {} -- (0,1.5) node[xi] {};
	\draw (-2,1) node[] {\tiny{$i$}};
	\draw (-2,2.5) node[] {\tiny{$\ell$}};
	\draw (2,1) node[] {\tiny{$k$}};
	\draw (0,2.5) node[] {\tiny{$j$}};
}

\DeclareSymbol{Xi4eabisc1tris}{1.5}{\draw (-1,2.5) node[xi] {} -- (-1,1) ; \draw[kernels2] (-1,1) node[xi] {} -- (0,0); 
	\draw (0,0)  -- (1,1) node[xi] {};
	\draw[kernels2] (0,0) node[not] {} -- (0,1.5) node[xi] {};
	\draw (-2,1) node[] {\tiny{i}};
	\draw (-2,2.5) node[] {\tiny{j}};
	\draw (2,1) node[] {\tiny{j}};
	\draw (0,2.5) node[] {\tiny{i}};
}

\DeclareSymbol{Xi4eabisc1quater}{1.5}{\draw (-1,2.5) node[xic] {} -- (-1,1) ; \draw[kernels2] (-1,1) node[xi] {} -- (0,0); 
	\draw (0,0)  -- (1,1) node[xic] {};
	\draw[kernels2] (0,0) node[not] {} -- (0,1.5) node[xi] {};
}

\DeclareSymbol{Xi4eabisc2}{1.5}{\draw (-1,2.5) node[xic] {} -- (-1,1) ; \draw[kernels2] (-1,1) node[xi] {} -- (0,0); 
	\draw (0,0)  -- (1,1) node[xic] {};
	\draw[kernels2] (0,0) node[not] {} -- (0,1.5) node[xi] {}; }

\DeclareSymbol{Xi4eabisc2l}{1.5}{\draw (-1,2.5) node[xiesf] {} -- (-1,1) ; \draw[kernels2] (-1,1) node[xi] {} -- (0,0); 
	\draw (0,0)  -- (1,1) node[xic] {};
	\draw[kernels2] (0,0) node[not] {} -- (0,1.5) node[xi] {}; }

\DeclareSymbol{Xi4eabisc2r}{1.5}{\draw (-1,2.5) node[xic] {} -- (-1,1) ; \draw[kernels2] (-1,1) node[xi] {} -- (0,0); 
	\draw (0,0)  -- (1,1) node[xiesf] {};
	\draw[kernels2] (0,0) node[not] {} -- (0,1.5) node[xi] {}; }

\DeclareSymbol{Xi4eabisc3}{1.5}{\draw (-1,2.5) node[xic] {} -- (-1,1) ; \draw[kernels2] (-1,1) node[xic] {} -- (0,0); 
	\draw (0,0)  -- (1,1) node[xi] {};
	\draw[kernels2] (0,0) node[not] {} -- (0,1.5) node[xi] {}; }

\DeclareSymbol{Xi4eb}{0}{
	\draw[kernels2] (0,2) node[xi] {} -- (-1,1) ; \draw[kernels2] (-2,2)  node[xi] {} -- (-1,1) ; \draw (-1,1)  node[not] {} -- (0,0); 
	\draw (0,0) node[xi] {}  -- (1,1) node[xi] {};
}

\DeclareSymbol{Xi4eab}{1.5}{\draw[kernels2] (-1,2.5) node[xi] {} -- (-1,1) ; \draw[kernels2] (-2,2)  node[xi] {} -- (-1,1) ; \draw (-1,1)  node[not] {} -- (0,0); 
	\draw[kernels2] (0,0)  -- (1,1) node[xi] {};
	\draw[kernels2] (0,0) node[not] {} -- (0,1.5) node[xi] {}; 
}

\DeclareSymbol{Xi4eabdis}{1.5}{\draw[kernels2] (-1,2.5) node[xies] {} -- (-1,1) ; \draw[kernels2] (-2,2)  node[xies] {} -- (-1,1) ; \draw (-1,1)  node[not] {} -- (0,0); 
	\draw[kernels2] (0,0)  -- (1,1) node[xies] {};
	\draw[kernels2] (0,0) node[not] {} -- (0,1.5) node[xies] {}; 
}

\DeclareSymbol{Xi4eabc1}{1.5}{\draw[kernels2] (-1,2.5) node[xic] {} -- (-1,1) ; \draw[kernels2] (-2,2)  node[xi] {} -- (-1,1) ; \draw (-1,1)  node[not] {} -- (0,0); 
	\draw[kernels2] (0,0)  -- (1,1) node[xic] {};
	\draw[kernels2] (0,0) node[not] {} -- (0,1.5) node[xi] {}; 
}

\DeclareSymbol{Xi4eabc2}{1.5}{\draw[kernels2] (-1,2.5) node[xi] {} -- (-1,1) ; \draw[kernels2] (-2,2)  node[xi] {} -- (-1,1) ; \draw (-1,1)  node[not] {} -- (0,0); 
	\draw[kernels2] (0,0)  -- (1,1) node[xic] {};
	\draw[kernels2] (0,0) node[not] {} -- (0,1.5) node[xic] {}; 
}

\DeclareSymbol{Xi4eabbis}{1.5}{\draw[kernels2] (-1,2.5) node[xi] {} -- (-1,1) ; \draw[kernels2] (-2,2)  node[xi] {} -- (-1,1) ; \draw[kernels2] (-1,1)  node[not] {} -- (0,0); 
	\draw (0,0)  -- (1,1) node[xi] {};
	\draw[kernels2] (0,0) node[not] {} -- (0,1.5) node[xi] {}; 
}

\DeclareSymbol{Xi4eabbisc1}{1.5}{\draw[kernels2] (-1,2.5) node[xic] {} -- (-1,1) ; \draw[kernels2] (-2,2)  node[xi] {} -- (-1,1) ; \draw[kernels2] (-1,1)  node[not] {} -- (0,0); 
	\draw (0,0)  -- (1,1) node[xic] {};
	\draw[kernels2] (0,0) node[not] {} -- (0,1.5) node[xi] {}; 
}

\DeclareSymbol{Xi4eabbisc1perm}{1.5}{\draw[kernels2] (-1,2.5) node[xi] {} -- (-1,1) ; \draw[kernels2] (-2,2)  node[xic] {} -- (-1,1) ; \draw[kernels2] (-1,1)  node[not] {} -- (0,0); 
	\draw (0,0)  -- (1,1) node[xic] {};
	\draw[kernels2] (0,0) node[not] {} -- (0,1.5) node[xi] {}; 
}

\DeclareSymbol{Xi4eabbisc2}{1.5}{\draw[kernels2] (-1,2.5) node[xi] {} -- (-1,1) ; \draw[kernels2] (-2,2)  node[xi] {} -- (-1,1) ; \draw[kernels2] (-1,1)  node[not] {} -- (0,0); 
	\draw (0,0)  -- (1,1) node[xic] {};
	\draw[kernels2] (0,0) node[not] {} -- (0,1.5) node[xic] {}; 
}

\DeclareSymbol{Xi2cbis}{0}{\draw[kernels2] (0,1) -- (0.8,2.2) node[xi] {};\draw[kernels2] (0,-0.25) node[not] {} -- (0,1); \draw[kernels2] (0,1) node[not] {} -- (-0.8,2.2) node[xi] {};}

\DeclareSymbol{Xi2cbis1}{0}{\draw (0,1) -- (-0.8,2.2) node[xi] {};\draw[kernels2] (0,-0.25) node[not] {} -- (0,1) node[xi] {}; }

\DeclareSymbol{Xi2Xbis}{-2}{\draw[kernels2] (0,-0.25)  -- (-1,1) ; \draw (-1,1) node[xix] {};
	\draw[kernels2] (0,-0.25) node[not] {} -- (1,1) node[xi] {};}

\DeclareSymbol{XXi2bis}{-2}{\draw[kernels2] (0,-0.25) -- (-1,1) node[xi] {};
	\draw[kernels2] (0,-0.25) node[X] {} -- (1,1) node[xi] {};}

\DeclareSymbol{I1XiIXi}{0}{\draw[kernels2] (0,-0.25) -- (1,1) node[xi] {};
	\draw (0,-0.25) node[not] {} -- (-1,1) node[xi] {};}

\DeclareSymbol{I1XiIXib}{0}{\draw  (0,-0.25) node[xi] {} -- (0,1) node[not] {};
	\draw[kernels2] (0,1) -- (0,2.25) ; \draw (0,2.25) node[xi]{}; }

\DeclareSymbol{I1XiIXic}{0}{
	\draw[kernels2] (0,0) -- (1,1) node[xi] {} ; 
	\draw[kernels2] (0,0) node[not] {}  -- (-1,1) node[not] {} -- (0,2) node[xi] {};
}

\DeclareSymbol{thin}{1.4}{\draw[pagebackground] (-0.3,0) -- (0.3,0); \draw  (0,0) -- (0,2);}
\DeclareSymbol{thin2}{1.4}{\draw[pagebackground] (-0.3,0) -- (0.3,0); \draw[tinydots]  (0,0) -- (0,2);}

\DeclareSymbol{thick}{1.4}{\draw[pagebackground] (-0.3,0) -- (0.3,0); \draw[kernels2]  (0,0) -- (0,2);}

\DeclareSymbol{thick2}{1.4}{\draw[pagebackground] (-0.3,0) -- (0.3,0); \draw[kernels2,tinydots]  (0,0) -- (0,2);}

\DeclareSymbol{Xi4ind}{2}{\draw (0,0) node[xi,label={[label distance=-0.2em]right: \scriptsize  $ i $}]  { } -- (-1,1) node[xi,label={[label distance=-0.2em]left: \scriptsize  $ j $}] {} -- (0,2) node[xi,label={[label distance=-0.2em]right: \scriptsize  $ k $}] {} -- (-1,3) node[xi,label={[label distance=-0.2em]left: \scriptsize  $ \ell $}] {};}

\DeclareSymbol{Xi4c1}{2}{\draw (0,0) node[xic] {} -- (-1,1) node[xi] {} -- (0,2) node[xic] {} -- (-1,3) node[xi] {};} 
\DeclareSymbol{IXi2ex}{0}{\draw (0,-0.25) node[xie] {} -- (-1,1) node[xi] {} ; \draw (0,-0.25)-- (1,1) node[xi] {};}
\DeclareSymbol{IXi2ex1}{0}{\draw (0,-0.25) node[xie] {} -- (-1,1) node[xi] {} -- (0,2) node[xi] {};}

\DeclareSymbol{Xi4b1}{0}{\draw(0,1.5) node[xic] {} -- (0,0); \draw (-1,1) node[xic] {} -- (0,0) node[xi] {} -- (1,1) node[xi] {};}

\DeclareSymbol{Xi4ec1}{0}{\draw (0,2) node[xi] {} -- (-1,1) node[xic] {} -- (0,0) node[xic] {} -- (1,1) node[xi] {};}
\DeclareSymbol{Xi4ec2}{0}{\draw (0,2) node[xic] {} -- (-1,1) node[xi] {} -- (0,0) node[xic] {} -- (1,1) node[xi] {};}
\DeclareSymbol{Xi4ec3}{0}{\draw (0,2) node[xic] {} -- (-1,1) node[xic] {} -- (0,0) node[xi] {} -- (1,1) node[xi] {};}

\DeclareSymbol{I1Xi4ac1}{2}{\draw[kernels2] (0,0) node[not] {} -- (-1,1) ; \draw[kernels2] (0,0) node[not] {} -- (1,1) node[xic] {} ;
	\draw (-1,1) node[xi] {} -- (0,2) node[xic] {} -- (-1,3) node[xi] {};}

\DeclareSymbol{I1Xi4ac2}{2}{\draw[kernels2] (0,0) node[not] {} -- (-1,1) ; \draw[kernels2] (0,0) node[not] {} -- (1,1) node[xic] {} ;
	\draw (-1,1) node[xi] {} -- (0,2) node[xi] {} -- (-1,3) node[xic] {};}

\DeclareSymbol{I1Xi4bp}{2}{\draw (0,0) node[not] {} -- (-1,1) node[xi] {} -- (0,2) ; \draw[kernels2] (0,2) node[not] {} -- (-1,3) node[xi] {};\draw[kernels2] (0,2)  -- (1,3) node[xi] {};
}

\DeclareSymbol{I1Xi4bc1}{2}{\draw (0,0) node[xic] {} -- (-1,1) node[xi] {} -- (0,2) ; \draw[kernels2] (0,2) node[not] {} -- (-1,3) node[xi] {};\draw[kernels2] (0,2)  -- (1,3) node[xic] {};
}

\DeclareSymbol{I1Xi4bc2}{2}{\draw (0,0) node[xic] {} -- (-1,1) node[xi] {} -- (0,2) ; \draw[kernels2] (0,2) node[not] {} -- (-1,3) node[xic] {};\draw[kernels2] (0,2)  -- (1,3) node[xi] {};
}

\DeclareSymbol{I1Xi4cp}{2}{\draw (0,0) node[not] {} -- (-1,1) node[not] {}; \draw[kernels2] (-1,1) -- (0,2) ; 
	\draw[kernels2] (-1,1) -- (-2,2) node[xi] {} ;
	\draw (0,2) node[xi] {} -- (-1,3) node[xi] {};}

\DeclareSymbol{I1Xi4cc1}{2}{\draw (0,0) node[xic] {} -- (-1,1) node[not] {}; \draw[kernels2] (-1,1) -- (0,2) ; 
	\draw[kernels2] (-1,1) -- (-2,2) node[xi] {} ;
	\draw (0,2) node[xic] {} -- (-1,3) node[xi] {};}

\DeclareSymbol{I1Xi4cc2}{2}{\draw (0,0) node[xic] {} -- (-1,1) node[not] {}; \draw[kernels2] (-1,1) -- (0,2) ; 
	\draw[kernels2] (-1,1) -- (-2,2) node[xi] {} ;
	\draw (0,2) node[xi] {} -- (-1,3) node[xic] {};}

\DeclareSymbol{I1Xi4abc1}{2}{\draw[kernels2] (0,0) node[not] {} -- (-1,1) ; \draw[kernels2] (0,0) node[not] {} -- (1,1) node[xic] {};\draw (-1,1) node[xi] {} -- (0,2) ; \draw[kernels2] (0,2) node[not] {} -- (-1,3) node[xic] {};\draw[kernels2] (0,2)  -- (1,3) node[xi] {}; }

\DeclareSymbol{I1Xi4abc2}{2}{\draw[kernels2] (0,0) node[not] {} -- (-1,1) ; \draw[kernels2] (0,0) node[not] {} -- (1,1) node[xic] {};\draw (-1,1) node[xi] {} -- (0,2) ; \draw[kernels2] (0,2) node[not] {} -- (-1,3) node[xi] {};\draw[kernels2] (0,2)  -- (1,3) node[xic] {}; }

\DeclareSymbol{R1}{0}{\draw (-1,1) node[xi] {} -- (0,0) node[not] {};
	\draw[kernels2] (0,1.5) node[xic] {} -- (0,0) -- (1,1) node[xic] {};}
\DeclareSymbol{R2}{0}{\draw (-1,1) node[xic] {} -- (0,0) node[not] {};
	\draw[kernels2] (0,1.5)  {} -- (0,0) -- (1,1)  {};
	\draw (0,1.5) node[xi] {};
	\draw (1,1) node[xic] {};
}
\DeclareSymbol{R3}{1}{\draw[kernels2] (-1,1.5)  {} -- (0,0) node[not] {} -- (1,1.5);
	\draw (-1,1.5) node[xi] {};
	\draw[kernels2] (0,3) {} -- (1,1.5) -- (2,3)  {};
	\draw  (0,3) node[xic] {} ;
	\draw (2,3) node[xic] {};}
\DeclareSymbol{R4}{1}{\draw[kernels2] (-1,1.5) node[xic] {} -- (0,0) node[not] {} -- (1,1.5);
	\draw[kernels2] (0,3) {} -- (1,1.5) -- (2,3) node[xic] {};
	\draw (0,3) node[xi] {};}

\DeclareSymbol{I1Xi4bcp}{2}{\draw (0,0) node[not] {} -- (-1,1) node[not] {}; \draw[kernels2] (-1,1) -- (0,2) ; 
	\draw[kernels2] (-1,1) -- (-2,2) node[xi] {} ; \draw[kernels2] (0,2) node[not] {} -- (-1,3) node[xi] {};\draw[kernels2] (0,2)  -- (1,3) node[xi] {};
}

\DeclareSymbol{I1Xi4bcc1}{2}{\draw (0,0) node[xic] {} -- (-1,1) node[not] {}; \draw[kernels2] (-1,1) -- (0,2) ; 
	\draw[kernels2] (-1,1) -- (-2,2) node[xi] {} ; \draw[kernels2] (0,2) node[not] {} -- (-1,3) node[xi] {};\draw[kernels2] (0,2)  -- (1,3) node[xic] {};
}

\DeclareSymbol{I1Xi4bcc2}{2}{\draw (0,0) node[xic] {} -- (-1,1) node[not] {}; \draw[kernels2] (-1,1) -- (0,2) ; 
	\draw[kernels2] (-1,1) -- (-2,2) node[xi] {} ; \draw[kernels2] (0,2) node[not] {} -- (-1,3) node[xic] {};\draw[kernels2] (0,2)  -- (1,3) node[xi] {};
} 

\DeclareSymbol{2I1Xi4bc1}{2}{\draw[kernels2] (0,0) node[not] {} -- (-1,1) ;
	\draw[kernels2] (0,0) -- (1,1);
	\draw (-1,1) node[xic] {} -- (-1,2.5) node[xi] {};
	\draw (1,1)  node[xic] {} -- (1,2.5) node[xi] {};
}

\DeclareSymbol{2I1Xi4bc2}{2}{\draw[kernels2] (0,0) node[not] {} -- (-1,1) ;
	\draw[kernels2] (0,0) -- (1,1);
	\draw (-1,1) node[xi] {} -- (-1,2.5) node[xic] {};
	\draw (1,1)  node[xic] {} -- (1,2.5) node[xi] {};
}

\DeclareSymbol{diff2I1Xi4bc2}{2}{\draw (0,0) node[diff] {} -- (-1,1) ;
	\draw (0,0) -- (1,1);
	\draw (-1,1) node[xi] {} -- (-1,2.5) node[xic] {};
	\draw (1,1)  node[xic] {} -- (1,2.5) node[xi] {};
}

\DeclareSymbol{2I1Xi4bc3}{2}{\draw[kernels2] (0,0) node[not] {} -- (-1,1) ;
	\draw[kernels2] (0,0) -- (1,1);
	\draw (-1,1) node[xic] {} -- (-1,2.5) node[xic] {};
	\draw (1,1)  node[xi] {} -- (1,2.5) node[xi] {};
}

\DeclareSymbol{Xi41}{0}{\draw (0,1) -- (0.8,2.2) node[xic] {};\draw (0,-0.25) node[xi] {} -- (0,1) node[xi] {} -- (-0.8,2.2) node[xic] {};} 

\DeclareSymbol{Xi42}{0}{\draw (0,1) -- (0.8,2.2) node[xi] {};\draw (0,-0.25) node[xic] {} -- (0,1) node[xi] {} -- (-0.8,2.2) node[xic] {};}

\DeclareSymbol{Xi4ca1}{0}{\draw (0,1) -- (-1,2.2) node[xic] {};\draw (0,-0.25) node[xi] {} -- (0,1) ; \draw[kernels2] (0,1) node[not] {} -- (1,2.2) node[xic] {};
	\draw[kernels2] (0,1) {} -- (0,2.7) node[xi] {};
}

\DeclareSymbol{Xi4ca2}{0}{\draw (0,1) -- (-1,2.2) node[xi] {};\draw (0,-0.25) node[xi] {} -- (0,1) ; \draw[kernels2] (0,1) node[not] {} -- (1,2.2) node[xic] {};
	\draw[kernels2] (0,1) {} -- (0,2.7) node[xic] {};
}

\DeclareSymbol{Xi4cap}{0}{\draw (0,1) -- (-1,2.2) node[xi] {};\draw (0,-0.25) node[not] {} -- (0,1) ; \draw[kernels2] (0,1) node[not] {} -- (1,2.2) node[xi] {};
	\draw[kernels2] (0,1) {} -- (0,2.7) node[xi] {};
}

\DeclareSymbol{Xi3a}{0}{
	\draw (-1,1)  node[xi] {} -- (0,0); 
	\draw (0,0) node[xi] {}  -- (1,1) node[xi] {};
}

\DeclareSymbol{Xi4ebc1}{0}{
	\draw[kernels2] (0,2) node[xi] {} -- (-1,1) ; \draw[kernels2] (-2,2)  node[xic] {} -- (-1,1) ; \draw (-1,1)  node[not] {} -- (0,0); 
	\draw (0,0) node[xic] {}  -- (1,1) node[xi] {};
}

\DeclareSymbol{Xi4ebc2}{0}{
	\draw[kernels2] (0,2) node[xi] {} -- (-1,1) ; \draw[kernels2] (-2,2)  node[xi] {} -- (-1,1) ; \draw (-1,1)  node[not] {} -- (0,0); 
	\draw (0,0) node[xic] {}  -- (1,1) node[xic] {};
}

\DeclareSymbol{Xi2cbispex}{0}{\draw[kernels2] (0,1) -- (0.8,2.2) node[xi] {};\draw (0,-0.25) node[xie] {} -- (0,1); \draw[kernels2] (0,1) node[not] {} -- (-0.8,2.2) node[xi] {};}

\DeclareSymbol{Xi2cbis1p}{0}{\draw (0,1) -- (-0.8,2.2) node[xi] {};\draw (0,-0.25) node[not] {} -- (0,1) node[xi] {}; }

\DeclareSymbol{Xi2Xp}{-2}{\draw (0,-0.25) node[not] {} -- (-1,1) node[xix] {};} 

\DeclareSymbol{I1XiIXib}{0}{\draw  (0,-0.25) node[xi] {} -- (0,1) node[not] {};
	\draw[kernels2] (0,1) -- (0,2.25) ; \draw (0,2.25) node[xi]{}; }

\DeclareSymbol{IXi2b}{0}{\draw  (0,-0.25) node[xi] {} -- (0,1) node[not] {};
	\draw (0,1) -- (0,2.25) ; \draw (0,2.25) node[xi]{}; }

\DeclareSymbol{IXi2bex}{0}{\draw  (0,-0.25) node[xi] {} -- (0,1) node[xie] {};
	\draw (0,1) -- (0,2.25) ; \draw (0,2.25) node[xi]{}; }

\def\1{\mathbf{\symbol{1}}}

\def\one{\mathbf{1}}

\DeclareSymbol{diff}{0}{
	\draw (0,0.5) node[diff] {};
}

\DeclareSymbol{diff1}{0}{
	\draw (0,0.5) node[diff1] {};
}

\DeclareSymbol{diff2}{0}{
	\draw (0,0.5) node[diff2] {};
}

\DeclareSymbol{geo}{0}{
	\draw (0,0) node[diff] {};
	\draw (0.3,0) node[diff] {};
}

\DeclareSymbol{generic}{0}{
	\draw (0,0.6) node[xi] {};
}

\DeclareSymbol{g}{0}{
	\draw (0,0.6) node[g] {};
}

\DeclareSymbol{Ito}{0}{
	\draw (0,0.6) node[xies] {};
}

\DeclareSymbol{Itob}{0}{
	\draw (0,0.6) node[xiesf] {};
}

\DeclareSymbol{greycirc}{0}{
	\draw (0,0.3) node[xi] {};
}

\DeclareSymbol{not}{0}{
	\draw (0,0.6) node[not] {};
	\draw[tinydots] (0,0.6) circle (0.8);
}

\DeclareSymbol{genericb}{0}{
	\draw (0,0.6) node[xic] {};
}

\DeclareSymbol{bluecirc}{0}{
	\draw (0,0.3) node[xic] {};
}

\DeclareSymbol{genericxix}{0}{
	\draw (0,0.6) node[xix] {};
}

\DeclareSymbol{genericX}{0}{
	\draw (0,0.6) node[X] {};
}

\DeclareSymbol{diffIto}{1}{
	\draw  (0,2.5) -- (0,0) ;
	\draw (0,-0.1) node[diff] {};
	\draw (0,2.5) node[xies] {};
}
\DeclareSymbol{Itodiff}{2}{
	\draw(0,2.9) -- (0,-0.2);
	\draw (0,2.9) node[diff] {};
	\draw (0,-0.1) node[xies] {};
}

\DeclareSymbol{diffgeneric}{1}{
	\draw  (0,2.5) -- (0,0) ;
	\draw (0,-0.1) node[diff] {};
	\draw (0,2.5) node[xi] {};
}

\DeclareSymbol{genericdiff}{2}{
	\draw(0,2.9) -- (0,-0.2);
	\draw (0,2.9) node[diff] {};
	\draw (0,-0.1) node[xi] {};
}

\DeclareSymbol{diffdot}{2}{
	\draw  (0,3) -- (0,-0.1) ;
	\draw (0,3) node[not] {};
	\draw (0,-0.1) node[diff] {};
}

\DeclareSymbol{diffdotmini}{0}{
	\draw  (0,0) -- (0,1.2) ;
	\draw (0,1.2) node[not] {};
	\draw (0,0) node[diffmini] {};
}

\DeclareSymbol{dotdiff}{2}{
	\draw[kernelsmod]  (0,3) -- (0,-0.1) ;
	\draw (0,3) node[diff] {};
	\draw (0,-0.1) node[not] {};
}

\DeclareSymbol{dotdiff1}{2}{
	\draw[kernelsmod]  (0,3) -- (0,-0.1) ;
	\draw (0,3) node[diff1] {};
	\draw (0,-0.1) node[not] {};
}

\DeclareSymbol{dotdiff1mini}{0}{
	\draw[kernelsmod]  (0,1.2) -- (0,0) ;
	\draw (0,1.2) node[diffmini] {};
	\draw (0,0) node[not] {};
}

\DeclareSymbol{dotdiff2}{2}{
	\draw (0,3) -- (0,-0.1) ;
	\draw (0,3) node[diff] {};
	\draw (0,-0.1) node[not] {};
}

\DeclareSymbol{dotdiff2mini}{0}{
	\draw (0,1.2) -- (0,0) ;
	\draw (0,1.2) node[diffmini] {};
	\draw (0,0) node[not] {};
}

\DeclareSymbol{dotdiffstraight}{0}{
	\draw  (0,3) -- (0,-0.1) ;
	\draw (0,3) node[diff] {};
	\draw (0,-0.1) node[not] {};
}

\DeclareSymbol{arbre1}{0}{
	\draw  (0,0) -- (1.5,1.5) ;
	\draw (1.5,1.5) node[not] {};
	\draw (0,0) node[not] {};
}

\DeclareSymbol{arbre2}{0}{
	\draw  (0,0) -- (1.5,1.5) ;
	\draw[kernelsmod] (0,0) -- (-1.5,1.5);
	\draw (1.5,1.5) node[not] {};
	\draw (0,0) node[not] {};
	\draw (-1.5,1.5) node[xi] {};
}

\DeclareSymbol{arbre3}{0}{
	\draw  (0,0) -- (1.5,1.5) ;
	\draw[kernelsmod] (1.5,1.5) -- (0,3);
	\draw (0,0) node[not] {};
	\draw (1.5,1.5) node[not] {};
	\draw (0,3) node[xi] {};
}

\DeclareSymbol{treeeval}{0}{
	\draw (0,0) -- (1,1);
	\draw (0,0) node[xi] {};
	\draw (1.25,1.25) node[xi] {};
	\draw (-0.6,0.6) node[]{\tiny{$i$}};
	\draw (0.65,1.85) node[]{\tiny{$j$}};
}

\DeclareSymbol{testeval}{0}{
	\draw (0,0) -- (1,1);
	\draw (0,0) -- (-1,1);
	\draw (0,0) node[xi] {};
	\draw (1.25,1.25) node[xi] {};
	\draw (-1.25,1.25) node[xi] {};
	\draw (-0.6,-0.6) node[]{\tiny{$i$}};
	\draw (0.65,1.85) node[]{\tiny{$j$}};
	\draw (-1.95,1.85) node[]{\tiny{$k$}};
}

\DeclareSymbol{treeeval2}{0}{
	\draw[kernelsmod] (-0.25,-1) -- (1,0.5) ;
	\draw[kernelsmod] (1,0.5) -- (-0.25,2);
	\draw (1,0.5) node[diff2] {};
	\draw (-0.25,-1) node[not] {};
	\draw (-0.25,2) node[xi] {};
	\draw (-0.6,1.2) node[]{\tiny{1}};
}

\DeclareSymbol{arbreact}{1}{
	\draw (0,0) node[not] {};
	\draw[kernelsmod] (0,0) -- (1,1);
	\draw[kernelsmod] (0,0) -- (-1,1);
	\draw (-1,1) node[xic] {};
	\draw  (0,2) -- (1,1) ;
	\draw (0,2) node[xic] {};
	\draw (1,1) node[xi] {};
}

\DeclareSymbol{arbreact1}{0}{
	\draw (0,-1.5) -- (0,0);
	\draw[kernelsmod] (0,0) -- (1,1);
	\draw[kernelsmod] (0,0) -- (-1,1);
	\draw  (0,2) -- (1,1) ;
	\draw (0,-1.5) node[diff] {};
	\draw (0,0) node[not] {};
	\draw (-1,1) node[xic] {};
	\draw (0,2) node[xic] {};
	\draw (1,1) node[xi] {};
}

\DeclareSymbol{arbreact2}{0}{
	\draw (0,-0.75) -- (-1,0.5); 
	\draw (0,-0.75) -- (1,0.5);
	\draw (0,1.5) -- (1,0.5);
	\draw (0,1.5) node[xic] {};
	\draw (1,0.5) node[xi] {};
	\draw (-1,0.5) node[xic] {};
	\draw (0,-0.75) node[diff] {};
}

\DeclareSymbol{arbreact3}{0}{
	\draw[kernelsmod] (0,-0.75) -- (-1,0.5); 
	\draw[kernelsmod] (0,-0.75) -- (1,0.5);
	\draw (0,1.75) -- (1,0.5);
	\draw (2,1.75) -- (1,0.5);
	\draw (0,1.75) node[xic] {};
	\draw (1,0.5) node[diff] {};
	\draw (-1,0.5) node[xic] {};
	\draw (2,1.75) node[xi] {};
	\draw (0,-0.75) node[not] {};
}

\DeclareSymbol{pre_im_I1Xitwo}{0}{
	\draw[kernels2] (0,-0.3) node[not] {} -- (-0.6,0.7) ;
	\draw[kernels2] (0,-0.3) -- (0.6,0.7);
	\draw (0,0.9) node[g] {};
}

\DeclareSymbol{pre_im_cI1Xi4ab}{2}{
	\draw[kernels2] (0,-1) node[not] {} -- (-0.6,0) ;
	\draw[kernels2] (0,-1) -- (0.6,0);
	\draw (0,0.2) node[g] {};
	\draw (0,0.6) -- (0,1.5);
	\draw[kernels2] (0,1.5) node[not] {} -- (-0.6,2.5) ;
	\draw[kernels2] (0,1.5) -- (0.6,2.5);
	\draw (0,2.7) node[g] {};
}

\DeclareSymbol{pre_im_I1Xi4acc2}{0}{
	\draw[kernels2] (-1,-0.5) node[not] {} -- (-1.6,0.5) ;
	\draw[kernels2] (-1,-0.5) -- (-0.4,0.5);
	\draw[kernels2] (-1,-0.5) -- (0.2,-1.5) node[not] {} ;
	\draw (-1,1.1) -- (-1,2);
	\draw[kernels2] (0.2,-1.5) -- (0.2,2);
	\draw (-1,0.7) node[g] {};
	\draw (-0.3,2.2) node[g] {};
}

\DeclareSymbol{pre_im_I1Xi4abcc2}{2}{
	\draw[kernels2] (0,-1) node[not] {} -- (-1,0) node[not] {};
	\draw[kernels2] (-1,1.2) node[not] {} -- (-1,0);
	\draw[kernels2] (-1,1.2) -- (-1.5,2.5);
	\draw[kernels2] (-1,1.2) -- (-0.5,2.5);
	\draw[kernels2] (-1,0) -- (0.7,2.5);
	\draw[kernels2] (0,-1) -- (1.5,2.5);
	\draw (-1,2.7) node[g] {};
	\draw (1,2.7) node[g] {};
}

\DeclareSymbol{pre_im_2I1Xi4c1}{2}{
	\draw[kernels2] (0,-0.5) node[not] {} -- (-1,0.5) node[not] {};
	\draw[kernels2] (0,-0.5) -- (1,0.5) node[not] {};
	\draw[kernels2] (-1,0.5) node[not] {}-- (-1.7,2);
	\draw[kernels2]  (-1,2) -- (1,0.5);
	\draw[kernels2] (-1,0.5) -- (1,2);
	\draw[kernels2] (1,0.5) -- (1.7,2);
	\draw (-1.2,2.2) node[g] {};
	\draw (1.2,2.2) node[g] {};
}

\DeclareSymbol{pre_im_Xi4eabisc2}{2}{
	\draw[kernels2] (1.2,-0.5) node[not] {} -- (-0.7,0.8) ;
	\draw[kernels2] (1.2,-0.5) -- (0.4,0.8);
	\draw (0,1.4)  -- (0,2.2);
	\draw (1.2,2.2) -- (1.2,-0.6);
	\draw (0,1) node[g] {};
	\draw (0.6,2.4) node[g] {};
}

\DeclareSymbol{pre_im_Xi4eabisc22}{2}{
	\draw (1.2,-0.5) node[not] {} -- (-0.7,0.8) ;
	\draw[kernels2] (1.2,-0.5) -- (0.4,0.8);
	\draw (0,1.4)  -- (0,2.2);
	\draw[kernels2] (1.2,2.2) -- (1.2,-0.6);
	\draw (0,1) node[g] {};
	\draw (0.6,2.4) node[g] {};
}

\DeclareSymbol{pre_im_Xi4eabisc222}{2}{
	\draw[kernels2] (0.4,-0.5) node[not] {} -- (-0.6,1) ;
	\draw[kernels2] (1.2,0) -- (0.3,1);
	\draw (0,1.1)  -- (0,2.5);
	\draw[kernels2] (1.2,2.5) -- (1.2,0) node[not] {} -- (0.4,-0.6);
	\draw (0,1.2) node[g] {};
	\draw (0.6,2.5) node[g] {};
}

\DeclareSymbol{pre_im_Xi4eabc2}{2}{
	\draw (0,-0.5) node[not] {} -- (-1,0.5) node[not] {};
	\draw[kernels2] (-1,0.5) -- (-1.5,2);
	\draw[kernels2] (0,-0.5)  -- (0.7,2);
	\draw[kernels2] (-1,0.5) -- (-0.5,2);
	\draw[kernels2] (0,-0.5) -- (1.5,2);
	\draw (-1,2.2) node[g] {};
	\draw (1,2.2) node[g] {};
}

\DeclareSymbol{pre_im_Xi4eabbisc2}{2}{
	\draw[kernels2] (0,-0.5) node[not] {} -- (-1,0.5) node[not] {};
	\draw[kernels2] (-1,0.5) -- (-1.5,2);
	\draw[kernels2] (0,-0.5)  -- (0.7,2);
	\draw[kernels2] (-1,0.5) -- (-0.5,2);
	\draw (0,-0.5) -- (1.5,2);
	\draw (-1,2.2) node[g] {};
	\draw (1,2.2) node[g] {};
}

\DeclareSymbol{pre_im_I1Xi4abcc1}{2}{
	\draw[kernels2] (0,-1) node[not] {} -- (-1,0) node[not] {};
	\draw[kernels2] (0,1.1) node[not] {} -- (-1,0);
	\draw[kernels2] (-1,0) -- (-1.5,2.5);
	\draw[kernels2] (0,1.1) node[not] {} -- (-0.5,2.5);
	\draw[kernels2] (0,1.1) -- (0.5,2.5);
	\draw[kernels2] (0,-1) -- (1.5,2.5);
	\draw (-1,2.7) node[g] {};
	\draw (1,2.7) node[g] {};
}

\DeclareSymbol{pre_im_Xi4eabc1}{2}{
	\draw (0,-0.5) node[not] {} -- (-1,0.5) node[not] {};
	\draw[kernels2] (-1,0.5) -- (-1.5,2);
	\draw[kernels2] (0,-0.5)  -- (-0.5,2);
	\draw[kernels2] (-1,0.5) -- (0.8,2);
	\draw[kernels2] (0,-0.5) -- (1.5,2);
	\draw (-1,2.2) node[g] {};
	\draw (1,2.2) node[g] {};
}

\DeclareSymbol{pre_im_Xi4ba1b}{2}{
	\draw[kernels2] (0,0) node[not] {}  -- (1.8,1.5);
	\draw[kernels2] (0,0) -- (0.8,1.5);
	\draw (0,-0.1) -- (-1.8,1.5);
	\draw (0,-0.1) -- (-0.8,1.5);
	\draw (-1,1.7) node[g] {};
	\draw (1,1.7) node[g] {};
}

\DeclareSymbol{pre_im_Xi4ba2}{2}{
	\draw (0,-0.1) node[not] {}  -- (1.8,1.5);
	\draw[kernels2] (0,0) -- (0.8,1.5);
	\draw (0,-0.1) -- (-1.8,1.5);
	\draw[kernels2] (0,0) -- (-0.8,1.5);
	\draw (-1,1.7) node[g] {};
	\draw (1,1.7) node[g] {};
}

\DeclareSymbol{pre_im_Xi4cabc2}{2}{
	\draw[kernels2] (0,-0.5) node[not] {} -- (-1,0.5) node[not] {};
	\draw[kernels2] (-1,0.5) -- (-1.5,2);
	\draw (-1,0.5)  -- (0.7,2);
	\draw[kernels2] (-1,0.5) -- (-0.5,2);
	\draw[kernels2] (0,-0.5) -- (1.5,2);
	\draw (-1,2.2) node[g] {};
	\draw (1,2.2) node[g] {};
}

\DeclareSymbol{pre_im_Xi4cabc1}{2}{
	\draw[kernels2] (0,-0.5) node[not] {} -- (-1,0.5) node[not] {};
	\draw (-1,0.5) -- (-1.5,2);
	\draw[kernels2] (-1,0.5)  -- (0.7,2);
	\draw[kernels2] (-1,0.5) -- (-0.5,2);
	\draw[kernels2] (0,-0.5) -- (1.5,2);
	\draw (-1,2.2) node[g] {};
	\draw (1,2.2) node[g] {};
}

\DeclareSymbol{pre_im_Xi4eabbisc1}{2}{
	\draw[kernels2] (0,-0.5) node[not] {} -- (-1,0.5) node[not] {};
	\draw[kernels2] (-1,0.5) -- (-1.5,2);
	\draw[kernels2] (0,-0.5)  -- (-0.5,2);
	\draw[kernels2] (-1,0.5) -- (0.8,2);
	\draw (0,-0.5) -- (1.5,2);
	\draw (-1,2.2) node[g] {};
	\draw (1,2.2) node[g] {};
}

\DeclareSymbol{pre_im_1}{0}{
	\draw[kernels2] (0,-0.5) node[not] {} -- (-0.6,0.5) ;
	\draw[kernels2] (0,-0.5) -- (0.6,0.5);
	\draw (0,1.1)  -- (-0.55,2);
	\draw (0,1.1)  -- (0.55,2);
	\draw (0,0.7) node[g] {};
	\draw (0,2.2) node[g] {};
}

\DeclareSymbol{disconnect}{0}{
	\draw[kernels2] (0,-0.5) node[not] {} -- (-0.6,0.5) ;
	\draw[kernels2] (0,-0.5) -- (0.6,0.5);
	\draw (-0.55,1.1)  -- (-0.55,2.3);
	\draw (0.55,2.3) -- (0.55,1.5) -- (1.2,1.5) -- (1.2,3.5) -- (0.55,3.5) -- (0.55,2.7);
	\draw (0,0.7) node[g] {};
	\draw (0,2.5) node[g] {};
}

\DeclareSymbol{pre_im_2}{2}{\draw[kernels2] (0,0) node[not] {} -- (-1,1) node[not] {};
	\draw[kernels2] (0,0) -- (1,1) node[not] {};
	\draw[kernels2] (-1,1) -- (-1.5,2.5);
	\draw[kernels2] (-1,1) -- (-0.5,2.5);
	\draw[kernels2] (1,1) -- (0.5,2.5);
	\draw[kernels2] (1,1) -- (1.5,2.5);
	\draw (-1,2.7) node[g] {};
	\draw (1,2.7) node[g] {};
}

\DeclareSymbol{CX_rec}{0}{
	\draw [black] (-0.3,1) to (-0.3,-0.3);
	\draw [black] (0.3,1) to (0.3,-0.3);
	\draw [black] (-0.3,1) to (-0.3,2.3);
	\draw [black] (0.3,1) to (0.3,2.3);
	\draw (0,1) node[rec] {};
}

\DeclareSymbol{CX_cerc}{0}{
	\draw [black] (0,1) to (0,-0.3);
	\draw (0,1) node[cerc] {};
}

\DeclareMathAlphabet{\mathpzc}{OT1}{pzc}{m}{it}

\def\eqref#1{(\ref{#1})}

\makeatletter 
\newcommand*{\bigcdot}{}
\DeclareRobustCommand*{\bigcdot}{%
	\mathbin{\mathpalette\bigcdot@{}}%
}
\newcommand*{\bigcdot@scalefactor}{.5}
\newcommand*{\bigcdot@widthfactor}{1.15}
\newcommand*{\bigcdot@}[2]{%
	\sbox0{$#1\vcenter{}$}
	\sbox2{$#1\cdot\m@th$}%
	\hbox to \bigcdot@widthfactor\wd2{%
		\hfil
		\raise\ht0\hbox{%
			\scalebox{\bigcdot@scalefactor}{%
				\lower\ht0\hbox{$#1\bullet\m@th$}%
			}%
		}%
		\hfil
	}%
}
\makeatother

\def\two{{\<generic>\kern0.05em\<genericb>}}
\def\twoI{{\<Ito>\kern0.05em\<Itob>}}

\begin{document}







\begin{center}
{\LARGE\sffamily{\textbf{Locality for singular stochastic PDEs}   \vspace{0.5cm}}}
\end{center}

\begin{center}
{\sf 
I. Bailleul \footnote{Univ Brest, CNRS UMR 6205, Laboratoire de Math\'ematiques de Bretagne Atlantique, France. {\it Email}: ismael.bailleul@univ-brest.fr}, 
Y. Bruned\footnote{IECL, UMR 7502, Universit\'e de Lorraine, France. {\it Email}:  yvain.bruned@univ-lorraine.fr} 
}
\end{center}

\vspace{1cm}

\begin{center}
\begin{minipage}{0.8\textwidth}
\renewcommand\baselinestretch{0.7} \scriptsize \textbf{\textsf{\noindent Abstract.}}  This work deals with singular stochastic PDEs driven by non-translation invariant differential operators. We describe the renormalized equation for a very large class of spacetime dependent renormalization schemes. Our approach bypasses in particular the use of decorated trees with extended decorations.
\end{minipage}
\end{center}

\bigskip


\section{Introduction}

Denote by $x=(x_0,x_1)\in\bbR\times\textbf{\textsf{T}}$ a typical spacetime point over the one dimensional periodic torus $\textbf{\textsf{T}}$. Let 
$$
\mcL^{p} v \defeq a^p(\cdot)\partial_{x_1}^2 v \qquad (1\leq p\leq p_0),
$$
stand for a finite family of second order differential operators with variable smooth coefficients. We consider systems of parabolic equations of the form
\begin{equation} \label{EqSPDE}
\big(\partial_{x_0} - \mcL^p\big)u_p = f_p(u)\xi + g_p(x,u,\partial_{x_1} u)   \qquad   (1\leq p\leq p_0),
\end{equation}
with $u \defeq (u_1,\dots,u_{p_0})$ and each $u_p$ taking values in a finite dimensional space $\bbR^{d_p}$, and with $\xi=(\xi_1,\dots,\xi_{l_0})$ an $l_0$-dimensional spacetime `noise'. The quantity $f_p(u)\xi$ is thus to be understood as $\sum_{l=1}^{l_0} f_p^l(u)\xi_l$. An initial condition $u(0)$ of positive H\"older regularity is given and fixed throughout that work. We consider the question of the renormalization of this system of non-translation invariant singular stochastic PDEs.

\ssk

Non-translation invariance refers to the fact that the operators $\mcL^p$ in the system \eqref{EqSPDE}  may have non-constant coefficients. Working in a non-translation invariant setting does not influence the algebraic side of the theory of regularity structures encoding the mechanics of the local description of modelled distributions/functions. We can use in the non-translation invariant setting the same algebraic structure, with no extended decorations, as in the translation invariant setting. The translation invariance of the operators involved in the equations studied so far using regularity structures manifests itself in two different ways. On a technical level one uses of a number of translation invariant kernels $q(x-y)$ rather than just kernels $q(x,y)$. On a more fundamental level, it manifests itself in the renormalization process in the fact that renormalization {\it constants} are involved  rather than renormalization {\it functions}, as expected in a non-translation invariant setting, as found for instance in the paracontrolled approach to singular stochastic PDEs developed in \cite{BB1, BB2, BB3}. In  a non-translation invariant setting the very formulation of the renormalization scheme of Bruned, Hairer and Zambotti \cite{BHZ} does not make obvious sense anymore. Denote by $\sf M = (\Pi, g)$ the canonical smooth admissible model associated with a given system of singular stochastic PDEs driven by a continuous noise, acting on the regularity structure $\big((T,\Delta),(T^+,\Delta^+)\big)$ associated with that system after Bruned, Hairer and Zambotti's work \cite{BHZ}. The renormalized admissible continuous models built in \cite{BHZ} are associated with some characters $\ell$ of an algebra to which one  associates some linear maps $M_\ell : T\rightarrow T$ and some renormalized interpretation operators
$$
{\sf \Pi}_\ell \defeq {\sf \Pi}\circ M_\ell.
$$
The map ${\sf \Pi}_\ell$ defines a unique admissible model. The transform ${\sf \Pi} \mapsto {\sf \Pi}\circ M$ is called a renormalization schmeme. There is no canonical way of associating a model to a non-conitnuous noise $\xi$. One then regularizes the noise into a continuous noise $\xi^\epsilon$ so it has an associated canonical model ${\sf \Pi}^\epsilon$. This model typically diverges as $\epsilon>0$ goes to $0$. One then look for an $\epsilon$-dependent character $\ell^\epsilon$ such that the model associated with the renormalized interpretation map ${\sf \Pi}^\epsilon\circ M_{\ell^\epsilon}$ is converging.

Such a renormalization scheme cannot work in a non-translation invariant setting, even by replacing $\ell$ by a character-valued function $\ell(x)$ of the state space variable $x$. Indeed such a map cannot account for some possible internal renormalization indexed by functions of the internal variables of the multiple integral expressions defining the renormalized interpretation map, not by the base/external variable $x$. Let us illustrate with an example why we need to have a local renormalisation. We consider the following integral
	\begin{equation*}
		\int K(y-x_1) G(x_1,x_2) F(x_2)   d x_1 d x_2.
	\end{equation*}
	We suppose that $ G $ is singular along the diagonal $x_1=x_2$ but it is not translation invariant. Therefore one has
	\begin{equation*}
		\int G(x_1,x_2) dx_2 \neq 	\int G(0,x_2) dx_2 
	\end{equation*}
	In order to make sense of the previous integral one has to renormalize it; in this case one has to subtract a Taylor jet to $F$ at the point $ x_1 $. This gives the renormalized integral 
	\begin{equation*}
	\int K(y-x_1) G(x_1,x_2) \big( F(x_2) - F(x_1) \big) d x_1 d x_2
	\end{equation*}
	where we have supposed that subtracting only $ F(x_1) $ is enough for counter-balancing the singularity of $ G(x_1,x_2) $. The counter-term is
	\begin{equation*}
		\int K(y-x_1) F(x_1) G(x_1,x_2)  d x_1 d x_2.
	\end{equation*}
	In that case the renormalisation function is given by 
	\begin{equation*}
		\int G(x_1,x_2) d x_2  
	\end{equation*}
	which is clearly not equal to $ 	\int G(0,x_2) d x_2  $. This situation will happen when one wants to have bounds on the moments of some stochastic iterated integrals. These moments will be given by some iterated integrals of the form above with kernels singular on some diagonals.   
	  
So a different approach of renormalization is needed in a non-translation invariant setting. It is one of the aims of the present work to provide such an alternative. This is first done in Section \ref{SectionLocalProduct} and Section \ref{SectionShortProof} in a translation invariant setting after we set the algebraic background in Section \ref{SectionDecoDeformedPreLie}. {\it Our approach involves in particular no extended decorations and does not rest on the pre-Lie mechanics} that was the workhorse of Bruned, Chandra, Chevyrev \& Hairer in their seminal work \cite{BCCH18}. Instead we rely on some important algebraic properties of a family of maps called preparation maps and a new right morphism property that they satisfy. This simplifies the general setting of regularity structures. We turn in Section \ref{SectionRenomalisationSchemes} to the general setting of systems of singular stochastic PDEs involving non-translation invariant differential operators. Section \ref{SectionAnalytic} describes the analytic ingredient needed to run the machinery of regularity structures in a non-translation invariant setting. The natural extension of our approach to the renormalized equation works perfectly in a non-translation invariant setting. This involves state space dependent preparation maps and their associated models, both of which are discussed in Section \ref{SubsectionStateSpaceDependentPreparationMaps}. Our main result reads as follows. The notions involved here are explained in the body of the text.

\begin{theorem} \label{ThmMain}
 Let $\xi$ be a continuous noise. Let $R : (\bbR\times\textbf{\textsf{T}})\times T\rightarrow T$ be a strong preparation map such that $R\tau=\tau$ for all polynomials and planted trees $\tau$ in $\bigoplus_{a\in\frak{T}^+\times\{0\}}\mathcal{I}_a(T)$. Denote by ${\sf M}^{\!R}$ the admissible model associated with $R$ and write ${\sf R}^{{\sf M}^{\!R}}$ for the corresponding reconstruction map. Let ${\sf u}^{\!R}$ stand for the modelled distribution solution of the regularity structure lift of system \eqref{EqSPDE} associated with ${\sf M}^{\!R}$, with initial condition $u(0)$. Then 
$$
u \defeq {\sf R}^{{\sf M}^{\!R}}( {\sf u}^{\!R})
$$ 
is a solution of the renormalized system
\begin{equation} \label{EqRenormalizedSystemIntro}
\big(\partial_{x_0} - \mcL^p\big) u_p = f_p(u)\xi + g_p(u,\partial_{x_1} u) + \sum_{l=0}^{l_0} \Upsilon_p\Big(\big(R(\cdot)^* - \textrm{\emph{Id}}\big) \zeta_l\Big)(u,\partial_{x_1} u)\,\xi_l \qquad (1\leq p\leq p_0),
\end{equation}
for some explicit functions $\Upsilon_p(\tau)(u,\partial_{x_1}u)$ indexed by $\tau\in T$.
\end{theorem}

 The reader will see from the proof below that this statement holds in a multidimensional space setting as well. Branched rough paths provide a particular example of simple regularity structures whose elements have no polynomial decorations. The class of renormalization maps of branched rough paths is indexed by a set of multi-pre-Lie morphisms. We prove in Section \ref{SectionRenormalizationRP} that this class is exactly parametrized by the set of preparation maps. 

We finish this introduction by a point of methodology. So far, in a singular SPDE setting, renormalization was mainly viewed as the task of finding some linear maps $M$ of the model space $T$ that turn any admissible model $\sf \Pi$ into an admissible model ${\sf \Pi}\circ M$. Our introduction of preparation maps has a different aim. It provides a mechanics for {\it building} inductively some admissible models.

\bigskip

\noindent \textbf{Notations --} {\it The letter $\zeta$ will be used exclusively for the noise symbol in a regularity structure. The letters $\sigma, \tau, \mu, \nu$ will denote some decorated trees.}

\ssk

{\it Define the parabolic distance $d$ on $\bbR\times\bbR$ by
$$
d(x,x') \defeq \sqrt{\vert x_0-x_0'\vert} + \vert x_1-x_1'\vert,
$$
for some arbitrary points $x=(x_0,x_1), x'=(x_0',x_1')$ in $\bbR^2$. For $k=(k_0,k_1)\in\bbN\times\bbN$ set
$$
\vert k\vert_{\frak{s}} \defeq 2k_0+k_1.
$$
}

\medskip

\section{Decorated trees}
\label{SectionDecoDeformedPreLie}

We recall from Bruned, Hairer \& Zambotti's work \cite{BHZ} and  Bruned \& Manchon's work \cite{BrunedManchon} what we need from the combinatorial setting used to build regularity structures associated with some system of singular stochastic PDEs. A fundamental morphism property of the $\star$ product introduced below is proved in Proposition \ref{star_morphism}. This statement will play later the role played by a fundamental multi-pre-Lie morphism introduced in \cite{BCCH18}. Our proof of the renormalized equation will in particular not involve any pre-Lie morphisms. This will be crucial for the results of Section \ref{SectionRenomalisationSchemes} that hold in the general setting of systems of singular stochastic PDEs with non-translation invariant differential operators, for which the renormalization strategy of \cite{BHZ} has no obvious counterpart.

\bigskip

Recall system \eqref{EqSPDE} with its noises $\xi_1,\dots,\xi_{l_0}$ and its operators $\mcL^1,\dots, \mcL^{p_0}$. Let
$$
 \mathfrak{T}_-=\big(\Labhom^1_-,\dots,\Labhom^{l_0}_-\big), \;\textrm{and}\; \mathfrak{T}_+ = \big(\Labhom^1_+,\dots,\Labhom^{p_0}_+\big)
$$ 
be some finite sets representing the noise types and the operator types, respectively. Denote by 
$$
\mathfrak{T} \defeq  \mathfrak{T}_-\cup\mathfrak{T}_+
$$
the set of all types. We consider some decorated trees $(\tau,\Labn,\Labe)$ where $\tau$ is a non-planar rooted tree with node set $N_\tau$ and edge set $E_\tau$. The maps $\Labn : N_\tau \rightarrow \N^{d+1}$ and $\Labe=\big(\frak{t}(\cdot), \frak{q}(\cdot)\big): E_\tau \rightarrow  \mathfrak{T} \times \N^{d+1}$ are node decorations and edge decorations, respectively. The $\N^{d+1}$-part $\frak{q}(e)$ in the edge decoration of an edge $e$ encodes some possible derivatives acting on the operator associated with the given edge type $\frak{t}(e)$. Let us stress that we describe here a setting where the noises are represented by edges. (There is a slightly different setting where noises are represented by a vertex type decoration.) We will frequently abuse notations and simply denote by $\tau$ a decorated tree, using a symbolic notation. 
\begin{enumerate}
   \item[$\bullet$] An edge decorated by  $ (\Labhom,q) \in \mathfrak{T} \times \N^{d+1} $  is denoted by $\CI_{(\Labhom,q)}$. The symbol $\CI_{(\Labhom,q)}$ is also viewed as  the operation that grafts a tree onto a new root via a new edge with edge decoration  $(\Labhom,q)$.
    \item[$\bullet$] A factor $ X^k $ encodes a single node decorated by $k \in \N^{d+1}$. Denote by 
    $$
    \big(\epsilon_1, \ldots, \epsilon_{d+1}\big)
    $$ 
    the canonical basis of $ \N^{d+1} $. For $1\leq i\leq d+1$, write $ X_i $ for $ X^{e_i} $. The element $ X^0 $ is identified with $ \one $.
\end{enumerate}
 We have $d=1$ and second order differential operators $\mcL^p$ in the particular setting of system \eqref{EqSPDE}. We describe below a setting adapted to deal with the general case where $d\geq 1$ and the operators $\mcL^p$ may be of order different from $2$. We require that every decorated tree $ \tau $ contains at most one edge decorated by $(\Labhom,q)$ with $ \Labhom \in {\color{red} \mathfrak{T}_-}  $ and any $q\in\N^{d+1}$, at each node. This encodes the fact that no product of two noises are involved in the analysis of the system \eqref{EqSPDE}. We suppose that these edges of noise type lead directly to leaves of the tree; we denote them by $\zeta_l$, for $1\leq l\leq l_0$. By convention we set
$$ 
\zeta_0 \defeq \one.
$$  
Any decorated tree $ \tau $ has a unique decomposition
\[
	\tau =  X^{k} \zeta_l \prod_{i=1}^{n} \CI_{a_i}(\tau_i) ,
\]
where $0\leq l\leq l_0$ and $\prod_i$ is the tree product, the $\tau_i$ are decorated trees and the $a_i$ belong to $ \mathfrak{T}_+  \times \N^{d+1}$, so no factor in the product is a noise symbol $\zeta_{l'}$.  The root of the particular tree $X^k\zeta_l$, with $1\leq l\leq l_0$, has for instance an $\frak{n}$-decoration equal to $k$. Here is an example of a decorated tree
	\begin{equation*}
	 X^{\beta} \zeta_2 \mathcal{I}_{b}(X^{\gamma}\zeta_3) 
	 \mathcal{I}_{a}(X^{\gamma}\zeta_1)
	 = 	\begin{tikzpicture}[scale=0.2,baseline=0.1cm]
			\node at (0,0)  [dot,label= {[label distance=-0.2em]below: \scriptsize  $  \beta    $} ] (root) {};
			\node at (0,5)  [dot,label= {[label distance=-0.2em]right: \scriptsize  $  \gamma   $} ] (center) {};
			\node at (0,9)  [dot,label= {[label distance=-0.2em]above: \scriptsize  $    $} ] (centerc) {};
			\node at (3,4)  [dot,label={[label distance=-0.2em]right: \scriptsize  $ \alpha $}] (right) {};
			\node at (-3,4)  [dot,label={[label distance=-0.2em]above: \scriptsize  $ $} ] (left) {};
			\node at (3,8)  [dot,label={[label distance=-0.2em]above: \scriptsize  $ $} ] (rightc) {};
			\draw[kernel1] (right) to
			node [sloped,below] {\small }     (root); 
			\draw[kernel1] (center) to
			node [sloped,below] {\small }     (root);
			\draw[kernel1] (center) to
			node [sloped,below] {\small }     (centerc);
			\draw[kernel1] (right) to
			node [sloped,below] {\small }     (rightc); 
			\draw[kernel1] (left) to
			node [sloped,below] {\small }     (root);
			\node at (3,6) [fill=white,label={[label distance=0em]center: \scriptsize  $ \zeta_1 $} ] () {};
			\node at (0,7) [fill=white,label={[label distance=0em]center: \scriptsize  $ \zeta_3 $} ] () {};
			\node at (-2,2) [fill=white,label={[label distance=0em]center: \scriptsize  $ \zeta_2 $} ] () {};
			\node at (2,2) [fill=white,label={[label distance=0em]center: \scriptsize   $ a  $} ] () {};
			\node at (0,2.5) [fill=white,label={[label distance=0em]center: \scriptsize  $ b  $} ] () {};
		\end{tikzpicture}
	\end{equation*}
where we have reprented the decoration on the edges and the nodes. For a noise edge, we have chosen to put $ \zeta_i $ as a decoration on the decorated tree. 

The algebraic symmetry factor $S(\tau)$ of a decorated tree $\tau= X^{k} \zeta_l \ \prod_{j=1}^m \mathcal{I}_{a_j}(\tau_j)^{\beta_j}$ is defined grouping terms uniquely in such a way that $(a_i,\tau_i) \neq (a_j,\tau_j)$ for $i \neq j$, and setting inductively
\begin{equation*}
S(\tau) = k!\,
\bigg(
\prod_{j=1}^{m}
S(\tau_{j})^{\beta_{j}}
\beta_{j}!
\bigg),
\end{equation*} 
 starting with $S(\zeta_l)=1$ for $0\leq l\leq l_0$. We define an inner product on the set of all decorated trees setting for all $\sigma, \tau$
\begin{align*} \label{inner_product}
\langle \sigma ,  \tau \rangle \defeq S(\tau)\,\textbf{\textsf{1}}_{\sigma=\tau}.
\end{align*}
We also set 
$$
\langle \sigma_1\otimes\sigma_2 ,  \tau_1\otimes\tau_2 \rangle \defeq \langle \sigma_1 , \tau_1 \rangle\,\langle \sigma_2 , \tau_2 \rangle.
$$
The linear span of the set of decorated trees will be denoted by $T$. A {\it planted tree} is a tree of the form $\CI_a(\sigma)$, for a decorated tree $\sigma$ and $a\in\mathfrak{T}_+\times\mathbb{N}^{d+1}$; we denote by $\CI(T)$ the linear span of planted trees.

\smallskip

We now associate some numbers to every decorated tree. Fix a scaling vector $ \frak{s} \in \N^{d+1} $ and a map 
$$
|\cdot|_{\s} : \mathfrak{T}  \rightarrow \mathbb{R}
$$ 
which is negative on the noise types $ \mathfrak{T}_- $ and positive on the operator types $ \mathfrak{T}_+ $. This map accounts for the regularity of the noises and the gain of regularity of some operators encoded in the Schauder-type estimates that they satisfy. We extend the map $|\cdot|_{\s}$ to $\mathfrak{T}\times\N^{d+1} $ setting 
$$
 |q|_{\s} \defeq \sum_{i=1}^{d+1} \s_i  q_i, \qquad \textrm{and}\qquad |(\frak{t}, q)|_{\s} \defeq |\frak{t}|_{\s}  - |q|_{\s}.
$$ 
The degree of a decorated rooted tree $(\tau, \Labn,\Labe)$  is defined by
\begin{equation}
 \textrm{deg}(\tau, \Labn,\Labe) \defeq \sum_{v \in N_{T}} \big|\Labn(v)\big| _{\s} + \sum_{e \in E_{T}}  \big|(\frak{t}(e),\frak{q}(e))\big| _{\s} .
\end{equation}
 For cross reference, note that the `Degree' is called `homogeneity' in Hairer's work \cite{Hai14}. We use the degree to introduce the space of `positive' decorated trees $T^+$. It is the linear span of trees of the form $X^k \prod_{i=1}^{n} \CI_{a_i}(\tau_i)$, where $\textrm{deg}( \CI_{a_i}(\tau_i)) > 0$ and $k\in\N^{d+1}$. We also consider the linear space $T^-$ spanned by the decorated trees with negative degree, and denote by $\mathbb{R}[ T^-]$ the linear space spanned by forests of elements of $T^-$.

\smallskip

Given $k\in\N^{d+1}$ denote by $\uparrow_v^k$ the derivation on decorated trees that adds $k$ to the decoration at the node $v$, so $\uparrow_v^k(\sigma\tau) = \uparrow_v^k(\sigma)\tau + \sigma\uparrow_v^k(\tau)$ for all $\sigma,\tau\in T$. For $a=(\frak{t},p)$, one writes $a-m$ for $(\frak{t},p-m)$. We introduce a family of pre-Lie products of grafting type setting for all decorated trees $\sigma, \tau\in T$, with $\tau=(\tau,\frak{n},\frak{e})$ and $a\in \mathfrak{T}_+  \times\N^{d+1}$,
\begin{equation*}
\sigma \curvearrowright_a \tau \defeq \sum_{v\in N_{\tau}}\sum_{m\in\N^{d+1}} {\Labn_v \choose m} \,\sigma  \curvearrowright^v_{a-m}(\uparrow_v^{-m} \tau),
\end{equation*}
where $\Labn_v=((\Labn_v)_1,\cdots,(\Labn_v)_{d+1})\in\N^{d+1}$ is the decoration at the node $v$, the binomial coefficient ${\Labn_v \choose m}$ equals $\prod_{1\leq i\leq d+1}{(\Labn_v)_i \choose m_i}$, and $\curvearrowright^v_{a-m}$ grafts $ \sigma $ onto $\tau$ at the node $ v $ with an edge decorated by $a -m$. The above sum is finite due to the binomial coefficient $ {\Labn_v \choose m} $ which is equal to zero if $m$ is greater than $ \Labn_v $, by convention.  It is forbidden to graft on the top of noise edges which are terminal edges. Here is an example of this grafting product
\begin{equation*}  
	\begin{tikzpicture}[scale=0.2,baseline=0.1cm]
		\node at (0,0)  [dot,label= {[label distance=-0.2em]below: \scriptsize  $  \alpha   $} ] (root) {};
		\node at (0,4)  [dot,label={[label distance=-0.2em]above: \scriptsize  $  $}] (center) {};
		\draw[kernel1] (center) to
		node [sloped,below] {\small }     (root); 
		\node at (0,2) [fill=white,label={[label distance=0em]center: \scriptsize  $ \zeta_1 $} ] () {};
	\end{tikzpicture}
	 \curvearrowright^a   
	 \begin{tikzpicture}[scale=0.2,baseline=0.1cm]
	 	\node at (0,0)  [dot,label= {[label distance=-0.2em]below: \scriptsize  $  \beta  $} ] (root) {};
	 	\node at (0,8)  [dot,label= {[label distance=-0.2em]above: \scriptsize  $     $} ] (center) {};
	 	\node at (2,4)  [dot,label={[label distance=-0.2em]right: \scriptsize  $ \gamma $}] (right) {};
	 	\node at (-2,4)  [dot,label={[label distance=-0.2em]above: \scriptsize  $ $} ] (left) {};
	 	\draw[kernel1] (right) to
	 	node [sloped,below] {\small }     (root);
	 	\draw[kernel1] (center) to
	 	node [sloped,below] {\small }     (right); \draw[kernel1] (left) to
	 	node [sloped,below] {\small }     (root);
	 	\node at (-1,2) [fill=white,label={[label distance=0em]center: \scriptsize  $ \zeta_2 $} ] () {};
	 	\node at (1,2) [fill=white,label={[label distance=0em]center: \scriptsize  $ b $} ] () {};
	 	\node at (1,6) [fill=white,label={[label distance=0em]center: \scriptsize  $ \zeta_3 $} ] () {};
	 \end{tikzpicture}
	 = \sum_{k \in \mathbb{N}^{d+1} }{\beta \choose k} \begin{tikzpicture}[scale=0.2,baseline=0.1cm]
		\node at (0,0)  [dot,label= {[label distance=-0.2em]below: \scriptsize  $  \beta - k   $} ] (root) {};
		\node at (0,5)  [dot,label= {[label distance=-0.2em]right: \scriptsize  $  \gamma   $} ] (center) {};
		\node at (0,9)  [dot,label= {[label distance=-0.2em]above: \scriptsize  $    $} ] (centerc) {};
		\node at (3,4)  [dot,label={[label distance=-0.2em]right: \scriptsize  $ \alpha $}] (right) {};
		\node at (-3,4)  [dot,label={[label distance=-0.2em]above: \scriptsize  $ $} ] (left) {};
		\node at (3,8)  [dot,label={[label distance=-0.2em]above: \scriptsize  $ $} ] (rightc) {};
		\draw[kernel1] (right) to
		node [sloped,below] {\small }     (root); 
		\draw[kernel1] (center) to
		node [sloped,below] {\small }     (root);
		\draw[kernel1] (center) to
		node [sloped,below] {\small }     (centerc);
		\draw[kernel1] (right) to
		node [sloped,below] {\small }     (rightc); 
		\draw[kernel1] (left) to
		node [sloped,below] {\small }     (root);
		\node at (3,6) [fill=white,label={[label distance=0em]center: \scriptsize  $ \zeta_1 $} ] () {};
		\node at (0,7) [fill=white,label={[label distance=0em]center: \scriptsize  $ \zeta_3 $} ] () {};
		\node at (-2,2) [fill=white,label={[label distance=0em]center: \scriptsize  $ \zeta_2 $} ] () {};
			\node at (2,2) [fill=white,label={[label distance=0em]center: \scriptsize  \qquad $ a - k $} ] () {};
		\node at (0,2.5) [fill=white,label={[label distance=0em]center: \scriptsize  $ b  $} ] () {};
	\end{tikzpicture} + \sum_{k \in \mathbb{N}^{d+1} }{\gamma \choose k} \begin{tikzpicture}[scale=0.2,baseline=0.1cm]
		\node at (0,0)  [dot,label= {[label distance=-0.2em]below: \scriptsize  $  \beta  $} ] (root) {};
		\node at (0,8)  [dot,label= {[label distance=-0.2em]above: \scriptsize  $    $} ] (center) {};
		\node at (2,4)  [dot,label={[label distance=-0.2em]right: \scriptsize  $ \gamma - k $}] (right) {};
		\node at (4,8)  [dot,label={[label distance=-0.2em]right: \scriptsize  $ \alpha $}] (rightr) {};
		\node at (2,12)  [dot,label={[label distance=-0.2em]right: \scriptsize  $ \ $}] (rightrl) {};
		\node at (-2,4)  [dot,label={[label distance=-0.2em]above: \scriptsize  $ $} ] (left) {};
		\draw[kernel1] (rightr) to
		node [sloped,below] {\small }     (rightrl);
		\draw[kernel1] (right) to
		node [sloped,below] {\small }     (root);
		\draw[kernel1] (right) to
		node [sloped,below] {\small }     (rightr);
		\draw[kernel1] (center) to
		node [sloped,below] {\small }     (right); \draw[kernel1] (left) to
		node [sloped,below] {\small }     (root);
		\node at (-1,2) [fill=white,label={[label distance=0em]center: \scriptsize  $ \zeta_2 $} ] () {};
		\node at (1,2) [fill=white,label={[label distance=0em]center: \scriptsize  $ b $} ] () {};
		\node at (3,6) [fill=white,label={[label distance=0em]center: \scriptsize \qquad $ a -k $} ] () {};
		\node at (3,10) [fill=white,label={[label distance=0em]center: \scriptsize  $ \zeta_1 $} ] () {};
		\node at (1,6) [fill=white,label={[label distance=0em]center: \scriptsize  $ \zeta_3 $} ] () {};
	\end{tikzpicture}. 
\end{equation*}
The pre-Lie products $\curvearrowright_a$ are non-commutative; they were first introduced in Bruned, Chandra, Chevyrev and Hairer's work \cite{BCCH18}. We recall one universal result that we will use in the sequel; it was first established in Corollary 4.23 of \cite{BCCH18}. It can be viewed as an extension of the universal result of Chapoton-Livernet \cite{ChaLiv} on pre-Lie algebras. (Such a result becomes immediate when one constructs $\curvearrowright_a$  as a deformation, as in Section 2.1 of \cite{BrunedManchon}. See also Foissy's work \cite{F2018} for the case with no deformation.)


\begin{proposition} \label{freeness} 
The space $T$ is freely generated by the elements 
\begin{equation*}
\Big\{ X^k \zeta_l; \, 1\leq l\leq l_0, \ k \in \N^{d+1}\Big\}
\end{equation*}
 and the operations $\big\{\hspace{-0.08cm}\curvearrowright_a\,;\,a\in \frak{T}_+\times \bbN^{d+1}\big\}$.
\end{proposition}


 Let us mention that we will not use in the sequel the pre-Lie character of the maps $\curvearrowright_a$; this property will only play an implicit role via the use of Proposition \ref{freeness} in the proof of Proposition \ref{star_morphism}. We define a product 
$$
\curvearrowright\,: \CI(T)\times T\rightarrow T
$$ 
setting for all $a\in \frak{T}_+\times\N^{d+1}, \sigma, \tau\in T$
\begin{equation} \label{def_pre_Lie}
\CI_a(\sigma ) \, \curvearrowright \, \tau \defeq \sigma \, \curvearrowright_a \, \tau.
\end{equation}
 We also set $\one\curvearrowright\tau \defeq \tau$.  We invite the reader to check that if we graft $n$ trees $\sigma_i$ along some grafting operators $\curvearrowright_{a_i}$ on $\tau$, {\it not on one another}, then the resulting element of $T$ does not depend on the order in which we grafted our trees. We can then extend the $\curvearrowright$ to products of planted trees $\prod_{i=1}^n \CI_{a_i}(\sigma_i)$, so
\begin{equation*}
	\prod_{i=1}^n \CI_{a_i}(\sigma_i) \curvearrowright\tau
\end{equation*}
means that we perform the silmutaneous grafting of the decorated trees $ \sigma_i $ onto $ \sigma $ along some edges decorated by the $ a_i $. Some nodes of $ \tau $ can received several grafting and we sum over all the possibilities. The order of the grafting does not matter. Below, we provide an example of computation
\begin{equation*}
 \begin{tikzpicture}[scale=0.2,baseline=0.1cm]
	\node at (0,0)  [dot,label= {[label distance=-0.2em]below: \scriptsize  $    $} ] (root) {};
	\node at (2,8)  [dot,label= {[label distance=-0.2em]above: \scriptsize  $     $} ] (center) {};
	\node at (-2,8)  [dot,label= {[label distance=-0.2em]above: \scriptsize  $ \     $} ] (centerl) {};
	\node at (2,4)  [dot,label={[label distance=-0.2em]right: \scriptsize  $ \beta $}] (right) {};
	\node at (-2,4)  [dot,label={[label distance=-0.2em]right: \scriptsize  $ \alpha $} ] (left) {};
	\draw[kernel1] (left) to
	node [sloped,below] {\small }     (centerl);
	\draw[kernel1] (right) to
	node [sloped,below] {\small }     (root);
	\draw[kernel1] (center) to
	node [sloped,below] {\small }     (right); \draw[kernel1] (left) to
	node [sloped,below] {\small }     (root);
	\node at (-1,2) [fill=white,label={[label distance=0em]center: \scriptsize  $ a $} ] () {};
	\node at (-2,6) [fill=white,label={[label distance=0em]center: \scriptsize  $ \zeta_2 $} ] () {};
	\node at (1,2) [fill=white,label={[label distance=0em]center: \scriptsize  $ b $} ] () {};
	\node at (2,6) [fill=white,label={[label distance=0em]center: \scriptsize  $ \zeta_3 $} ] () {};
\end{tikzpicture}		\curvearrowright 	\begin{tikzpicture}[scale=0.2,baseline=0.1cm]
		\node at (0,0)  [dot,label= {[label distance=-0.2em]below: \scriptsize  $  k   $} ] (root) {};
		\node at (0,4)  [dot,label={[label distance=-0.2em]above: \scriptsize  $  $}] (center) {};
		\draw[kernel1] (center) to
		node [sloped,below] {\small }     (root); 
		\node at (0,2) [fill=white,label={[label distance=0em]center: \scriptsize  $ \zeta_1 $} ] () {};
	\end{tikzpicture}
	= \sum_{k_1, k_2 \in \mathbb{N}^{d+1} }{\gamma \choose k_1, k_2} \begin{tikzpicture}[scale=0.2,baseline=0.1cm]
		\node at (0,0)  [dot,label= {[label distance=-0.2em]below: \scriptsize  $  k - k_1 - k_2   $} ] (root) {};
		\node at (0,5)  [dot,label= {[label distance=-0.2em]right: \scriptsize  $  \alpha   $} ] (center) {};
		\node at (0,9)  [dot,label= {[label distance=-0.2em]above: \scriptsize  $    $} ] (centerc) {};
		\node at (3,4)  [dot,label={[label distance=-0.2em]right: \scriptsize  $ \beta $}] (right) {};
		\node at (-3,4)  [dot,label={[label distance=-0.2em]above: \scriptsize  $ $} ] (left) {};
		\node at (3,8)  [dot,label={[label distance=-0.2em]above: \scriptsize  $ $} ] (rightc) {};
		\draw[kernel1] (right) to
		node [sloped,below] {\small }     (root); 
		\draw[kernel1] (center) to
		node [sloped,below] {\small }     (root);
		\draw[kernel1] (center) to
		node [sloped,below] {\small }     (centerc);
		\draw[kernel1] (right) to
		node [sloped,below] {\small }     (rightc); 
		\draw[kernel1] (left) to
		node [sloped,below] {\small }     (root);
		\node at (3,6) [fill=white,label={[label distance=0em]center: \scriptsize  $ \zeta_3 $} ] () {};
		\node at (0,7) [fill=white,label={[label distance=0em]center: \scriptsize  $ \zeta_2 $} ] () {};
		\node at (-2,2) [fill=white,label={[label distance=0em]center: \scriptsize  $ \zeta_1$} ] () {};
		\node at (2,2) [fill=white,label={[label distance=0em]center: \scriptsize  \qquad $ b_k $} ] () {};
		\node at (0,2.5) [fill=white,label={[label distance=0em]center: \scriptsize  $a_k$} ] () {};
	\end{tikzpicture} .
\end{equation*}
where $ a_k = a - k_1 $ and $ b_k = b - k_2 $. In this example only one node is allowed for the grafting which is the root decorated by $k$. One can observe the symmetry between $k_1$ and $k_2$ that illustrate that the order of the grafting does not matter.

Following Bruned and Manchon's construction in \cite{BrunedManchon}, for  $ B \subset N_{\tau} $, consider the derivation map $ \uparrow^{k}_{B}$ defined as
\begin{equation}
\uparrow^{k}_{B} \tau = \sum_{\sum_{v \in B} k_v  = k } \prod_{v \in B}\uparrow_v^{k_v}  \tau.
\end{equation}
 In the sequel, we will use the following short hand notation: $\uparrow^k_{N_\tau} \tau = \uparrow^k \tau$. Below, we provide an example of computation
	\begin{equation*}
	\uparrow^k	 \begin{tikzpicture}[scale=0.2,baseline=0.1cm]
			\node at (0,0)  [dot,label= {[label distance=-0.2em]below: \scriptsize  $  \beta  $} ] (root) {};
			\node at (0,8)  [dot,label= {[label distance=-0.2em]above: \scriptsize  $     $} ] (center) {};
			\node at (2,4)  [dot,label={[label distance=-0.2em]right: \scriptsize  $ \gamma $}] (right) {};
			\node at (-2,4)  [dot,label={[label distance=-0.2em]above: \scriptsize  $ $} ] (left) {};
			\draw[kernel1] (right) to
			node [sloped,below] {\small }     (root);
			\draw[kernel1] (center) to
			node [sloped,below] {\small }     (right); \draw[kernel1] (left) to
			node [sloped,below] {\small }     (root);
			\node at (-1,2) [fill=white,label={[label distance=0em]center: \scriptsize  $ \zeta_2 $} ] () {};
			\node at (1,2) [fill=white,label={[label distance=0em]center: \scriptsize  $ b $} ] () {};
			\node at (1,6) [fill=white,label={[label distance=0em]center: \scriptsize  $ \zeta_3 $} ] () {};
		\end{tikzpicture}  = \sum_{k = k_1 + k_2} {k \choose k_1, k_2}
	\begin{tikzpicture}[scale=0.2,baseline=0.1cm]
		\node at (0,0)  [dot,label= {[label distance=-0.2em]below: \scriptsize  $  \beta + k_1  $} ] (root) {};
		\node at (0,8)  [dot,label= {[label distance=-0.2em]above: \scriptsize  $     $} ] (center) {};
		\node at (2,4)  [dot,label={[label distance=-0.2em]right: \scriptsize  $ \gamma + k_2 $}] (right) {};
		\node at (-2,4)  [dot,label={[label distance=-0.2em]above: \scriptsize  $ $} ] (left) {};
		\draw[kernel1] (right) to
		node [sloped,below] {\small }     (root);
		\draw[kernel1] (center) to
		node [sloped,below] {\small }     (right); \draw[kernel1] (left) to
		node [sloped,below] {\small }     (root);
		\node at (-1,2) [fill=white,label={[label distance=0em]center: \scriptsize  $ \zeta_2 $} ] () {};
		\node at (1,2) [fill=white,label={[label distance=0em]center: \scriptsize  $ b $} ] () {};
		\node at (1,6) [fill=white,label={[label distance=0em]center: \scriptsize  $ \zeta_3 $} ] () {};
	\end{tikzpicture}.
	\end{equation*}   
We define the product
$$
\star : T\times T\rightarrow T
$$
for all $\sigma = X^k \prod_{i} \CI_{a_i}(\sigma_i) \in T$ and $\tau \in T$ by the formula
\begin{equation}
\sigma \star \tau \defeq  \, \uparrow^k_{N_\tau} \left(  \Big(\prod_i \mathcal{I}_{a_i}(\sigma_i)\Big)\curvearrowright \tau \right).
\end{equation}
 Let us emphasize that the operator $\uparrow^k_{N_\tau}$ only acts on the nodes of $\tau$.  One has for instance
\begin{equation} \label{EqExampleStarProduct}
X^k\,\zeta_l\prod_{i=1}^n\CI_{a_i}(\tau_i) = \left(X^k\prod_{i=1}^n\CI_{a_i}(\tau_i)\right)\star \zeta_l,
\end{equation}
 and 
$$
\Big\{X^{k} \prod_{j=1}^{n} \CI_{a_j}(\tau_j)\Big\} \star (X^m \zeta_l) = \zeta_l \hspace{-0.2cm} \sum_{r_1,\dots,r_n \in \N^{d+1}}  X^{m -  \sum_j r_j + k} \prod_{j=1}^n {m \choose r_j} \CI_{a_j-r_i}(\tau_j).
$$
Note that our product $\star$ corresponds to the $\star_2$ product in \cite{BrunedManchon}. 
One has also
	\begin{equation*}
		\begin{aligned}
	\begin{tikzpicture}[scale=0.2,baseline=0.1cm]
	\node at (0,0)  [dot,label= {[label distance=-0.2em]below: \scriptsize  $  k  $} ] (root) {};
	\node at (0,8)  [dot,label={[label distance=-0.2em]above: \scriptsize  $  $}] (centerc) {};
		\node at (0,4)  [dot,label={[label distance=-0.2em]right: \scriptsize  $ \alpha  $}] (center) {};
	\draw[kernel1] (center) to
	node [sloped,below] {\small }     (root);
	\draw[kernel1] (center) to
	node [sloped,below] {\small }     (centerc); 
	\node at (0,6) [fill=white,label={[label distance=0em]center: \scriptsize  $ \zeta_1 $} ] () {};
		\node at (0,2) [fill=white,label={[label distance=0em]center: \scriptsize  $ a $} ] () {};
\end{tikzpicture}
\star   
\begin{tikzpicture}[scale=0.2,baseline=0.1cm]
	\node at (0,0)  [dot,label= {[label distance=-0.2em]below: \scriptsize  $  \beta  $} ] (root) {};
	\node at (0,8)  [dot,label= {[label distance=-0.2em]above: \scriptsize  $     $} ] (center) {};
	\node at (2,4)  [dot,label={[label distance=-0.2em]right: \scriptsize  $ \gamma $}] (right) {};
	\node at (-2,4)  [dot,label={[label distance=-0.2em]above: \scriptsize  $ $} ] (left) {};
	\draw[kernel1] (right) to
	node [sloped,below] {\small }     (root);
	\draw[kernel1] (center) to
	node [sloped,below] {\small }     (right); \draw[kernel1] (left) to
	node [sloped,below] {\small }     (root);
	\node at (-1,2) [fill=white,label={[label distance=0em]center: \scriptsize  $ \zeta_2 $} ] () {};
	\node at (1,2) [fill=white,label={[label distance=0em]center: \scriptsize  $ b $} ] () {};
	\node at (1,6) [fill=white,label={[label distance=0em]center: \scriptsize  $ \zeta_3 $} ] () {};
\end{tikzpicture}
& = \sum_{k = k_1+k_2} {k \choose k_1,k_2} \bigg(  \sum_{m \in \mathbb{N}^{d+1} }{\beta \choose m} \begin{tikzpicture}[scale=0.2,baseline=0.1cm]
	\node at (0,0)  [dot,label= {[label distance=-0.2em]below: \scriptsize  $  \beta +k_1 - m   $} ] (root) {};
	\node at (0,5)  [dot,label= {[label distance=-0.2em]left: \scriptsize  $  \gamma + k_2   $} ] (center) {};
	\node at (0,9)  [dot,label= {[label distance=-0.2em]above: \scriptsize  $    $} ] (centerc) {};
	\node at (3,4)  [dot,label={[label distance=-0.2em]right: \scriptsize  $ \alpha $}] (right) {};
	\node at (-3,4)  [dot,label={[label distance=-0.2em]above: \scriptsize  $ $} ] (left) {};
	\node at (3,8)  [dot,label={[label distance=-0.2em]above: \scriptsize  $ $} ] (rightc) {};
	\draw[kernel1] (right) to
	node [sloped,below] {\small }     (root); 
	\draw[kernel1] (center) to
	node [sloped,below] {\small }     (root);
	\draw[kernel1] (center) to
	node [sloped,below] {\small }     (centerc);
	\draw[kernel1] (right) to
	node [sloped,below] {\small }     (rightc); 
	\draw[kernel1] (left) to
	node [sloped,below] {\small }     (root);
	\node at (3,6) [fill=white,label={[label distance=0em]center: \scriptsize  $ \zeta_1 $} ] () {};
	\node at (0,7) [fill=white,label={[label distance=0em]center: \scriptsize  $ \zeta_3 $} ] () {};
	\node at (-2,2) [fill=white,label={[label distance=0em]center: \scriptsize  $ \zeta_2 $} ] () {};
	\node at (2,2) [fill=white,label={[label distance=0em]center: \scriptsize  \qquad $ a - m $} ] () {};
	\node at (0,2.5) [fill=white,label={[label distance=0em]center: \scriptsize  $ b  $} ] () {};
\end{tikzpicture} \\ & + \sum_{m \in \mathbb{N}^{d+1} }{\gamma \choose m} \begin{tikzpicture}[scale=0.2,baseline=0.1cm]
	\node at (0,0)  [dot,label= {[label distance=-0.2em]below: \scriptsize  $  \beta + k_1  $} ] (root) {};
	\node at (0,8)  [dot,label= {[label distance=-0.2em]above: \scriptsize  $    $} ] (center) {};
	\node at (2,4)  [dot,label={[label distance=-0.2em]right: \scriptsize  $ \gamma + k_2 - m $}] (right) {};
	\node at (4,8)  [dot,label={[label distance=-0.2em]right: \scriptsize  $ \alpha $}] (rightr) {};
	\node at (2,12)  [dot,label={[label distance=-0.2em]right: \scriptsize  $ \ $}] (rightrl) {};
	\node at (-2,4)  [dot,label={[label distance=-0.2em]above: \scriptsize  $ $} ] (left) {};
	\draw[kernel1] (rightr) to
	node [sloped,below] {\small }     (rightrl);
	\draw[kernel1] (right) to
	node [sloped,below] {\small }     (root);
	\draw[kernel1] (right) to
	node [sloped,below] {\small }     (rightr);
	\draw[kernel1] (center) to
	node [sloped,below] {\small }     (right); \draw[kernel1] (left) to
	node [sloped,below] {\small }     (root);
	\node at (-1,2) [fill=white,label={[label distance=0em]center: \scriptsize  $ \zeta_2 $} ] () {};
	\node at (1,2) [fill=white,label={[label distance=0em]center: \scriptsize  $ b $} ] () {};
	\node at (3,6) [fill=white,label={[label distance=0em]center: \scriptsize \qquad $ a - m $} ] () {};
	\node at (3,10) [fill=white,label={[label distance=0em]center: \scriptsize  $ \zeta_1 $} ] () {};
	\node at (1,6) [fill=white,label={[label distance=0em]center: \scriptsize  $ \zeta_3 $} ] () {};
\end{tikzpicture} \bigg).
\end{aligned}
\end{equation*}

It has been proved in Section 3.3 of \cite{BrunedManchon} that this product is associative; this can be obtained by applying the Guin-Oudom procedure \cite{Guin1,Guin2} to a well-chosen pre-Lie product. When $\sigma \in T^+$ and $\tau, \mu\in T$ one has from Theorem 4.2 in \cite{BrunedManchon}
\begin{equation} \label{def_coaction}
\langle  \sigma \star \tau ,  \mu \rangle  = \langle  \tau \otimes \sigma ,  \Delta \mu \rangle 
\end{equation}
where 
$$
\Delta : T \rightarrow T \otimes T^+ 
$$ 
is a co-action first introduced in Hairer's seminal work \cite{Hai14} -- see also \cite{BHZ} and \cite{RSGuide} where it plays a prominent role. So the restriction to $T^+\times T$ of the product $\star$ is the $\langle\cdot,\cdot\rangle$-dual of the splitting map $\Delta$. The restriction to $ T^-\times T$ of the product $\star$ is also the dual of the splitting map $\delta_r$ introduced by Bruned in Section 4 of \cite{BrunedRecursive}, in formula (23) and Definition 4.2 therein, in the sense that one has for all $(\tau,\sigma)\in T^-\times T$ and $\mu\in T$ the identity
$$
\langle  \sigma \star \tau ,  \mu \rangle  = \langle  \tau \otimes \sigma ,  \delta_r \mu \rangle.
$$
 One can take this relation as a definition of the map $\delta_r$. The index $r$ stands here for the word `root'. Note here as in \eqref{def_coaction} the different orders of $\sigma$ and $\tau$ in the left and right hand sides of the equalities.

\smallskip

 Recall $(\epsilon_i)_{1\leq i\leq d+1}$ stands for the canonical basis of $\N^{d+1}$. With a view on the system \eqref{EqSPDE} assume we are given a family 
$$
(\Upsilon_p^l)_{1\leq p\leq p_0, 0\leq l\leq l_0}
$$ 
of functions of a finite number of abstract variables
\begin{equation} \label{variables}
	(Z_a)_{a\in( \frak{T}_+ \cup\{*\})\times\N^{d+1}}.
	\end{equation}
   We impose that the $Z_{(*,\epsilon_i)}$ are all equal to $1$ and $Z_{(*,k)}=0$ for $k\in\N^{d+1}$ not in the canonical basis. In the particular setting of the system \eqref{EqSPDE} the functions $\Upsilon_p^l((Z_a)_a)$ are place holders for the functions $f_p^l$, when $1\leq l\leq l_0$, and the function $g_p$ when $l=0$. The variables $Z_a$ play the roles of $x$ and $u$ and its derivatives, with $Z_*$ in the role of $x$ and $Z_{\frak{t}^p_+}$ in the role of $u_p$, and $Z_{(\frak{t}^p_+,q)}$ in the role of $\partial^q u_p$, for instance. This interpretation explains our choice for the variables $Z_{(*,k)}$. We define the partial derivative operators 
$$
D_a = D_{Z_a}
$$ 
with respect to the variable $Z_a$ in the usual way and set
$$ 
\partial^{ \epsilon_i} \defeq \sum_a Z_{a+ \epsilon_i}\,D_a.
$$
 where $a + \epsilon_i$ is defined by $ (\mathfrak{t}, p +  \epsilon_i)$ for $a = (\mathfrak{t}, p)$. These operators commute since
\begin{equation*} \begin{split}
\partial^{\epsilon_i}\partial^{\epsilon_j} &= \sum_{a_1,a_2}Z_{a_1+\epsilon_1}D_{a_1}\big(Z_{a_2+\epsilon_2}D_{a_2}\big)   \\
&= \sum_{a_1,a_2} \Big(Z_{a_1+\epsilon_1}\delta_{a_1}(a_2+\epsilon_2) D_{a_2} + Z_{a_1+\epsilon_1}Z_{a_2+\epsilon_2}D_{a_1}D_{a_2}\Big)   \\
&= \sum_a Z_{a+\epsilon_1+\epsilon_2}D_a + \sum_{a_1,a_2} Z_{a_1+\epsilon_1}Z_{a_2+\epsilon_2}D_{a_1}D_{a_2}
\end{split} \end{equation*}
is symmetric in $\epsilon_i, \epsilon_j$. So one can define for $k=\big(k_1,\dots,k_{d+1}\big)\in\N^{d+1}$
\begin{equation*}
\partial^k \defeq \prod_{i=1}^{d+1} (\partial^{\epsilon_i})^{k_i}.
\end{equation*}
These operators will act on functions that depend on finitely many variables $Z_a$ so their action will not involve infinite sums, unlike what the definition of $\partial^{\epsilon_i}$ may suggest. 

We define inductively a family $\Upsilon = (\Upsilon_p)_{1\leq p\leq p_0}$ of  functions of the variables $(Z_a)_a$, indexed by $T$, setting for $ \tau = X^{k} \zeta_{l}\prod_{j=1}^{n} \CI_{a_j}(\tau_j) $, with $a_j = (\Labhom_{l_j},k_j)$ and $1\leq l_j\leq p_0$
\begin{equation} \label{def_upsilon} \begin{split}
\Upsilon_p(\zeta_l) &\defeq \Upsilon^l_p    \\
\Upsilon_p(\tau) &\defeq \Big\{\partial^k D_{a_1} ... D_{a_n} \Upsilon_p(\zeta_l)\Big\}\,\prod_{j=1}^n \Upsilon_{l_j}(\tau_j)
\end{split} \end{equation}
 for all $1\leq p\leq p_0$. The maps $ \Upsilon_p $ are extended to $T$ by linearity. The following fact plays a key role in the proof of Proposition \ref{star_morphism}.

\begin{lemma} \label{LemBailleulHoshino}
One has for all $m\in\N^{d+1}$ and $a=(\frak{t},q)\in(\frak{T}_+ \cup\{*\})\times\N^{d+1}$ the identity
\begin{equation} \label{commutation derivatives}
\sum_{r_1,\dots,r_n \in \N^{d+1}} \partial^{m-\sum_{j_1} r_{j_1}} \prod_{j=1}^n  {m \choose r_j}  D_{a_j - r_j} = \Big\{\prod_{j=1}^n D_{a_j} \Big\} \partial^{m}.
\end{equation}
\end{lemma}

\begin{proof}
 One can check that for an element $\epsilon_i$ of the canonical basis of $ \mathbb{N}^{d+1} $, one has
\begin{equation*}
	\sum_{r_1,\dots,r_n \in \N^{d+1}} \partial^{\epsilon_i - \sum_{j_1} r_{j_1}} \prod_{j=1}^n  {\epsilon_i \choose r_j}  D_{a_j - r_j} = \Big\{\prod_{j=1}^n D_{a_j} \Big\} \partial^{\epsilon_i}.
\end{equation*}
This is a consequence of the fact that
\begin{equation*}
	D_{a_j} \partial^{\epsilon_i} = \partial_{\epsilon_i} D_{a_j} + D_{a_j - \epsilon_i}
\end{equation*}
coming from the definitions of $ \partial^{\epsilon_i} $ and $ D_{a_j} $. We proceed by recurrence and we suppose true for $m \in \mathbb{N}^{d+1}$. One has 
\begin{equation*}
	\begin{aligned}
	\Big\{\prod_{j=1}^n D_{a_j} \Big\} \partial^{m + \epsilon_i} & = \sum_{r_1,\dots,r_n \in \N^{d+1}} \partial^{m-\sum_{j_1} r_{j_1}} \prod_{j=1}^n  {m \choose r_j}  D_{a_j - r_j} \partial^{\epsilon_i}
	\\  & = \sum_{r_1,\bar{r}_1,\dots,r_n, \bar{r}_n \in \N^{d+1}} \partial^{m-\sum_{j_1} r_{j_1}} \prod_{j=1}^n  {m \choose r_j} {\epsilon_i \choose \bar{r}_j}  D_{a_j - r_j-\bar{r}_{j}} 
	\\ & = \sum_{r_1,\dots,r_n \in \N^{d+1}} \partial^{m-\sum_{j_1} r_{j_1}} \prod_{j=1}^n  {m + \epsilon_i \choose r_j}  D_{a_j - r_j}
	\end{aligned}
\end{equation*}
where we have used the Chu-Vandermonde identity that tells us that
\begin{equation*}
	\sum_{r_j+\bar{r}_j = k_j}  {m \choose r_j} {\epsilon_i \choose \bar{r}_j} =   {m + \epsilon_i \choose k_j}.
	\end{equation*}
\end{proof}

The next statement is a morphism property of the map $(\tau\in T)\mapsto \Upsilon(\tau)$ for the $\star$ product.


\begin{proposition} \label{star_morphism}
For every $1\leq p\leq p_0$, for every $ \sigma\in T$, and $\tau = X^{k} \prod_{j=1}^{n} \CI_{a_j}(\tau_j) \in T$ with $a_j = (\Labhom_{l_j}, q_j)$ and $1\leq l_k\leq p_0$, one has
\begin{equation}\label{star2_morphism}
\Upsilon_p\bigg(\Big\{X^{k} \prod_{j=1}^{n} \CI_{a_j}(\tau_j)\Big\} \star \sigma\bigg) = \Big\{\partial^k D_{a_1} ... D_{a_n} \Upsilon_p(\sigma)\Big\}\,\prod_{j=1}^n \Upsilon_{l_j}(\tau_j).
\end{equation}
\end{proposition}


We see from identity \eqref{EqExampleStarProduct} that the defining inductive property \eqref{def_upsilon} is a particular case of \eqref{star2_morphism}. In the simple case where $\tau = \mcI_{b_1}(\sigma_1)$ with $b_1 = (\Labhom_b,q_b)$ formula \eqref{star2_morphism} reads
$$
 \Upsilon_p\big(\mcI_{b_1}(\sigma_1)\star\sigma\big) = D_{b_1}\hspace{-0.1cm}\Upsilon_p(\sigma) \, \Upsilon_b(\sigma_1).
$$


\begin{proof}
We proceed by induction on $\vert\sigma\vert_\zeta$. It follows from Proposition \ref{freeness} that it suffices to prove the statement for $\sigma=X^m\zeta_l$ and $\sigma=\mcI_b(\sigma_1)\curvearrowright\sigma_2$ with $\sigma_1$ and $\sigma_2$ satisfying \eqref{star2_morphism} for all $\tau$ as above.

-- The case $ \sigma = \zeta_l $ is part of the definition \eqref{def_upsilon}. For $\sigma = X^m \zeta_l$, we have
$$
\Big\{X^{k} \prod_{j=1}^{n} \CI_{a_j}(\tau_j)\Big\} \star (X^m \zeta_l) = \zeta_l \hspace{-0.2cm} \sum_{r_1,\dots,r_n \in \N^{d+1}}  X^{k+m -  \sum_{j_1} r_{j_1}} \prod_{j=1}^n {m \choose r_j} \CI_{a_j - r_j}(\tau_j).
$$
Using Lemma \ref{LemBailleulHoshino} and the linearity of $ \Upsilon_p $ one has
\begin{equation*} \begin{split}
\Upsilon_p\bigg(\Big\{X^{k} \prod_{j=1}^{n} &\CI_{a_j}(\tau_j)\Big\} \star (X^m \zeta_l)\bigg)   \\
& = \Big\{\prod_{j=1}^{n} \Upsilon_{l_j}(\tau_j)\Big\} \sum_{r_1,\dots, r_{n}     \in \N^{d+1}} \partial^{k+m-\sum_{j_1} r_{j_1}} \prod_{j=1}^n  {m \choose r_j} D_{a_j - r_j} \Upsilon_p (\zeta_l)   \\
&= \Big\{\prod_{j=1}^{n} \Upsilon_{l_j}(\tau_j) \Big\} \, \Big\{\partial^k \prod_{j=1}^n D_{a_{j}} \partial^m \Upsilon_p(\zeta_l) \Big\}  \\
&= \Big\{\prod_{j=1}^{n} \Upsilon_{l_j}(\tau_j) \Big\} \, \Big\{\partial^k \prod_{j=1}^n D_{a_{j}} \Upsilon_p(X^{m}\zeta_l)\Big\},
\end{split} \end{equation*}
which is indeed \eqref{star2_morphism} in the particular case where $ \sigma = X^{m} \zeta_l $. 

-- We now assume that $ \sigma  = \CI_{b_1}(\sigma_1)\curvearrowright\sigma_2  = \CI_{b_1}(\sigma_1) \star \sigma_2 $ with $b_1 = (\Labhom_b,q_b)$  and $\sigma_1$ and $\sigma_2$ satisfying  \eqref{star2_morphism} for all $\tau$. Pick $\tau = X^{k} \prod_{j=1}^{n} \CI_{a_j}(\tau_j) \in T$ with $ a_j = (\Labhom_{l_j}, q_j)$. One has from the associativity of $\star$ the equality
\begin{equation*} \begin{split}
& \tau  \star \sigma  = \big( \tau \star \CI_{b_1}(\sigma_1) \big) \star \sigma_2
\\ & = \left\{ \sum_{I \subset \lbrace 1,...,n \rbrace}\sum_{k_1 + k_2 = k}  {k \choose k_1,k_2} X^{k_1} \prod_{i \in I} \CI_{a_i}(\tau_i) \, \CI_{b_1}\left( \bigg(X^{k_2} \prod_{j\in I^c}  \CI_{a_j}(\tau_j) \bigg) \star \sigma_1 \right) \right\} \star \sigma_2.
\end{split} \end{equation*}
 The idendity for $ \tau \star \CI_{b_1}(\sigma_1) $ comes from the fact that
\begin{equation*}
	\tau \star \CI_{b_1}(\sigma_1) = \uparrow^k_{N_\sigma} \left(  \Big(\prod_{j=1}^n \mathcal{I}_{a_j}(\tau_j)\Big)\curvearrowright \CI_{b_1}(\sigma_1) \right)
	\end{equation*}
Then, one grafts some of the decorated trees $ \tau_j$ at the root of $ \sigma$ via the term $ \prod_{i \in I} \CI_{a_i}(\tau_i) $.
The rest of the decorated trees are grafted on the other nodes given. This term is given by $  \bigg(\prod_{j\in I^c}  \CI_{a_j}(\tau_j) \bigg) \curvearrowright \sigma_1$. One is summing over all the possibilities of such splitting. We do the same reasoning for justifying the splitting of $ X^{k} $ into $ X^{k_1} $ and $ X^{k_2} $. One has to be careful to not forget the combinatorial factor $ {k \choose k_1,k_2} $.

 We apply the induction hypothesis to get for $\Upsilon_p(\tau \star \sigma)$ the formula
\begin{equation*} \begin{split}
&\sum_{I \subset \lbrace 1,...,n \rbrace}\sum_{k_1 + k_2 = k}  {k \choose k_1,k_2} \bigg\{\Big(\partial^{k_1}\hspace{-0.05cm}\prod_{i\in I} D_{a_i} D_{b_1}\Big) \Upsilon_p(\sigma_2)\bigg\} \prod_{i\in I}\Upsilon_{l_i} (\tau_i) \, \Upsilon_b\Big(\big(X^{k_2} \prod_{j\in I^c}  \CI_{a_j}(\tau_j) \big) \star \sigma_1\Big)   \\ 
&= \underset{k_1 + k_2 = k}{\sum_{I \subset \lbrace 1,...,n \rbrace}}  {k \choose k_1,k_2} \bigg\{\Big(\partial^{k_1}\hspace{-0.05cm}\prod_{i\in I} D_{a_i} D_{b_1}\Big) \Upsilon_p(\sigma_2)\bigg\} \prod_{i\in I}\Upsilon_{l_i} (\tau_i) \, \Big\{\Big(\partial^{k_2}\hspace{-0.1cm}\prod_{j\in I^c} D_{a_j}\Big) \Upsilon_b(\sigma_1)\Big\} \prod_{j\in I^c} \Upsilon_{l_j}(\tau_j)  \\
&=  \prod_{j =1}^{n} \Upsilon_{l_j}(\tau_j) \sum_{I \subset \lbrace 1,...,n \rbrace} \sum_{k_1 + k_2 = k}  {k \choose k_1,k_2} \Big\{\Big(\partial^{k_1} \hspace{-0.07cm}\prod_{i\in I} D_{a_i} D_{b_1}\Big) \Upsilon_p(\sigma_2)\Big\} \,  \Big\{\Big(\partial^{k_2} \hspace{-0.13cm}\prod_{j\in I^c} D_{a_j}\Big) \Upsilon_b(\sigma_1)\Big\}.
\end{split} \end{equation*}
Using the fact that $ \partial^{k} $ and $ D_{a_i} $ satisfy the Leibniz rule one then gets
\begin{equation*} \begin{split}
\Upsilon_p(\tau \star \sigma) & = \prod_{j =1}^{n} \Upsilon_{l_j}(\tau_j) \, \Big\{\partial^{k} \prod_{i=1}^n D_{a_i} \Big\} \Big(\Upsilon_b(\sigma_1) D_{b_1}\hspace{-0.1cm} \Upsilon_p(\sigma_2)\Big)   \\
 &= \prod_{j =1}^{n} \Upsilon_{l_j}(\tau_j) \, \Big\{\partial^{k} \prod_{i=1}^n D_{a_i} \Big\} \Upsilon_p\big(\CI_{b_1}(\sigma_1) \curvearrowright\sigma_2\big)   \\ 
 &= \prod_{j =1}^{n} \Upsilon_{l_j}(\tau_j) \, \Big\{\partial^{k} \prod_{i=1}^n D_{a_i} \Big\} \Upsilon_p(\sigma),
\end{split} \end{equation*}
 as we have from the induction assumption that $\sigma_2$ satisfies the conclusion of the proposition, so $\Upsilon_p\big(\CI_{b_1}(\sigma_1) \curvearrowright\sigma_2\big) = \Upsilon_b(\sigma_1) D_{b_1}\hspace{-0.1cm} \Upsilon_p(\sigma_2)$. This concludes the proof.
\end{proof}


 Lemma \ref{LemBailleulHoshino} was first noticed in the proof of Proposition 33 in Bailleul and Hoshino's work \cite{RSGuide}. It lead the authors to a simple proof of the fact that for $a=(\frak{t}_j,p_a)$ and all $p$
\begin{equation*}
 \Upsilon_p\big(\mathcal{I}_{a}( \sigma) \curvearrowright \tau\big)  = \Upsilon_j(\sigma) \,D_a \Upsilon_p(\tau);
\end{equation*}
a special case of \eqref{star2_morphism}. This was a huge simplification in comparison to the original proof given in Bruned, Chandra, Chevyrev \& Hairer's work \cite{BCCH18}, where the authors had to go through an extended space of rooted trees in Section $4$ therein. Identity \eqref{star2_morphism} was observed in the simpler context of rough differential equations in Lemma 3.4 of Bonnefoi, Chandra, Moinat \& Weber's work \cite{BCMW}. The $\star$ product happens to be the adjoint of the Butcher-Connes-Kreimer coproduct in that setting.

\medskip

\section{Preparation maps and their associated models}
\label{SectionLocalProduct}

Preparation maps were introduced by Bruned in \cite{BrunedRecursive} as the building block of a family of renormalization procedures encompassing the BHZ renormalization procedure from \cite{BHZ}. Their dual action on the equation was not investigated so far. We introduce in the present section a notion of strong preparation map characterized by a morphism property for the $\star$ product.

\medskip

\subsection{Strong preparation maps} For $\tau\in T$ denote by $\vert\tau\vert_\zeta$ the number of noise symbols that appear in $\tau$. Recall from \cite{BrunedRecursive} the following definition.

\begin{definition}
A \textbf{\textsf{preparation map}} is a linear map $R : T\rightarrow T$ such that $R(X^k\tau)=X^kR(\tau)$, for all $k\in\N^{d+1}, \tau\in T$, and $R\mcI_a = \mcI_a$ for all $a$, and 
\begin{equation} \label{Commutation_R}
\left( R\otimes \textrm{Id} \right) \Delta = \Delta R.
\end{equation}
and for each $ \tau \in T $ there exist finitely many $\tau_i \in T$ and constants $\lambda_i$ such that one has
\begin{equation} \label{analytical}
R(\tau) = \tau + \sum_i \lambda_i \tau_i, \quad\textrm{with}\quad \textrm{deg}(\tau_i) \geq \textrm{deg}(\tau) \quad\textrm{and}\quad |\tau_i|_\zeta < |\tau|_\zeta.
\end{equation}
\end{definition}

Preparation maps act non-trivially on trees with multiple edges at the root. This accounts for the fact that such trees represent products of analytical quantities, some of which need to be renormalized to be given sense. The `deformed product' provided by $R(\tau)$ for such trees $\tau$ makes precisely that. Preparation maps were named for that reason `{\it local product renormalization maps}' in Chandra, Moinat \& Weber's work \cite{CMW} on the $ \phi^{4}_{4-\delta}$ equation, as well as in Bruned's work \cite{Br1111} on the renormalization of branched rough paths. A preparation map is in particular a perturbation of the identity by elements that are not more singular ($\textrm{deg}(\tau_i) \geq \textrm{deg}(\tau)$) and defined with strictly less noises ($|\tau_i|_\zeta < |\tau|_\zeta$). Note that the linear map $R-\textrm{Id}$ is nilpotent as a consequence of the condition \eqref{analytical}. Recall that the recentering operators $\Gamma_{yx}$ in a regularity structure are given by some maps of the form $\Gamma_{yx}=(\textrm{Id}\otimes {\sf g}_{yx})\Delta$, for some character ${\sf g}_{yx}$ on $T^+$ and arbitrary points $x,y$ in the state space.  Identity \eqref{Commutation_R} encodes the fact that the recentering operators and the preparation map commute
$$
\Gamma_{zy}R = R\Gamma_{zy}
$$
for all $z,y$ since 
$$
\Gamma_{zy}(R\tau) = (\textrm{Id}\otimes {\sf g}_{zy})\Delta R(\tau) =  (\textrm{Id}\otimes {\sf g}_{zy})(R\otimes \textrm{Id})\Delta\tau = (R\otimes {\sf g}_{zy})\Delta\tau = R\Gamma_{zy}(\tau).
$$   
The next statement is a direct consequence of the duality relation \eqref{def_coaction} between the product $\star$ and the splitting map $\Delta$.


\begin{proposition} \label{prop_equivalence}
Identity \eqref{Commutation_R} is equivalent to having
\begin{equation} \label{Commutation_R*Translation}
R^{*} \left( \sigma \star \tau \right)  =  \sigma \star (R^*(\tau))
\end{equation}
for all $\sigma \in T^+ $ and $\tau \in T$.
\end{proposition}


\begin{proof}
We use the duality relation \eqref{def_coaction} between $\Delta$ and the restriction to $T\otimes T^+$ of the $\star$ product to write for $\mu,\nu\in T$ and $\sigma\in T^+$
$$
\langle  \mu \otimes \sigma , \Delta R \nu  \rangle = \langle \sigma \star \mu , R \nu  \rangle = \big\langle R^* (\sigma \star \mu),  \nu  \big\rangle.
$$
The result is thus a consequence of the identity
$$
\big\langle  \mu \otimes \sigma , \left( R \otimes \textrm{Id} \right) \Delta \nu \big\rangle = \big\langle  (R^*(\mu)) \otimes \sigma ,  \Delta\nu  \big\rangle =  \big\langle  \sigma \star (R^*(\mu)) ,  \nu  \big\rangle.
$$
\end{proof}


\begin{definition}
A \textbf{\textsf{strong preparation map}} is a preparation map satisfying the identity \eqref{Commutation_R*Translation} for all $\sigma\in T$ and $\tau\in T$ -- and not only for $\sigma\in T^+$.
\end{definition}


Taking some specific $\sigma$ in \eqref{Commutation_R*Translation} for strong preparation maps yields some special identities. For $ \sigma = \CI_{a}(\sigma_1)$ the identity \eqref{Commutation_R*Translation} reads
\begin{equation*} \label{condition_1}
R^{*}  \big( \mathcal{I}_{a}(\sigma_1) \curvearrowright \, \tau \big) =   \mathcal{I}_{a}(\sigma_1)\curvearrowright(R^*(\tau)),
\end{equation*}
that is
\begin{equation} \label{condition_11}
R^*\big(\sigma_1\curvearrowright_a \, \tau\big) = \sigma_1\curvearrowright_a( R^*(\tau))
\end{equation}
Another interesting case is when $ \sigma $ is equal to $X^k$. In that case one has for all $k\in\N^{d+1}$ and $\tau\in T$ 
\begin{equation} \label{condition_22}
 R^*(X^k \star \tau) = (X^k \star R^* \tau)  = \uparrow^k R^{*} \tau
\end{equation}
Note that the universality property of $T$ stated in Proposition \ref{freeness} implies that the identities \eqref{condition_11} and \eqref{condition_22} characterize the map $R^*$ once its values on the generators $X^k\zeta_l$ are given.

\medskip

Denote by $\mcB^-$ the canonical basis of $T^-$ and recall from \cite{BHZ} or \cite{RSGuide} that $\mathbb{R}[T^-]$ is equipped with a Hopf algebra structure. We follow \cite{BrunedRecursive} and define for any character $\ell$ of $\mathbb{R}[T^-]$ and all $\sigma \in T$ 
\begin{equation} \label{def_R}
R^{*}_\ell(\sigma) \defeq \sum_{ \tau \in \mathcal{B}^-} \frac{\ell( \tau)}{S( \tau)} \,(  \sigma\star  \tau).
\end{equation}
This definition corresponds to the dual of its usual definition 
$$
R_\ell(\mu) = \sum_{ \tau\in\mcB^-,  \sigma\in\mcB}\frac{\ell( \tau)}{S( \tau)}\,( \tau\otimes \sigma, \delta_r\mu) \sigma,
$$
 where the map $\delta_r$  can be defined form the duality relation $\langle\mu , R_\ell^*(\sigma)\rangle = \langle R(\mu) , \sigma\rangle$.


\begin{proposition} \label{PropStrongPrepMapsREll}
 Assume that $\ell$ is null on planted trees and trees of $T^-$ of the form $X^k\sigma$. The map $R_\ell$ is a strong preparation  map.
\end{proposition}


\begin{proof} 
From definition \eqref{def_R} one has for any $\mu,\tau\in T$ 
\begin{equation*}
R^*_\ell (\mu \star \tau) = \sum_{\sigma \in \mathcal{B}^-} \frac{\ell(\sigma)}{S(\sigma)} \,( \mu \star \tau) \star \sigma.
\end{equation*}
By using the associativity of $\star$ one gets 
\begin{equation*}
R^*_\ell( \mu \star \tau) = \sum_{\sigma \in \mathcal{B}^-} \frac{\ell(\sigma)}{S(\sigma)} \,  \mu \star (\tau\star \sigma) =  \mu \star R^{*}_\ell(\tau),
\end{equation*}
so Proposition \ref{prop_equivalence} ensures that the commutation condition \eqref{Commutation_R} is in particular satisfied. The condition $\textrm{deg}(\tau_i)\geq \textrm{deg}(\tau)$ comes from the fact that we are summing over decorated trees with negative degree in the definition of $R^*_\ell$. The condition on $|\cdot|_\zeta$ is a consequence of the the fact that the $\delta_r$ map breaks any decorated tree into a sum of tensor products of trees with a strictly smaller number of noises.
\end{proof}

 The BPHZ renormalization map from \cite{BHZ, CH} is a map built from $R_\ell$ for a particular choice of character $\ell$ on $T^-$. Although elementary the next statement will play a crucial role in the proof of our main result, Theorem \ref{ThmMain}, in the next section.


\begin{proposition} \label{PropCombined}
Let $R$ be a strong preparation map. For every $1\leq p\leq p_0$, for every $\tau =  X^{k} \zeta_{l}\prod_{j=1}^n \CI_{a_j}(\tau_j)$ with $a_j=( \frak{t}_+^{l_j} , p_j)\in\frak{T}^+\times\N^{d+1}$, one has
$$
 \Upsilon_p(R^*(\tau)) = \Big\{ \partial^k D_{a_1} \dots D_{a_n}  \Upsilon_p(R^*(\zeta_l))\Big\}\,\prod_{j=1}^{n}  \Upsilon_{l_j}(\tau_j)
$$
\end{proposition}


\begin{proof}
Writing $\tau = \left( X^{k} \prod_{j=1}^{n} \CI_{a_j}(\tau_j) \right) \star \zeta_{l}$ and using the property \eqref{Commutation_R*Translation} one gets
$$
R^{*} \left( \left( X^{k} \prod_{j=1}^n \CI_{a_j}(\tau_j) \right) \star \zeta_{l} \right) = \left(  X^{k} \prod_{j=1}^n \CI_{a_j}(\tau_j) \right) \star \left( R^*(\zeta_l) \right).
$$
Identity \eqref{star2_morphism} in Proposition \ref{star_morphism} then yields
$$
 \Upsilon_p\bigg(X^{k} \prod_{j=1}^n \CI_{a_j}(\tau_j)  \star \left( R^*(\zeta_l) \right) \bigg) = \Big\{ \partial^k D_{a_1} \dots D_{a_n} \big( \Upsilon_p(R^*(\zeta_l))\Big\}\,\prod_{j=1}^{n}  \Upsilon_{l_j}(\tau_j).
$$
\end{proof}

\subsection{Models associated with a preparation map.}
We recall now from \cite{BrunedRecursive} the construction of a continuous admissible model associated with a preparation map. This allows to emphasize in \eqref{EqFeatureOne} and \eqref{EqRelationPiMPiMCirc}  below a fundamental feature of its associated reconstruction map.

\ssk

 Assume we are given some kernels $(K^p)_{1\leq p\leq p_0}$ with a polynomial singularity at $0$ satisfying Assumption (5.1) in Hairer's seminal work \cite{Hai14}, and some continuous noises $(\xi_l)_{1\leq l\leq l_0}$ on the state space. Following \cite{BrunedRecursive} one can associate to a preparation map $R$ an admissible model ${\sf M}^{\!R}$ on $T$. It is defined from a side family $\big((\Pi_x^{\!R\times}\tau)(\cdot)\big)_{x,\tau}$ of continuous functions on the state space satisfying
$$
\big(\Pi_x^{\!R\times} \one\big)(y) = 1  , \quad \big(\Pi_x^{\!R\times} \zeta_l\big)(y)  = \xi_l(y) , \quad \big(\Pi_x^{\!R\times} X_i\big)(y) = y_i-x_i, 
$$
{\it the multiplicativity condition}
$$
\big(\Pi_x^{\!R\times}(\sigma\tau)\big)(y) = \big(\Pi_x^{\!R\times}\sigma\big)(y) \, \big(\Pi_x^{\!R\times}\tau\big)(y),
$$
and the condition
\begin{equation} \label{EqConditionAlmostAdmissibility}  \begin{aligned}
\Big(\Pi_x^{\!R\times} \CI_{a} (\tau)\Big)(y) & = \hspace{-0.05cm}\Big( D^q K^p * \Pi^{\!R\times}_{x}(R\tau) \Big)(y) \\ & - \hspace{-0.4cm}\sum_{|r|_{\s} \leq \textrm{deg}( \CI_{a}\tau)} \frac{(y-x)^{r}}{r!} \Big(D^{q + r} K^p * \Pi_x^{\!R\times}(R\tau)\Big)(x),
\end{aligned} \end{equation}
for  $a=(\frak{t}_+^p,q)$. Define for all $x$ and $\tau$ a  continuous function on the state space
$$
\big(\Pi^{\!R}_x\tau\big)(\cdot) \defeq \Big(\Pi^{\!R\times}_x(R\tau)\Big)(\cdot).
$$
Bruned gave in Proposition 3.16 of \cite{BrunedRecursive} an explicit contruction of an  $R$-dependent admissible model ${\sf M}^{\!R} = ({\sf \Pi}, {\sf g}) =({\sf \Pi}^R, {\sf g}^R)$ on $T$, with values in the space of  continuous functions, such that the operators $\Pi^{\!R}_x$ are indeed associated with  the model ${\sf M}^{\!R}$ in the sense that one has for all $\tau\in T$ and $x$
$$
\Pi^{\!R}_x\tau  = \big({\sf \Pi}^{\!R}\otimes({\sf g}^{\!R}_x)^{-1}\big)\Delta\tau.
$$
 (You will find some details in Section \ref{SubsectionStateSpaceDependentPreparationMaps} in a more general setting.) Since the model ${\sf M}^{\!R}$ takes values in the space of continuous functions its associated reconstruction operator ${\sf R}^{{\sf M}^{\!R}}$ is given by the explicit formula
$$
\big({\sf R}^{{\sf M}^{\!R}} {\sf v}\big)(x) = \big({\sf \Pi}_x^{\sf g}{\sf v}(x)\big)(x)
$$
for any modelled distribution $\sf v$ with positive regularity, so 
\begin{equation} \label{EqFeatureOne}
\big({\sf R}^{{\sf M}^{\!R}} {\sf v}\big)(x) = \Big(\Pi_x^{\!R\times}\big(R{\sf v}(x)\big)\Big)(x).
\end{equation}
Since $R\CI_a = \CI_a$ one has in particular
$$
\Pi^{\!R}_x(\CI_a\tau) = \Pi_x^{\!R\times}(\CI_a\tau).
$$
{\it The point here is that $\Pi_x^{\!R\times}$ is multiplicative while $\Pi_x^{\!R}$ is not multiplicative.} Denote by $T_X\subset T$ the linear space spanned by the polynomials in $T$. Modelled distribution $\sf u$ with values in the subspace $\CI(T)\oplus T_X$ of $T$ satisfy in that case the identity
\begin{equation} \label{EqRelationPiMPiMCirc}
\big({\sf R}^{{\sf M}^{\!R}} {\sf u}\big)(x) = \Big(\Pi_x^{\!R\times}{\sf u}(x)\Big)(x).
\end{equation}
It follows further from the relation \eqref{EqConditionAlmostAdmissibility} that the model $\sf M$ is admissible.


\section{A short proof for the renormalized equation for translation invariant operators}
\label{SectionShortProof}

 For $\gamma>0$ denote by $\mathcal{Q}_\gamma$ the natural projection from $T$ to the linear subspace $T_{<\gamma}$ of elements of $T$ of degree less than $\gamma$. Recall that for a smooth function $G$, with derivatives $G^{(k)}$, and a modelled function-like distribution ${\sf v}$ of regularity $\gamma$ the formula 
\begin{equation*} \label{ref_composition}
\widehat{G}({\sf v}) =  \mathcal{Q}_\gamma \sum_{k\in\N^{d+1}} \frac{G^{(k)}(\langle{\sf v},\one\rangle)}{k!}\,\big({\sf v} - \langle{\sf v},\one\rangle\one\big)^k
\end{equation*}
defines a modelled distribution $\widehat{G}({\sf v})$ of regularity $\gamma$. This section contains the statement and proof of  a particular case of our main result, Theorem \ref{ThmMain}  when the operators $\mcL^p$ are translation invariant. Theorem \ref{ThmMainTranslation} below describes the autonomous dynamics satisfied by ${\sf R}^{{\sf M}^{\!R}}({\sf u}^{\!R})$ when ${\sf M}^{\!R}$ is the model constructed from a preparation map $R$ and ${\sf u}^{\!R}$ is the solution to the lift of system \eqref{EqSPDE} to its associated regularity structure 
\begin{equation} \label{EqLiftedSystem}
{\sf u}^{\!R}_p = \mathcal{K}^{{\sf M}^{\!R}}_p\Big(\mathcal{Q}_{\gamma-2} \big(\widehat{F}_p({\sf u}^{\!R} , D{\sf u}^{\!R})\zeta\big)\Big) + \mathcal{P}_{p,\gamma} \big(u_i(0)\big),
\end{equation}
 where $1\leq p\leq p_0$. The real number $\gamma$ is here strictly bigger than $3/2$. The notation $\widehat{F}_p({\sf u}^{\!R} , D{\sf u}^{\!R})\zeta$ is a shorthand notation for $\widehat{f}_p({\sf u}^{\!R})\zeta + \widehat{g}_p({\sf u}^{\!R})(D{\sf u}^{\!R})^2$. The operator $\mathcal{P}_{p,\gamma}$ stands here for the projection on the subspace of elements of degree less than $\gamma$ of the natural lift in the polynomial regularity structure $T_X$ of the map $x=(x_0,x_1)\mapsto e^{-x_0\mcL^p}(u(0))(x_1)$. The operator $\mathcal{Q}_{\gamma-2}$ is the natural projection from $T$ to the linear subspace $T_{<\gamma-2}$ of elements of $T$ of degree less than $\gamma-2$. The ${\sf M}^{\!R}$-dependent maps $\mathcal{K}^{{\sf M}^{\!R}}_p$ are the regularity structure lifts of the operators $(\partial_{x_0} - \mcL^p)^{-1}$. From the definition of the operators $\mathcal{K}^{{\sf M}^{\!R}}_p$, for a solution ${\sf u}^{\!R} = \big({\sf u}_1^{\!R},\dots,{\sf u}_{p_0}^{\!R}\big)$ of the system \eqref{EqLiftedSystem} each map ${\sf u}^{\!R}_p$ takes values in $\CI_{(\frak{t}_+^p,0)}(T)\oplus T_X$. Write
$$
{\sf u}^{\!R}_p =: \sum u_{p,\sigma}\sigma
$$
for a sum over trees $\sigma$ in the canonical basis of $\CI_{(\frak{t}_+^p,0)}(T)\oplus T_X$ -- monomials are seen as trees with just one vertex. 

One can get a general expansion for the derivatives of $ {\sf u}^{\!R}_{p} $ as 
\begin{equation*}
	 {\sf u}^{\!R}_{a} =  \sum u_{p,\sigma} D_a \sigma = \sum u_{a,\sigma} \sigma
\end{equation*}
where  $ a = (\frak{t}^p_+,q) $.
 Then, the variable $ Z_{a} $  introduced in \eqref{variables} corresponds to
\begin{equation*} 
	 u_{a, X^k} = \langle {\sf u}^{\!R}_{a}, X^{k} \rangle.
\end{equation*}
 It follows from \cite[Lem. 3.21]{BCCH18} that
\begin{equation*}
	\label{coherence_upsilon}
\mathcal{Q}_{\gamma-2}\big(\widehat{F}_p({\sf u}^{\!R} , D{\sf u}^{\!R})\zeta\big) = \sum_{\textrm{deg}(\tau)<\gamma-2} \frac{ \Upsilon_p(\tau)( u , \partial_{x_1}  u)}{S(\tau)} \,\tau;
\end{equation*}
a fact called {\it coherence}. In fact, one has for $ \sigma \in \CI_{(\frak{t}_+^p,0)}(T)  $
\begin{equation*}
	u_{a,\sigma} = u_{p, \sigma} = \frac{ \Upsilon_p(\tau)( u , \partial_{x_1}  u)}{S(\tau)}
\end{equation*}
and 
\begin{equation*}
	u_{a, X^{k}} = u_{p, X^{k+q}}.
\end{equation*}
 Set 
$$
u\defeq{\sf R}^{{\sf M}^{\!R}} ({\sf u}^{\!R}) = (u_1,\dots,u_{p_0}).
$$

\begin{theorem} \label{ThmMainTranslation}
 Let $\xi=(\xi_1,\dots,\xi_{l_0})$ stand for a continuous function. Let $R$ be a strong preparation map with associated admissible model ${\sf M}^{\!R}$. Then $u$ is a solution of the renormalized $p_0$-dimensional system 
\begin{equation} \label{EqRenormalizedSystem}
(\partial_t -  \mcL^p) u_p =  f_p(u)\xi + g_p(u)(\partial_{x_1}u)^2 + \sum_{l=0}^{l_0}  \Upsilon_p\Big(\big(R^* - \textrm{\emph{Id}}\big) \zeta_l\Big)(u,\nabla u)\,\xi_l.
\end{equation}
\end{theorem}


\begin{proof}
As we are working with an admissible model we have
\begin{equation*} \begin{split}
(\partial_t- \mcL^p) u_p &= {\sf R}^{{\sf M}^{\!R}}\Big(\mathcal{Q}_{\gamma-2}\big(\widehat{F_p}({\sf u}^{\!R} , D{\sf u}^{\!R})\zeta\big)\Big)   \\
&= {\sf R}^{{\sf M}^{\!R}}\bigg(  \sum_{\textrm{deg}(\tau)<\gamma-2} \frac{ \Upsilon_p(\tau)(u , \partial_{x_1} u)}{S(\tau)} \,\tau \bigg) \eqdef {\sf R}^{{\sf M}^{\!R}} ({\sf v}_p)
\end{split} \end{equation*}
for a sum over the canonical basis of $T$, from the coherence condition \eqref{coherence_upsilon}. One has by construction
$$
{\sf R}^{{\sf M}^{\!R}}({\sf v}_p)(x) = \Big( \Pi_x^{\!R\times}\big(R{\sf v}_p(x)\big)\Big)(x)
$$
and 
\begin{equation*} \begin{split}
R{\sf v}_p &= \sum_{\textrm{deg}(\tau)<\gamma-2} \langle R{\sf v}_p,\tau\rangle \tau =  \sum_{\textrm{deg}(\tau)<\gamma-2} \langle {\sf v}_p\,,\,R^*\tau\rangle \tau = \sum_{\textrm{deg}(\tau)<\gamma-2} \frac{ \Upsilon_p(R^*\tau)(u , \partial_{x_1} u)}{S(\tau)}\,\tau.
\end{split} \end{equation*}
Writing each $\tau$ in its canonical form $\tau=X^k\zeta_l\prod_{j=1}^n\CI_{a_j}(\tau_j)$, with $a_j=(\frak{t}_+^{l_{j}},k_j)$, and using the morphism property
$$
 \Upsilon_p(R^*\tau)(u , \partial_{x_1} u) = \Big\{\partial^kD_{a_1}\cdots D_{a_n}  \Upsilon_p(R^*\zeta_l)(u , \partial_{x_1} u)\Big\} \prod_{j=1}^n  \Upsilon_{l_j}(\tau_j)(u , \partial_{x_1} u)
$$
of the $\Upsilon_p$ maps obtained in Proposition \ref{PropCombined} one gets
\begin{equation*} \begin{aligned}
R {\sf v}_p & = \sum_{ \sigma = X^k\zeta_l\prod_{j=1}^n\CI_{a_j}(\tau_j)}   \frac{k! \prod_{j=1}^n S(\tau_j)}{S\big(X^k\zeta_l\prod_{j=1}^n\CI_{a_j}(\tau_j)\big)} \, \frac{X^k}{k!} \prod_{j=1}^n  \frac{ \Upsilon_{l_j}(\tau_j)(u , \partial_{x_1} u)}{S(\tau_j)}\, \CI_{a_j}(\tau_j)   \\ 
&\quad \Big\{\partial^k \prod_{j=1}^n D_{a_j} \Upsilon_p(R^*\zeta_l)(u , \partial_{x_1} u) \zeta_l\Big\},
\end{aligned} \end{equation*}
that is $R{\sf v}_p$ is equal to
$$
\sum_{ \sigma = X^k\zeta_l\prod_{j=1}^m\CI_{a_j}(\tau_j)^{\beta_j}} \frac{X^k}{k!} \prod_{j=1}^m \frac{1}{\beta_j !} \left(  \frac{ \Upsilon_{l_j}(\tau_j)(u , \partial_{x_1} u)}{S(\tau_j)}\,\CI_{a_j}(\tau_j) \right)^{\beta_j} \Big\{\partial^k \prod_{j=1}^m (D_{a_{j}})^{\beta_{j}}  \Upsilon_p(R^*\zeta_l)(u , \partial_{x_1} u) \zeta_l,\Big\}
$$
for some integer $m\leq n$, grouping together the terms that coincide and using the recursive definition of the symmetry factor $S(\cdot)$. Recall now that if  $G$ stands for a function of some variables $Z_{a_1},\dots, Z_{a_r}$ that are a qfinite subset of the $Z_a$ given in \eqref{variables}. The Fa\`a di Bruno formula from Lemma A.1 in \cite{BCCH18} states that
$$
\frac{\partial^{k}G}{k!} =  \sum_{a_1,...,a_r} \, \sum_{k = \sum_{j=1}^r \beta_j k_j} \; \prod_{j=1}^r \frac{1}{\beta_j !} \left(\frac{Z_{a_j + k_j}}{k_j!}\right)^{\beta_j} \prod_{j=1}^r (D_{b_j})^{\beta_j} G.
$$
 where the $ k_j$ in the previous sum are non-zero. We can apply this identity to the functions $G = \Upsilon_p(R^*\zeta_l)$ of the variables $Z_{b_i} = \langle {\sf u}^{\!R}_{b_i}, \one \rangle$.  
One obtains
\begin{equation*} \begin{aligned}
&R{\sf v}_p   \\
&= \sum_{\sigma = X^k\zeta_l\prod_{j=1}^m\CI_{a_j}(\tau_j)^{\beta_j}}
\sum_{a_{m+1},...,a_{m+r}} \, \sum_{k = \sum_{j=m+1}^{m+r} \beta_j k_j} \frac{X^k}{k!} \prod_{j=1}^{m} \frac{1}{\beta_j !} \left(  \frac{\Upsilon_{l_j}(\tau_j)(u , \partial_{x_1} u)}{S(\tau_j)}\,\CI_{a_j}(\tau_j) \right)^{\beta_j}   \\ 
&\quad  \prod_{j=m+1}^{r} \frac{1}{\beta_j !} \left(  \frac{\langle {\sf u}^{\!R}_{a_j}, X^{k_j} \rangle}{k_j!}\, X^{k_j} \right)^{\beta_j} \Big\{ \prod_{j=1}^{m+m'} (D_{a_{j}})^{\beta_{j}}  \Upsilon_p(R^*\zeta_l)(u , \partial_{x_1} u) \zeta_l,\Big\}.
\end{aligned} \end{equation*}
We notice that for $a=(\frak{t}_+^{l},k)$ one has
\begin{equation*}
	{\sf u}^{\!R}_{a} - \langle {\sf u}^{\!R}_{a}, \one \rangle = \sum_{k \neq 0} \frac{\langle {\sf u}^{\!R}_{a}, X^{k} \rangle}{k!}\, X^{k} + \sum_{\tau} \frac{{\Upsilon}_{l}(\tau)(u , \partial_{x_1} u)}{S(\tau)}\,\CI_{a_j}(\tau) \tau
\end{equation*}
Using the binomial identity, one gets
\begin{equation*} \begin{split}
R{\sf v}_p &= \sum_{l,n} \underset{\beta_1,...,\beta_n}{\sum_{a_1,...,a_n}} \, \prod_{j=1}^n \frac{1}{\beta_j !} \left({\sf u}^{\!R}_{a_j} - \langle {\sf u}^{\!R}_{a_j}, \one \rangle \right)^{\beta_j} \prod_{j=1}^n (D_{a_{j}})^{\beta_{j}} \Upsilon_p (R^*\zeta_l)(u , \partial_{x_1} u) \zeta_l   \\
&= \sum_{l=0}^{l_0}  \widehat{\Upsilon_p(R^{*}\zeta_l)}({\sf u}^{\!R} , D{\sf u}^{\!R}) \zeta_l.
\end{split} \end{equation*}
 where $ \widehat{\Upsilon_p(R^{*}\zeta_l)} $ has to be understood in the sense of \eqref{ref_composition}. Using the multiplicativity of $\Pi_x^{M^{\!\times}_{\!R}}$ and identity \eqref{EqRelationPiMPiMCirc} identifying ${\sf R}^{{\sf M}^{\!R}}({\sf u}^{\!R})(x)$ we obtain
\begin{equation*} \begin{split}
\big((\partial_{x_0} - \mcL^p)u_p\big)(x) &= {\sf R}^{{\sf M}^{\!R}}({\sf v}_p)(x)  =  \Pi_x^{\!R\times}\big(R{\sf v}_p(x)\big)(x)   \\
&= \sum_{l=0}^{l_0} \Pi_x^{\!R\times}\Big( \widehat{\Upsilon_p\big(R^{*}\zeta_l)}({\sf u}^{\!R}(x), D{\sf u}^{\!R}(x)\big) \zeta_l\Big)(x)   \\
&= \sum_{l=0}^{l_0}  \Upsilon_p(R^*\zeta_l)\Big(\big( \Pi_x^{\!R\times}{\sf u}^{\!R}(x)\big)(x), \partial_{x_1} \big( \Pi_x^{\!R\times}{\sf u}^{\!R}(x)\big)(x)\Big) \,  \Pi_x^{\!R\times}\zeta_l   \\
&= \sum_{l=0}^{l_0}  \Upsilon_p(R^*\zeta_l)\big(u(x) , \partial_{x_1} u(x)\big)\xi_l.
\end{split}\end{equation*} 
 This identity is equivalent to \eqref{EqRenormalizedSystem}.
\end{proof}

 We invite the reader to note right now that the above proof works verbatim if $R$ were $x$-dependent and we knew that there is indeed a model ${\sf M}^{\!R}$ with reconstruction operator ${\sf R}^{{\sf M}^{\!R}}$ and some multiplicative maps $\Pi_x^{\!R\times}$ such that 
\begin{equation} \label{EqStructuralRelation}
{\sf R}^{{\sf M}^{\!R}}({\sf v})(x) = \Pi_x^{\!R\times}\big(R(x){\sf v}(x)\big)(x)
\end{equation}
for all modelled distribution $\sf v$. We will see in Section \ref{SectionRenomalisationSchemes} that this is indeed the case.   

In the particular case where the preparation map $R$ is of BPHZ form \eqref{def_R} a direct computation shows that the renormalized system \eqref{EqRenormalizedSystem} takes the form
$$
 (\partial_t - \mcL^p) u_p = f_p(u)\xi + g_p(u)(\partial_{x_1}u)^2 + \sum_{\tau \in \mathcal{B}^{-}\setminus \{\one \} } \ell(\tau) \, \frac{ \Upsilon_p(\tau)(u,\partial_{x_1}u)}{S(\tau)} 
$$
if we further assume that $R_{\ell}^*(\zeta_l) = \zeta_l$ for all $l\neq 0$ -- recall $\zeta_0=\textbf{\textsf{1}}$. This assumption accounts for the fact that we never need to subtract a multiple of one of the noises $(\zeta_l)_{1\leq l\leq l_0}$ to any tree-indexed quantity in our renormalization algorithm. We note that Chandra, Moinat \& Weber used in \cite{CMW} a similar strategy to get back the renormalized equation in the particular case of the $\Phi^4_{4-\delta}$ equation.

\medskip

\noindent \textbf{\textsf{Remarks --}} \textbf{\textsf{1.}} {\it The pattern of proof of the renormalized equation devised in \cite{BCCH18} has its roots in Bruned, Chevyrev, Friz \& Preiss' work \cite{BCFP} and rests heavily on the formalism and properties of multi-pre-Lie algebras. The renormalized equation is obtained from the fact that one can associate to a BHZ preparation map $R_\ell$ a renormalization map  $M_\ell$ whose dual has a pre-Lie morphism property. Writing $\Pi_x$ for the canonical interpretation map associated with a continuous noise, and $\Pi_x^\ell$ for the interpretation map of the renormalized model, the core of the mechanics of \cite{BCCH18} boils down to the identity
\begin{equation*} \begin{split}
\Pi^\ell_x({\sf v}_p) &= \sum_{\textrm{\emph{deg}}(\tau)<\gamma-2} \frac{ \Upsilon_p(\tau)(u,\partial_{x_1}u)}{S(\tau)}  \, (\Pi_x M_\ell \tau)   \\
			      &= \sum_{\textrm{\emph{deg}}(\tau)<\gamma-2} \frac{ \Upsilon_p^{M_{\ell}}(\tau)(u,\partial_{x_1}u)}{S(\tau)}  \, (\Pi_x  \tau) 
\end{split} \end{equation*}
in which the functions $ \Upsilon_p^{M_\ell}$ are constructed from the pre-Lie property. This pattern of proof cannot work in a non-translation invariant setting as one cannot build  some renormalized models for which $\Pi^{\!R}_x = \Pi_x\circ M_{\!R}(x)$.   }

\smallskip

\textbf{\textsf{2.}} {\it The proof of Theorem \ref{ThmMainTranslation} has been used in many places as it provides a shorter proof than the original proof given in \cite[Thm 2.2, Eq (2.21)]{BCCH18}. This proof is very robust when one considers a different framework. For example, this resulf has been used for quasilinear equations in \cite[Prop. 24]{BailleulHoshinoKusuoka} and \cite[Thm 3.4]{BGN24}, where a variant of the proof is given with non-commutative derivatives. In a geometric context  Hairer \& Singh provide a proof of the renormalised equation in \cite[Prop. 14.9]{HairerSingh} which follows the main steps of the present proof. Ideas of this proof were also used for composition and substitution for Regularity Structures B-series \cite{Br23} and Multi-indice B-series \cite{BEF24}.}

\section{Renormalization schemes and renormalized equation in a non-translation invariant setting}
\label{SectionRenomalisationSchemes}

\subsection{The analytic setting for non-translation invariant operators}
\label{SectionAnalytic}

Recall the operators
$$
\mcL^p v \defeq a_p(\cdot)\partial_{x_1}^2v
$$ 
have some smooth coefficients. The analysis of the regularity structure lift of the resolvent operator $(\partial_t - \Delta)^{-1}$ done by Hairer in Section 5 of \cite{Hai14} rests entirely on a notion of $\beta$-regularizing singular kernel whose definition takes profit of the translation-invariance of the heat operator -- see e.g. Assumption 5.1 in Section 5 of \cite{Hai14}. A slightly different approach is needed in a translation invariant setting to build the resolvent of the operator $\partial_t-\mcL^p$ and one can for instance follow the alternative heat kernel approach of Bailleul \& Hoshino's `{\it Tourist Guide}' \cite{RSGuide} which is fully developed in Hoshino's work \cite{HoshinoSemigroup}. 

\ssk

Define the non-positive symmetric fourth order elliptic differential operator on $\bbR\times\bbR$
$$
\mcG^p \defeq \partial_{x_0}^2 - (\mcL^p)^2,
$$
and denote by $K^p(t,x,x')$ the fundamental solution of the operator $\partial_t - \mcG^p$. We learn from Theorem 5.12 in \cite{HoshinoSemigroup} that all we need to build a setting where to prove a version of the multilevel Schauder estimates is that $K^p$ satisfies some (anisotrope) Gaussian type estimates. These estimates ensure in particular that one has the uniform estimate
\begin{equation} \label{EqAnalyticCondition}
\int_{\bbR^2} \big\vert \partial_x^n K^p(t,x,x')\big\vert\,d(x,x')^c\,dx \lesssim t^{\frac{c-\vert n\vert_{\frak{s}}}{4}} \qquad(\forall x'\in\bbR^2, \,n\in\bbN\times\bbN, \,c\in\bbR_+).
\end{equation}
This estimate gives in some technical computations (e.g. in Section 5 of \cite{HoshinoSemigroup}) the estimates provided by Hairer's decompositions of the $\beta$-regularizing singular kernels as a sum of kernels with support in some annuli. The fact that the heat kernel of $\mcG^p$ satisfies the above mentioned Gaussian type estimate is folklore in the PDE litterature. One can find all the details fully spelled out in Appendix A of Bailleul, Hoshino \& Kusuoka's work \cite{BailleulHoshinoKusuoka}.

\subsection{State space dependent preparation maps and their associated models}
\label{SubsectionStateSpaceDependentPreparationMaps}

 Preparation maps can be state space dependent. We define this notion in this section and show that one can associate to such preparation maps some admissible models when the noise $\xi$ is continuous. They satisfy the relation \eqref{EqStructuralRelation} so it follows from the proof of Theorem \ref{ThmMainTranslation} and the remark following it that the renormalized equation associated with system \eqref{EqSPDE} is indeed given by Theorem \ref{ThmMain}.

\begin{definition}
A (state space dependent) \textbf{\textsf{preparation map}} is a map 
$$
R : (\mathbb{R}\times\textbf{\textsf{T}}) \times T \rightarrow T
$$ 
such that all the maps $R(x,\cdot)$ fix the polynomials and $R(x,\cdot)\mcI_a = \mcI_a$ for all $a$, and such that 
\begin{itemize}
   \item[$\bullet$] for each $ \tau \in T $ there exist finitely many $\tau_i \in T$ and some continuous real-valued bounded functions $\lambda_i$ on $\mathbb{R}\times\textbf{\textsf{T}}$ such that for every $ x \in \mathbb{R}\times\textbf{\textsf{T}}$
\begin{equation*} \label{EqAnalytical}
R(x,\tau) = \tau + \sum_i \lambda_i(x) \tau_i, \quad\textrm{with}\quad \textrm{deg}(\tau_i) \geq \textrm{deg}(\tau) \quad\textrm{and}\quad |\tau_i|_{\zeta} < |\tau|_{\zeta},
\end{equation*} 

   \item[$\bullet$] one has for all $x\in\mathbb{R}\times\textbf{\textsf{T}}$
 \begin{equation*} \label{EqCommutationRDelta}
 \big( R(x,\cdot) \otimes \textrm{Id}\big) \Delta = \Delta R(x,\cdot).
 \end{equation*}
\end{itemize}
A \textbf{\textsf{strong preparation map}} is a preparation map satisfying 
\begin{equation*} \label{Commutation_R*}
R(x,\cdot)^{*} \left( \sigma \star \tau \right)  =  \sigma \star \left( R(x,\cdot)^{*} \tau \right)
\end{equation*}
 for all $\sigma\in T$ and $\tau\in T$ and all $x\in\bbR\times\textbf{\textsf{T}}$. 
\end{definition}

 As an example one sees from Proposition \ref{PropStrongPrepMapsREll} that for any $x$-dependent character $\ell(x,\cdot)$ of the algebra $T^-$ that is null on planted trees and trees of the form $X^k\sigma$ the map
\begin{equation*}
R_\ell(x) \defeq \big(\ell(x,\cdot)\otimes\textrm{Id}\big)\delta_r
\end{equation*}
is a strong preparation map. Recall that the co-action $\Delta : T \rightarrow T \otimes T^+$ satisfies in our setting the induction relation
\begin{equation*} \label{coaction} \begin{split} 
\Delta(\bullet) &\defeq \bullet\otimes \textsf{\textbf{1}}, \quad \textrm{ for }\bullet\in\big\{\textsf{\textbf{1}}, X_i, \zeta\big\},   \\
\Delta(\mcI_a\tau) &\defeq (\mcI_a\otimes\textrm{Id})\Delta + \sum_{\vert k +m\vert_{\frak{s}} < \textrm{deg}(\mcI_a\tau)} \frac{X^k}{k!}\otimes \frac{X^m}{m!}\mcI^+_{a+k+m}(\tau).
\end{split}\end{equation*}
 Denote by $S^+$ the antipode of the coproduct $\Delta^+: T^+\rightarrow T^+\otimes T^+$ and let $\xi$ stand for a continuous noise. Set 
$$
 {\sf \Pi}^{\!R}(\zeta) =  {\Pi}^{\!R\times}(\zeta) = \xi,
$$
and define inductively the maps ${\sf \Pi}^{\!R}$ and $\Pi^{\!R\times}$ by the relations
\begin{equation*} \label{defPiR}
{\sf \Pi}^{\!R}  = {\Pi}^{\!R\times} R, \quad {\Pi}^{\!R\times}(\tau \bar \tau) = ({\Pi}^{\!R\times}\tau) ({\Pi}^{\!R\times} \bar \tau), \quad  {\Pi}^{\!R\times}(\mathcal{I}_a\tau) = D^{a} K *  ({\sf \Pi}^{\!R} \tau).
\end{equation*}
It follows from this definition and the fact that $R\mcI_a=\mcI_a$ that the map $ {\sf \Pi}^{\!R}$ satisfies the admissibility condition ${\sf \Pi}^{\!R}(\mcI_a\tau) = D^{a} K *  ({\sf \Pi}^{\!R} \tau)$. Set
\begin{equation*} \begin{aligned}
{\sf \Pi}_x^{\!R} \tau = \left(  {\sf \Pi}^{\!R} \otimes  ({\sf g}_x^{\!R})^{-1} \right) \Delta
\end{aligned} \end{equation*}
with
\begin{equation*} \label{EqRelationAdmissibleModel}
({\sf g}_x^{\!R})^{-1}\big(\mathcal{I}^+_a\tau\big) \defeq - \big(D^{a} K * {\sf \Pi}^{\!R}_x\tau\big)(x)
\end{equation*}
and
\begin{equation*} \label{reexpansion_map}
\widehat{{\sf g}^{\!R}_{yx}} \defeq \left(  \textrm{Id} \otimes {\sf g}_{yx}^{\!R} \right) \Delta
\end{equation*}
and
$$
{\sf g}^{\!R}_{yx} \defeq \Big(  ({\sf g}_y^{\!R})^{-1} \circ S^{+} \otimes  ({\sf g}_x^{\!R})^{-1} \Big) \Delta^{\!+}.
$$
An induction on $\textrm{deg}(\tau)$ shows that ${\sf \Pi}_x^{\!R}$ is the solution of the induction relation
\begin{equation*} \begin{aligned}
\big({\sf \Pi}_x^{\!R} \tau\big)(y) &= \Big(\Pi^{\!R\times}_x \big(R(y)\tau\big)\Big)(y),   \\
\Pi^{\!R\times}_x (\tau \bar \tau) &= (\Pi^{\!R\times}_x \tau) \, (\Pi^{\!R\times}_x\bar \tau),  \\ 
\big(\Pi^{\!R\times}_x(\mathcal{I}_a\tau)\big)(y) &= \big(D^q K^p * {\sf \Pi}_x^{\!R} \tau\big)(y) - \sum_{|k|_{\frak{s}} \leq \textrm{deg}(\mathcal{I}_a\tau) }\frac{(y-x)^k}{k!}  \big(D^{q +k} K^p * {\sf \Pi}_x^{\!R} \tau\big)(x),
\end{aligned} \end{equation*}
for $a=(\frak{t}_+^p,q)$. This was done in Section 3.2 of \cite{BrunedRecursive} in the case where the preparation map $R$ does not depend on its $\bbR^{2}$ argument. The proof carries over verbatim to the present setting. It turns out to be useful to have the following recursive identity.

\begin{lemma} \label{LemGfromPi}
One has the identity
\begin{equation*} \label{EqRecursiveG}
\widehat{{\sf g}^{\!R}_{yx}}\big(\mcI_a(\sigma)\big) = \mcI_a\Big(\widehat{{\sf g}^{\!R}_{yx}}(\sigma)\Big) - \sum_{\vert k\vert_{\frak{s}} < \textrm{\emph{deg}}(\mcI_a(\sigma))} \frac{(X+x-y)^k}{k !}\,{\sf \Pi}_x^{\!R}\Big(\mcI_{a+k}\big(\widehat{{\sf g}^{\!R}_{yx}}(\sigma)\big)\Big)(y).
\end{equation*}
\end{lemma}

Its proof follows the same lines as the proof of Lemma 12 in the authors' previous work \cite{BailleulBruned2}. 
    
\begin{proposition} \label{PropRenormalizedModel}
The pair ${\sf M}^{\!R} = \big({\sf g}^{\!R},{\sf \Pi}^{\!R}\big)$ defines a continuous admissible model.
\end{proposition}

\begin{proof}
We proceed in two steps.

{\bf 1.} We prove that
$$
\big\vert \big({\sf \Pi}_x^{\!R}\tau\big)(y) \big\vert \lesssim \vert y-x\vert_{\frak{s}}^{\textrm{deg}(\tau)}
$$
by induction on $\textrm{deg}(\tau)-\vert\tau\vert_\zeta$. Since we subtract in the definition $\Pi_x^{\!R\times}(\mcI_a\sigma)$ the correct Taylor expansion we have
$$
\big\vert\Pi_x^{\!R\times}(\mcI_a\sigma)(y)\big\vert \lesssim \vert y-x\vert_{\frak{s}}^{\textrm{deg}(\mcI_a\sigma)}.
$$
For $\tau=X^k\zeta_l\prod_{i=1}^n \mcI_{a_i}(\sigma_i)$, with $R(y)(\tau) = \tau + \sum \lambda_j(y)\tau_j$, the continuity of $\zeta_l$ and the multiplicativity of $\Pi_x^{\!R\times}$, together with the inductive assumption, give
\begin{equation*} \begin{split}
\big\vert \big({\sf \Pi}_x^{\!R}\tau\big)(y) \big\vert &\lesssim \big\vert \big( \Pi_x^{\!R\times}\tau\big)(y) \big\vert  + \sum_j \big\vert \big(\Pi_x^{\!R\times}\tau_j\big)(y) \big\vert   \\
&\lesssim \vert y-x\vert_{\frak{s}}^{\vert k\vert + \sum_j \textrm{deg}(\mcI_{a_i}(\sigma_i))} + \sum_j \vert y-x\vert_{\frak{s}}^{\textrm{deg}(\tau_j)} \lesssim \vert y-x\vert_{\frak{s}}^{\textrm{deg}(\tau)}.
\end{split} \end{equation*}

\ssk

{\bf 2. } We use the identity of Lemma \ref{LemGfromPi} to get the bounds on ${\sf g}_{yx}(\tau)$ inductively on $\textrm{deg}(\tau)$, using the bounds on ${\sf \Pi}_x$ obtained in the first point. This type of argument was first implemented in the proof of Proposition 3.16 in \cite{BrunedRecursive}; we recall the details here. 

For $\beta$ one of the finitely many values of $\textrm{deg}(\tau)$ less than a fixed constant $\beta_1$ denote by $\nu_\beta$ the canonical projection of an element $\nu\in T$ on the elements of degree $\beta$. The bounds on ${\sf g}_{yx}$ are equivalent to the inequality 
$$
\big\vert \big\{\widehat{{\sf g}^{\!R}_{yx}}(\tau)\big\}_\beta \big\vert \lesssim \vert y-x\vert_{\frak{s}}^{\textrm{deg}(\tau)-\beta}
$$
for all $\beta<\textrm{deg}(\tau)\leq \beta_1$. By multiplicativity of $\widehat{{\sf g}^{\!R}_{yx}}$ it suffices to prove this estimate for $\tau$ of the form $\mcI_a(\sigma)$. If $\beta$ is not an integer
$$
\big\{\widehat{{\sf g}^{\!R}_{yx}}(\mcI_a(\sigma))\big\}_\beta = \big\{\mcI_a\big(\widehat{{\sf g}^{\!R}_{yx}}(\sigma)\big)\big\}_\beta
$$
and we get from the induction assumption $\widehat{{\sf g}^{\!R}_{yx}}(\sigma) = \sigma + \sum_j c_{yx}^j\sigma_j$, with $\vert c_{yx}^j\vert \lesssim \vert y-x\vert_{\frak{s}}^{\textrm{deg}(\sigma)-\textrm{deg}(\sigma_j)}$, so the conclusion follows in that case. If $\beta=n<\textrm{deg}(\tau)$ is an integer we use Lemma \ref{LemGfromPi}, step 1 and the induction assumption to see that
\begin{equation*} \begin{split}
\big\vert \big\{\widehat{{\sf g}^{\!R}_{yx}}(\tau)\big\}_n \big\vert &\lesssim \sum_{n\leq \vert r\vert < \textrm{deg}(\tau)} \vert y-x\vert_{\frak{s}}^{\vert r\vert-n} \big\vert {\sf \Pi}_x^{\!R}\big(\mcI_{a+r}(\widehat{\sf g}^{\!R}_{yx}(\sigma))\big)(y)\big\vert   \\
&\lesssim \sum_{n\leq \vert r\vert < \textrm{deg}(\tau)} \vert y-x\vert_{\frak{s}}^{\vert r\vert-n} \sum_\beta \big\vert \big\{\widehat{\sf g}^{\!R}_{yx}(\sigma)\big\}_\beta \big\vert \vert y-x\vert_{\frak{s}}^{2-\vert a\vert-\vert r\vert +\beta} \lesssim \vert y-x\vert_{\frak{s}}^{\vert\tau\vert-n}.
\end{split} \end{equation*}   
\end{proof}

Since the model ${\sf M}^{\!R}$ takes values in the space of continuous functions its associated reconstruction operator ${\sf R}^{{\sf M}^{\!R}}$ is given by
$$
\big({\sf R}^{{\sf M}^{\!R}} {\sf v}\big)(x) = \big({\sf \Pi}_x^{\!R}{\sf v}(x)\big)(x)
$$
for any modelled distribution $\sf v$ with positive regularity, so 
$$
\big({\sf R}^{{\sf M}^{\!R}} {\sf v}\big)(x) = \Big(\Pi^{\!R\times}_x\big(R(x){\sf v}(x)\big)\Big)(x).
$$
 As we saw in the proof of Theorem \ref{ThmMainTranslation} the multiplicative character of the map $\Pi^{\!R\times}_x(\cdot)(x)$ is the crucial feature of this factorization of the reconstruction map ${\sf R}^{{\sf M}^{\!R}}$  that leads to the proof of Theorem \ref{ThmMain}.



\subsection{BHZ renormalization}
\label{SubsectionBHZ}

In a translation-invariant setting one can recast the renormalization scheme introduced by Bruned, Hairer and Zambotti in \cite{BHZ} in the setting of preparation maps by associating to any character $\ell$ of the algebra $T^-$ that is null on planted trees and trees of the form $X^k\sigma$ a strong preparation map $R_\ell$ whose dual $R_\ell^*$ is given by 
$$
R_\ell^*(\tau) = \sum_{\sigma\in\mcB^-} \frac{\ell(\sigma)}{S(\sigma)}\,(\tau\star\sigma),
$$
where $\mcB^-$ stands for the canonical basis of $T^-$. If the random noise $\xi$ in \eqref{EqSPDE} takes values in the space of continuous functions and has a translation-invariant distribution the BHZ renormalization corresponds to taking
$$
\ell(\tau) =\bbE\big[{\sf \Pi}(\widetilde{\mathcal{A}}_-\tau)(0)\big] ,
$$
where $ \widetilde{\mathcal{A}}_- $ is a `twisted antipode' associated with a co-action $ \delta : T \rightarrow T^{-} \otimes T $. This character is the BPHZ character introduced in Section 6.3 of \cite{BHZ} -- see also Section 7 of \cite{RSGuide}. The renormalized system \eqref{EqRenormalizedSystem} then takes the form
$$
(\partial_{x_0} - \mcL^p) u_p = f_p(u)\,\xi + g(x,u,\partial_{x_1}u) + \sum_{\tau\in\mcB^-\backslash\{\one\}} \ell(\tau) \, \frac{\Upsilon_p(\tau)\big(u,\partial_{x_1} u\big)}{S(\tau)} \qquad (1\leq p\leq p_0).
$$
In the general non-translation invariant setting, and  for any choice of strong preparation map $R$ such that $R(x)^*\zeta_l=\zeta_l$ for all $l\neq 0$ and all $x\in\bbR\times\textbf{\textsf{T}}$ the renormalized system takes the form 
$$
(\partial_{x_0} - \mcL^p) u_p = f_p(u)\,\xi + g_p(x,u,\partial_{x_1}u) +\Upsilon_p\big(R^*(\one)-\one\big)\big(u,\partial_{x_1} u\big) \qquad (1\leq p\leq p_0).
$$
In particular for any choice of state space dependent character $\ell$ on $T^-$ that is null on planted trees and trees of the form $X^k\sigma$, one has indeed $R_\ell(x)^*\zeta_l = \zeta_l$ for all $l\neq 0$ and all $x\in\bbR\times\textbf{\textsf{T}}$, and the renormalized system takes the form 
$$
(\partial_{x_0} - \mcL^p) u_p = f_p(u)\,\xi + g_p(x,u,\partial_{x_1}u) + \sum_{\tau\in\mcB^-\backslash\{\one\}} \ell(\cdot,\tau) \, \frac{\Upsilon_p(\tau)\big(u,\partial_{x_1} u\big)}{S(\tau)} \qquad (1\leq p\leq p_0).
$$
 As a final remark note that one has the following result in the translation invariant setting.
	
\begin{proposition}  For every $ \tau \in \mathcal{B}_- \backslash\{\one\} $
\begin{equation} \label{BPHZ}
\ell(\tau) = - \bbE\big[({\widetilde{\sf \Pi}}^{\!R_{\ell}}\tau)(0)\big] 
\end{equation}
where 
\begin{equation*}
	({\widetilde{\sf \Pi}}^{\!R_{\ell}}\tau)(0) = 
	({\sf \Pi}^{\!R_{\ell}}  \tau )(0)  - \ell(\tau)
\end{equation*}
with ${\sf \Pi}^{\!R_{\ell}}$ the renormalized smooth model associated with $ R_{\ell} $.
\end{proposition}
\begin{proof}
	We first recall that one has from \cite[Prop. 3.12]{BrunedRecursive} and \cite[Cor. 4.5]{BrunedRecursive}
	\begin{equation*}
		({\sf \Pi}^{\!R_{\ell}}\tau)(0) =   ({\sf \Pi} M_{\ell}\tau)(0).
	\end{equation*}
Then
\begin{equation*}
 \mathbb{E} \big[ ({\sf \Pi} M_{\ell}\tau)(0) \big] = 0
\end{equation*}
due to the choice of the character $ \ell $ (see \cite[Thm. 6.18, Eq. (6.26)]{BHZ}) which gives the result.	
\end{proof}

 This suggest to choose in a non-translation invariant setting a character defined inductively for all $x$ and all $\tau$ by the identity
\begin{equation}\label{suitable_C}
\ell(x,\tau) := - \bbE\big[(\widetilde{\sf \Pi}^{\!R_{\ell}}\tau)(x)\big].
\end{equation}  
Such a definition makes sense as one can rewrite $ (\widetilde{\sf \Pi}^{\!R_{\ell}}\tau)(x) $ as 
\begin{equation*}
	\begin{aligned}
	(\widetilde{\sf \Pi}^{\!R_{\ell}}\tau)(x)  &= ({\sf \Pi}^{\!R_{\ell}}\tau)(x) - \ell(x,\tau) 
	\\ & = \Big(\Pi^{\!R_{\ell}\times}_x\big(R_{\ell(x, \cdot)} \tau \big)\Big)(x)  - \ell(x,\tau)
	\\ & = \Big(\Pi^{\!R_{\ell}\times}_x\big(\widetilde{R}_{\ell(x, \cdot)} \tau \big)\Big)(x) 
	\end{aligned}
	\end{equation*}
where one has
\begin{equation*}
	\widetilde{R}_{\ell(x, \cdot)} \tau = \big(\ell(x, \cdot) \otimes \id \big) (\delta_r \tau - \tau \otimes \one).
\end{equation*} 
It is an immediate consequence from the definition \eqref{suitable_C} that one has
\begin{equation*}
	\bbE\big[({\sf \Pi}^{\!R_{\ell}}\tau)(x)\big]  = 0.
	\end{equation*}
Indeed, one has
\begin{equation*}
	\bbE\big[({\sf \Pi}^{\!R_{\ell}}\tau)(x)\big] =  \bbE\big[(\widetilde{\sf \Pi}^{\!R_{\ell}}\tau)(x)\big] + \ell(x,\tau)  =  - \ell(x,\tau) + \ell(x,\tau) = 0.
\end{equation*}


\medskip

\noindent \textbf{\textsf{Remarks --}} \textbf{\textsf{1.}} {\it How far are we from having a framework for dealing with systems of real-valued functions on a $2$ or $3$-dimensional closed manifold satisfying singular stochastic PDEs of the form \eqref{EqSPDE} with the differential of $u$ in the role of $\partial_{x_1}u$? The proof of the estimate \eqref{EqAnalyticCondition} is robust enough to be adapted to a closed manifold setting. The algebraic results of Section \ref{SectionRenomalisationSchemes} are insensitive to the fact that we could be on a manifold, and we know from Dahlqvist, Diehl \& Driver's work \cite{DDD} what notion of model to take. The only point where analysis enters the scene is when using the reconstruction theorem. The proof of the latter given in the `{\it Tourist guide}' \cite{RSGuide} works verbatim on a manifold setting once one has the estimates \eqref{EqAnalyticCondition}. Alternatively one can use the manifold version of Caravenna \& Zambotti's approach to the reconstruction theorem \cite{CaravennaZambotti} given by Rinaldi and Sclavi in \cite{RinaldiSclavi}. However the present setting is not sufficient to deal with more general systems of singular stochastic PDEs like the scalar $\Phi^4_3$ equation on a closed manifold. In that setting some more geometrical background needs to be introduced (jets in particular) to deal with local expansions of order higher than $1$. See the work \cite{HairerSingh} of Hairer \& Singh on the subject.}

\ssk

\textbf{\textsf{2.}} {\it The group $G^-$ of characters of the algebra $T^-$ comes equipped with a convolution product derived from the extraction/contraction coproduct $ \Delta^{\!-} : T^{-}  \rightarrow T^{-}  \otimes  T^{-} $:
\begin{align*} 
\ell \circ \bar \ell \defeq \left(\ell \otimes \bar{\ell} \,\right) \Delta^{\!-}, \quad \ell^{-1}(\cdot) = \ell( S^{-} \cdot),
\end{align*}
where $S^-$ is the antipode associated to $ \Delta^{\!-} $. In the translation invariant setting of all previous works, this group has an explicit representation in the space of endomorphisms on $ T $ given by
$$
M_{\ell} \defeq \left( \ell \otimes \textrm{\emph{Id}} \right) \delta,
$$
for which $\delta : T \rightarrow T^- \otimes T $ is a co-action and
$$
M_{\ell} \circ M_{\bar{\ell}} = M_{\ell \circ \bar{\ell}}.
$$
Although we cannot build such a representation of the group $G^-$ in a non-translation invariant setting we can associate 
$$
R_\ell(x) = \big(\ell(x,\cdot)\otimes\textrm{\emph{Id}}\big)\delta_r
$$ 
to any state space dependent character $\ell$ and any state space point $x$. For any two characters $\ell, \bar\ell$ on $\mcT^-$ set
\begin{align*}
\ell \circ \bar\ell \defeq \left( \ell \otimes \bar\ell \right) \delta_r,
\end{align*}
pointwise in $x$. One then has
$$
R_{\ell} \circ R_{\ell'} = R_{\ell \circ \bar\ell}.
$$   }

\section{Renormalization maps for branched rough paths}
\label{SectionRenormalizationRP}

The setting of subcritical singular stochastic PDEs devised in \cite{BHZ} defines {\it renormalization maps} over a given regularity structure as pairs of maps $(M,M^+)$ that turn any admissible models $\sf (\Pi, g)$ into an admissible model $({\sf \Pi}\circ M, {\sf g}\circ M^+)$. Branched rough paths provide an example of regularity structure with no polynomial decorations and for which $T=T^+$ and $M=M^+$ -- we refer the reader to Gubinelli's original article \cite{Gub10}, Hairer and Kelly's work \cite{HK} or \cite{BaiSimple,CassWeidner} for some background on this class of rough paths. This definition of a renormalization map takes it roots in the work \cite{BCFP} of Bruned, Chevyrev, Friz \& Preiss where the role of pre-Lie structures for rough differential equations was first noticed.

We give in this section an alternative picture of renormalization maps on branched rough paths and prove in Theorem \ref{adjoint_M} that they are all obtained from a preparation map. It is not clear yet whether this statement holds for an arbitrary regularity structure. However as part of the arguments below hold in full generality part of this section is written for a general regularity structure. Things are only specialized to branched rough paths in Corollary \ref{R_translation} where the generators are $ \zeta_{l} $ (there is no monomial decoration and $ \uparrow^k $ is not used). Moreover, in this case, we assume that $ \frak{T}^+ $ is a singleton. One has only one grafting product. Indeed, convolution with several kernels is replaced by integration in time.


\begin{definition} \label{DefnGoodMultiPreLieMorphism}
A \textbf{\textsf{good multi-pre-Lie morphism on $T$}} is a map $A : T\rightarrow T$ such that one has for all $\sigma, \tau\in T$ and $k\in\N^{d+1}$, and all $a\in\frak{T}^+\times\N^{d+1}$,
\begin{equation} \label{EqPreLieMorphism} \begin{split}
A(\sigma\curvearrowright_a\tau) &= (A\sigma)\curvearrowright_a(A\tau),   \\
A\hspace{-0.05cm}\uparrow^k\,&=\,\uparrow^k\hspace{-0.1cm}A.
\end{split} \end{equation}
\end{definition}
In the particular case of branched rough paths there are no polynomial decoration and good multi-pre-Lie morphisms are just multi-pre-Lie morphisms, that is they satisfy only the first identity in \eqref{EqPreLieMorphism}. We associate to a strong preparation map $R$ a linear map 
$$
M^{\!\times}_{\!R} : T\rightarrow T
$$ 
defined by the requirement that $M^{\!\times}_{\!R}(\one)  = \one$, that $M^{\!\times}_{\!R}$ is multiplicative, by the data of the $M^{\!\times}_{\!R}(\zeta_l)$, for $1\leq l\leq l_0$, and by the induction relation
\begin{equation} \label{EqConstructionRecipeBis}
M^{\!\times}_{\!R}\big(\CI_{(\Labhom,k)}\tau\big) = \CI_{(\Labhom,k)}\big(M^{\!\times}_{\!R}(R\tau)\big)
\end{equation}
for all $\tau\in T$ and $(\frak{t},k)\in\mathcal{L}^+\times\N^{d+1}$. Following Section 3.1 of \cite{BrunedRecursive} define
\begin{equation} \label{EqConstructionRecipe}
M_{\!R} \defeq M^{\!\times}_{\!R} R.
\end{equation}


\begin{theorem} \label{adjoint_M} 
\begin{itemize}
	\item[$\bullet$] Let $R$ be a strong preparation map. The map $M_{\!R}^*$ is a good multi-pre-Lie morphism satisfying $M_{\!R}^*(\zeta_l) = R^*(\zeta_l)$, for all $1\leq l\leq l_0$.
	\item[$\bullet$] Let $ M $ be a good multi-pre-Lie morphism then it is  obtained from a preparation map $R$ using the construction \eqref{EqConstructionRecipeBis}-\eqref{EqConstructionRecipe}. The preparation map $ R $ is defined  for $ \tau = \zeta_l \bar{\tau} $ with $ \bar{\tau} = \prod_{i=1}^n \CI_{a_i}(\tau_i) $ by
	\begin{equation} \label{def_preparation_map}
		R^*(\zeta_l) \defeq M^*(\zeta_l), \quad R^*(\tau)  = R^{*} (\bar{\tau} \star \zeta_l )  :=   \bar{\tau} \star R^{*} \zeta_l.
	\end{equation}   
\end{itemize}
\end{theorem}


\begin{proof}
 We start with $R$ a strong preparation map and we proceed in two steps.

-- We show that the map $(M^{\!\times}_{\!R})^*$ is also multiplicative. One has
\begin{equation*} \begin{split}
\big\langle (M^{\!\times}_{\!R})^{*} \CI_a(\sigma),  \CI_{a}(\tau)\big \rangle &=  \big\langle  \CI_a(\sigma), M^{\!\times}_{\!R} \CI_{a}(\tau) \big\rangle = \big\langle  \CI_a(\sigma),  \CI_{a}(M_{\!R} \tau) \big\rangle   \\ 
&= \langle  \sigma,  M_{\!R} \tau  \rangle = \langle  M_{\!R}^* \sigma,   \tau \rangle = \big \langle  \CI_a(M_{\!R}^*\sigma),  \CI_{a}(\tau) \big\rangle,
\end{split} \end{equation*}
which implies $(M^{\!\times}_{\!R})^{*} \CI_a(\sigma)  =  \CI_a(M_{\!R}^{*}\sigma)$. Denote by $\Delta_d$ the deconcatenation coproduct
\begin{equation*}
\Delta_d \CI_{a}(\tau) = \CI_{a}(\tau) \otimes \one + \one \otimes \CI_{a}(\tau), \quad \Delta_d X^k = X^k \otimes \one + \one\otimes X^k,
\end{equation*}
extended multiplicatively not to the tree product but to the product between a decorated tree and a decorated tree with no polynomial decorations at the root. It follows then from the multiplicativity of $M^{\!\times}_{\!R}$ and the identity $M^{\!\times}_{\!R} \CI_a(\sigma)  =  \CI_a(M_{\!R}\sigma)$, that
\begin{equation*}
\big\langle (M^{\!\times}_{\!R})^{*} \CI_{a}(\sigma) \mu,  \tau \big\rangle = \big\langle  \CI_{a}(\sigma) \mu, M^{\!\times}_{\!R} \tau \big\rangle = \big\langle  \CI_{a}(\sigma) \otimes \mu, \Delta_d M^{\!\times}_{\!R} \tau \big\rangle,
\end{equation*}
 where we have used the fact that $ \Delta_d $ is the adjoint of the joint root product $ \mathcal{M} $. We have $ \CI_{a}(\sigma) \mu = \mathcal{M}(\CI_{a}(\sigma) \otimes \mu) $.
One has also
\begin{equation*}
\Delta_d M^{\!\times}_{\!R} = \left( M^{\!\times}_{\!R} \otimes M^{\!\times}_{\!R} \right)\Delta_d.
\end{equation*}
Then, we get
\begin{equation*} \begin{split}
\big\langle  \CI_{a}(\sigma) \otimes \mu, \Delta_d M^{\!\times}_{\!R} \tau \big\rangle &= \big\langle  \CI_{a}(\sigma) \otimes \mu, \left( M^{\!\times}_{\!R} \otimes  M^{\!\times}_{\!R} \right)\Delta_d \tau \big\rangle   \\ 
&= \big\langle  (M^{\!\times}_{\!R})^{*} \CI_{a}(\sigma) \otimes (M^{\!\times}_{\!R})^{*} \mu, \Delta_d \tau \big\rangle = \big\langle  (M^{\!\times}_{\!R})^{*} \CI_{a}(\sigma)  (M^{\!\times}_{\!R})^{*} \mu,  \tau \big\rangle
\end{split} \end{equation*}
which concludes the proof that $(M^{\!\times}_{\!R})^{*}$ is multiplicative.

\smallskip

-- It will be useful for our purpose to decompose the grafting map $\curvearrowright_a$ into the sum of a grafting map at the root and a grafting map outside the root
$$
\curvearrowright_a = \curvearrowright_a^\textrm{root} + \curvearrowright_a^\textrm{non-root},
$$
with, for $\tau=X^k\prod_{j=1}^n\CI_{a_i}(\tau_j)$ and $a_i\in\frak{T}\times\N^{d+1}$,
$$
\sigma\curvearrowright_a^\textrm{root}\tau \defeq \sum_{m\in\N^{d+1}} {k\choose m}\,X^{k-m}\,\CI_{a-m}(\sigma)\prod_{j=1}^n\CI_{a_j}(\tau_j)
$$
and
$$
\sigma\curvearrowright_a^\textrm{non-root}\tau \defeq X^k \sum_{i=1}^n \CI_{a_i}(\sigma\curvearrowright_a\tau_i)\prod_{j\neq i}\CI_{a_j}(\tau_j).
$$
We proceed by induction on the size of the trees appearing in the product. In the induction hypothesis, we include the two following identities for $ \sigma, \tau \in \CT $
 \begin{equation} \label{ident_1}
 M_{\!R}^{*} \left( \sigma \curvearrowright_a \tau \right) =   (M_{\!R}^{*} \sigma) \curvearrowright_a  (M_{\!R}^{*}  \tau)
 \end{equation}
 and
 \begin{equation} \label{ident_2}
 (M^{\!\times}_{\!R})^{*} \left( \sigma \curvearrowright_a \tau \right) =   (M_{\!R}^{*} \sigma) \curvearrowright_a   \big((M^{\!\times}_{\!R})^{*} \tau\big).
 \end{equation}
  The base cases are given when one applies $  (M^{\!\times}_{\!R})^{*} $ and $  (M_{\!R})^{*}$ to the generators  $ X^k \zeta_{l} $ . One has
\begin{equation*}
	 (M^{\!\times}_{\!R})^{*}  X^k \zeta_{l} =  X^k \zeta_{l}, \quad (M^{\!\times}_{\!R})^{*} X^k \zeta_{l} =  R^{*}  (M^{\!\times}_{\!R})^{*}  X^k \zeta_{l} = R^{*} X^k \zeta_{l} =  \uparrow^k R^{*} \zeta_{l} 
\end{equation*}
where for the last identity we have used that
\begin{equation*}
R^{*}	X^k \zeta_{l} =  R^{*} (	X^k \star \zeta_{l}) =  X^k \star R\zeta_{l} = \uparrow^k R\zeta_{l}
	\end{equation*}
that is a consequence of the right morphism property of the strong preparation map with the product $ \star $ (see \eqref{Commutation_R*Translation}).
Let $ \sigma, \tau \in \CT $, one has
\begin{equation*} \begin{split}
M_{\!R}^{*} \left( \sigma \curvearrowright_a \tau \right) & =  R^* (M^{\!\times}_{\!R})^{*}\left( \sigma\curvearrowright_a \tau \right)   \\ 
&=  R^* \left( M_{\!R}^{*} \sigma \curvearrowright_a (M^{\!\times}_{\!R})^{*} \tau \right)   \\ 
&=  \left( M_{\!R}^* \sigma  \curvearrowright_a R^{*}(M^{\!\times}_{\!R})^{*} \tau \right) = \left( M_{\!R}^{*} \sigma\curvearrowright_a M_{\!R}^{*} \tau \right)
\end{split} \end{equation*}
where we have applied the induction hypothesis given by \eqref{ident_2} on $  (M^{\!\times}_{\!R})^{*} $ and the righ-morphism property of $ R^{*} $, given by the assumption that $R$ is a {\it strong} preparation map. We consider $ \tau = X^k \bar \tau $ where $ \bar \tau = \prod_{i=1}^n \CI_{a_i}(\tau_i) $. The multiplicativity property of $(M^{\!\times}_{\!R})^*$ and the fact that $(M^{\!\times}_{\!R})^*\CI_a(\sigma)  =  \CI_a(M_{\!R}^*\sigma)$ yield
\begin{equation*} \begin{split}
(M^{\!\times}_{\!R})^{*}\left( \sigma\curvearrowright_a^{\root}\tau \right) & = (M^{\!\times}_{\!R})^{*}  \sum_{\ell \in \N^{d+1}} {k \choose \ell} X^{k-\ell} \CI_{a-\ell}(\sigma) \bar \tau 
\\ &  =    \sum_{\ell \in \N^{d+1}} {k \choose \ell} X^{k-\ell} \CI_{a-\ell}(M_{\!R}^{*}\sigma) (M^{\!\times}_{\!R})^{*} \bar \tau 
 \\ & = M_{\!R}^{*} \sigma\curvearrowright_a^{\root}(M^{\!\times}_{\!R})^{*} \tau 
\end{split} \end{equation*}
For the grafting outside the root, we use the induction hypothesis. One has

\begin{equation*} \begin{split}
(M^{\!\times}_{\!R})^{*}\left( \sigma\curvearrowright_a^{\nonroot}\tau \right)  &= (M^{\!\times}_{\!R})^{*}  X^k  \CI_{a_i}\left( \sigma\curvearrowright_a  \tau_i \right) \prod_{j \neq i} \CI_{a_j}(\tau_j) 
\\&  = (M^{\!\times}_{\!R})^{*}( X^k)   \CI_{a_i}\left( M_{\!R}^{*} (\sigma\curvearrowright_a  \tau_i) \right) \prod_{j \neq i} \CI_{a_j}(M_{\!R}^{*} \tau_j) 
\\ & =   X^k \CI_{a_i}\left( M^{*}_R \sigma\curvearrowright_a  M_{\!R}^{*} \tau_i \right) \prod_{j \neq i} \CI_{a_j}(M_{\!R}^{*} \tau_j)
\\ & = M_{\!R}^* \sigma \curvearrowright_a^{\nonroot}(M^{\!\times}_{\!R})^{*}\tau.
\end{split} \end{equation*}
where we have used the induction hypothesis \eqref{ident_1} on $ \sigma $ and $ \tau_i $. The proof that $M_R^*\uparrow^k = \uparrow^{k} M_R^*$, works the same by decomposing insertion of polynomial decorations at the root and outside the root. For the second part of the theorem, let $ M $ be a good multi-pre-Lie morphism. The map $ R$ defined uniquely by \eqref{def_preparation_map} is actually a strong preparation map. Then, from the first part $M_{R}$ built out of $R$ is a multi-pre-Lie morphism which coincide with $M$ a multi-pre-Lie morphism on the generator. Therefore, one has $ M_R = M $ which concludes the proof.
\end{proof}

We now specialize the previous result to the setting of branched rough paths. Recall $T=T^+$ is in that case the Butcher-Connes-Kreimer Hopf algebra and all preparation maps are strong in that setting. We can then apply Theorem \ref{adjoint_M} to obtain the following fact.

\begin{corollary} \label{R_translation}
In the Butcher-Connes-Kreimer setting the family of (automatically good) multi-pre-Lie morphisms are in bijection with the set of preparation maps.
\end{corollary}


 In this particular setting Theorem \ref{ThmMain} takes the following form. Let ${\sf X}$ stand for a branched rough path over $\bbR^e$, where $e\geq 2$. In a regularity structure terminology $\sf X$ provides a model over the regularity structure $(T^+,T^+)$ with ${\sf g}_{ts}={\sf \Pi}_s(\cdot)(t)={\sf X}_{ts}$. Denote by $\zeta_0, \zeta_1,\dots,\zeta_e$ the generators of $T^+$. If we are given some regular enough vector field $V_0,\dots, V_e$ on $\bbR^d$ the map 
$$
V : \zeta_i \mapsto V_i
$$
has a unique extension as a pre-Lie morphism from $(T^+,\curvearrowright)$ to the the set of vector fields on $\bbR^d$ with the pre-Lie bracket $(W_1,W_2) \mapsto (DW_2)(W_1)$. Given a preparation map $R$ on $T^+$ the path $(z_t)$ is a solution of the differential equation
$$
dz_t = V(z_t) d{\sf X}^{\!R}_t
$$
driven by the renormalized branched rough path ${\sf X}^{\!R} = {\sf X}\circ M_{\!R}$ if and only if
$$
dz_t = V(z_t)d{\sf X}_t + V\big((M_{\!R}^*-\textrm{Id})\one\big)(z_t) dt.
$$

\bigskip

\noindent {\bf Funding:} I.B. acknowledges support from the CNRS and PIMS and the ANR-16-CE40-0020-01 grant. 	Y.B. gratefully acknowledges funding support from the European Research Council (ERC) through the ERC Starting Grant Low Regularity Dynamics via Decorated Trees (LoRDeT), grant agreement No.\ 101075208.




\end{document}